\input amstex
\documentstyle{amsppt}
\NoBlackBoxes
\magnification 1100
\define\cal{\Cal}
\define\br{\Bbb R}
\define\bz{\Bbb Z}
\define\ho{\operatorname{holink}}
\define\map{\operatorname{Map}}
\define\mapst{\operatorname{{Map}_s}}
\define\inclusion{\operatorname{inclusion}}
\define\proj{\operatorname{proj}}
\define\diam{\operatorname{diam}}
\define\lub{\operatorname{lub}}
\define\inr{\operatorname{int}}
\define\im{\operatorname{Im}}
\define\cyl{\operatorname{cyl}}
\define\ocyl{\overset\circ\to\cyl}
\define\M{\operatorname{M}}
\define\oc{\overset\circ\to{\text{\rm c}}}
\define\id{\operatorname{id}}
\define\cl{\operatorname{cl}}
\define\tp{\operatorname{TOP}}
\define\btop{\operatorname{BTOP}}
\define\level{\operatorname{level}}
\define\rel{\text{\rm rel\ }}
\define\SN{\operatorname{SN}}
\define\st{\operatorname{st}}
\define\maf{\operatorname{MAF}}
\define\ov{\overline}
\define\topbfri{\tp^b(F\times\br^i)}
\define\cbfri{\text{\rm C}^b(F\times\br^i)}
\define\topbf{\tp^b(F}
\define\cbf{\text{\rm C}^b(F}
\define\hcobf{h\text{\rm cob}(F)}
\define\ihcobf{\text{\rm I}h\text{cob}(F)}
\define\ihcobfpr{\text{\rm I}h\text{cob}(F')}
\define\Wh{\operatorname{Wh}}

\topmatter
\title Neighborhoods in Stratified Spaces \\ with Two Strata\endtitle
\rightheadtext{STRATIFIED SPACES}
\author  Bruce Hughes,
Laurence R. Taylor, Shmuel Weinberger \\ and Bruce Williams\endauthor
\leftheadtext{B. Hughes,
L. Taylor, S. Weinberger and B. Williams}
\address Department of Mathematics, Vanderbilt University, Nashville,
TN 37240 \endaddress 
\email hughes\@math.vanderbilt.edu \endemail
\address Department of Mathematics, University of Notre Dame, Notre Dame,
IN 46556\endaddress 
\email taylor.2\@nd.edu \endemail
\address Department of Mathematics, University of Chicago,
Chicago, IL 60637\endaddress 
\email shmuel\@math.uchicago.edu \endemail
\address Department of Mathematics, University of Notre Dame, Notre Dame,
IN 46556\endaddress 
\email williams.4\@nd.edu \endemail

\thanks The first author was supported in part by NSF Grant DMS--9504759
\endthanks
\thanks The second and fourth authors were supported in part
by NSF Grant DMS--9505024 \endthanks
\thanks The third author was supported in part by NSF Grant DMS--9504913
\endthanks


\date 27 August 1998 
\enddate

\toc
\head 1. Introduction\page{1}\endhead
\head 2. Definitions and the main results\page{3}\endhead
\head 3. The topology of the teardrop\page{7}\endhead
\head 4. The teardrop of an approximate fibration\page{14}\endhead
\head 5. Spaces of stratified neighborhoods and manifold approximate
fibrations\page{22}\endhead
\head 6. Homotopy near the lower stratum\page{29}\endhead
\head 7. Higher classification of stratified neighborhoods\page{37}\endhead
\head 8. Examples of exotic stratifications\page{43}\endhead
\head 9. Extensions of isotopies and $h$--cobordisms\page{55}\endhead
\specialhead {} References\page{60}\endspecialhead
\endtoc

\abstract
We develop a theory of tubular neighborhoods for the lower strata
in manifold stratified spaces with two strata. In these topologically
stratified spaces, manifold approximate fibrations and teardrops
play the role that fibre bundles and mapping cylinders play in
smoothly stratified spaces. Applications include the classification
of neighborhood germs, the  construction of
exotic stratifications, a multiparameter isotopy extension theorem
and an $h$--cobordism extension theorem.
\endabstract

\endtopmatter
\document

\head 1. Introduction
\endhead
The question that motivates this paper is a basic one:  suppose that
one has a locally flat topological submanifold of a manifold, what kind of
geometric structure describes the neighborhood?

For smooth manifolds the entirely satisfactory answer is given by the
tubular neighborhood theorem which identifies neighborhood germs with
vector bundles.  In the piecewise linear category, one has the theory of block
bundles \cite{51}.  For the topological category, the situation
is much messier: essentially one can classify the neighborhoods without
really describing them  (see \cite{52}).

The answer that we give is in terms of a variant of the notion of a
fiber bundle, the manifold approximate fibration (MAF).  
While fiber bundles are
maps with identifications of the inverse images of points, MAFs are
essentially maps with identifications of the inverse images of open balls.
At the level of definitions, they are to fiber bundles what cell-like 
maps are to
homeomorphisms.  However, unlike the cell-like case, they cannot always be
approximated by bundles (or even block bundles) and represent a genuinely
more general notion.  Happily, though, one has a good control of the theory
of MAFs, see \cite{29}, \cite{30}.

A special case of our theorem asserts that the (space of) $d+n$
dimensional locally flat germ neighborhoods of an $n$--manifold 
$M^n$ are (is
homotopy equivalent to the space of) MAFs mapping to $M\times\br$, 
with the inverse
images of small balls in $M\times\br$ homeomorphic to $S^{d-1}\times\br^{n+1}$.
One should think
of a MAF mapping to $M\times\br$ as having as domain a deleted 
neighborhood of $M$ and as consisting
of two pieces: the first is the projection of generalized tubular
neighborhood bundle, and the second is the radial direction, e.g. something
like distance from the submanifold.  We call this structure a `teardrop
neighborhood.'

Actually, though, our paper is written in more generality.  It gives an
analysis of neighborhoods of the singular stratum of a stratified space as in
\cite{48} which has only two strata.  This means that our results apply, for
instance to quotients of semifree group actions, and leads to new results
for these.

The description of germ neighborhoods is good enough to recover and
reprove Quinn's isotopy and homogeneity theorems, and go rather further:
we obtain multiparameter isotopy extension theorems, which lead to local
contractibility of homeomorphism groups for such spaces.

Another important application is to complete (in the two stratum
case) the $h$--cobordism theorem given in \cite{48}.  That paper provides an
invariant whose vanishing is necessary and sufficient for a stratified 
$h$--cobordism to be a product.  We give the realization: any element
in the appropriate Whitehead group can be realized by a stratified
$h$--cobordism.

The picture we give of stratified spaces, when combined with the
analysis of MAFs in \cite{29} and the stable homeomorphism groups in 
\cite{63}, is
more than fine enough to be used to give an independent 
proof of the two stratum case of
the stratified classification results in \cite{62}.  
However, the current approach
is more directly geometric, which has at least two important advantages.
The first is that the analysis is done here unstably: i.e. without first
crossing with Euclidean spaces and then removing them.   We apply this to
give examples of stratified spaces (even topological locally linear
orbifolds) where no amount of Euclidean stabilization allows one to
construct a block bundle neighborhood of the lower stratum.

The other main advantage is that of canonicity, which is important
for the multiparameter results discussed above, and also plays a key role in
relating the splitting results for spaces of MAFs over Hadamard manifolds
proven in \cite{32}, and the Novikov rigidity results proven by Ferry and
Weinberger (see \cite{16}, \cite{17}) 
for stratified spaces with nonpositively curved
strata.  These seemingly different results are essentially equivalent after
taking a loop space.

Finally, these results form the bottom of an induction that leads to
extensions of all of the theorems and applications mentioned above to
general stratified spaces with an arbitrary number of strata (see \cite{26},
\cite{27}).


\head 2. Definitions and the main results
\endhead
Quinn \cite{48} has proposed a setting for the study of
those spaces admitting purely topological
stratifications as distinct from the smooth
stratifications of Whitney \cite{65}, Thom \cite{58}, Mather \cite{40} and
others (cf\. \cite{18}).  In this paper we consider spaces
$X$ containing a manifold $B$ such that the pair
$(X,B)$ is a manifold homotopically stratified set in
the sense of Quinn.  We call $X$ a manifold stratified space
with two strata.  Roughly, this means that $X\setminus B$ is a
manifold, $B$ satisfies a tameness condition in $X$,
and there is a good homotopy model for a normal
fibration of $B$ in $X$.

We begin by recalling the definitions relevant to the
manifold stratified spaces.  Most of these concepts
can be found in Quinn \cite{48} and Weinberger \cite{62}, but our
terminology is not consistent with either source.
Moreover, since we are only dealing with stratified
spaces with two strata, our definitions are
specialized to that case.

Let $(X,A)$ be a pair of spaces so that $A\subseteq X$.  
Then $X$ is said to
have two {\it strata\/}: the lower (or bottom) stratum
$A$ and the top stratum $X\setminus A$.  If $(Y,B)$ is another
pair, then a map $f: (X,A)\to (Y,B)$ is said to be
{\it strict\/}, or {\it stratum-preserving\/}, 
if $f(X\setminus A)\subseteq Y\setminus B$ and $f(A)\subseteq B$.
The subspace $A$ of $X$ is said
to be {\it forward tame\/} if there exists a neighborhood
$N$ of $A$ in $X$ and a strict map $H: (N\times I,\
A\times I\cup N\times \{ 0 \} )\to (X,A)$ such that 
$H(x,t) = x$ for all $(x,t)\in A\times I$  and $H(x,1)
= x$ for all $x\in N$.  In this case, $H$ is called a
{\it nearly strict deformation\/} of $N$ into $A$.

Let $\mapst ((X,A), (Y,B))$ denote the space of strict
maps with the compact-open topology.  The {\it homotopy
link\/} of $A$ in $X$ is 
$$\ho(X,A) = \mapst (([0,1],\{ 0\} ),\ (X,A)).$$  
Evaluation at 0 defines a map $q:
\ho(X,A)\to A$ which should be thought of as a
model for a normal fibration of $A$ in $X$.  
A point inverse $q^{-1}(x)$ is the {\it local homotopy link\/}
(or {\it local holink\/}) at $x\in A$.
In the
case that $X$ is an $n$-manifold and $A$ is a locally
flat submanifold of dimension $i$, then Fadell proved
that $q:\ho(X,A)\to A$ is a fibration with
homotopy fibre $S^{n-i-1}$ and used the homotopy link as a
substitute in the topological category for tubular
neighborhoods in the differential category (see \cite{12},
\cite{44}, \cite{24}, \cite{25}, \cite{28\rm, App\. B}.)

The pair $(X,A)$ is said to be a {\it homotopically
stratified pair\/} if $A$ is forward  tame in $X$ and if $q:
\ho(X,A)\to A$ is a fibration.  If in addition,
the fibre of $q:\ho(X,A)\to A$ is finitely dominated, then $(X,A)$ is said to
be  {\it homotopically stratified with finitely dominated
local holinks}. 
(When we say that the fibre of $q$
is finitely dominated and $A$ is not path connected, we mean that each
fibre of $q$ is finitely dominated.)
If 
the strata $A$ and $X\setminus A$ are manifolds 
(without boundary), $X$ is a
locally compact separable metric space, and $(X,A)$ is 
homotopically stratified with finitely dominated
local holinks, then $(X,A)$
is a {\it manifold stratified pair}.

We now define the set of equivalence classes of neighborhoods which is the
main object of study in this paper.  Let $B$ be an
$i$-manifold (without boundary) and let $n\ge 0$ be a
fixed integer.  A {\it germ of a stratified
neighborhood\/} of $B$ is an equivalence class
represented by a manifold stratified pair $(X,B)$ with
$\dim(X\setminus B) = n$.  Two such pairs $(X,B)$ and $(Y,B)$ are {\it
germ equivalent} provided that there exist open
neighborhoods $U$ and $V$ of $B$ in $X$ and $Y$,
respectively, and a homeomorphism $h: U\to V$ such
that $h|B = \id_B$.  In this paper we will classify
stratified neighborhoods of $B$ up to germ equivalence
(provided $n\ge 5$).  The basic
construction which makes this possible is now described.

Let $p :X\to Y\times \br$ be a map.  The {\it
teardrop\/} of $p$, denoted $X\cup_p Y$, is the space
with underlying set the disjoint union $X\amalg Y$ and natural
topology defined in \S3 below.
We are interestested in  those maps $p$ with the property
that $(X\cup_p Y, Y)$ is a manifold stratified or homotopically
stratified pair.

Recall that an {\it approximate fibration\/} is a map
with the approximate homotopy lifting property (see Definition 4.5) and
that a map $p: X\to Y$ is a {\it manifold approximate
fibration\/} if $p$ is an approximate fibration, $p$ is
proper, and $X$ and $Y$ are manifolds (without
boundary) (see e\.g\. \cite{29}).  Two maps $p: X\to Y$
and $p' : X'\to Y$ are {\it controlled
homeomorphic\/} if there is a homeomorphism $h: \cyl
(p)\to \cyl (p')$ between mapping cylinders such
that $h|Y= \id_Y $ which is {\it level\/} in the sense
that $h$ commutes with the natural projections to
$[0,1]$.  In \cite{29} manifold approximate fibrations
over $Y$ with total space of dimension greater than
four are classified up to controlled homeomorphism.

The main results can now be stated.  Let $n\ge 5$ be a
fixed integer and let $B$ be a closed
manifold.
In the general setting  of
manifold stratified pairs $(X,B)$, neighborhoods of 
$B$ in $X$ need not have nice geometric structure.
For example, $B$ need not be locally conelike in $X$
and $B$ may even fail to have mapping cylinder
neighborhoods (locally or globally).  However, 
the first theorem says that the lower stratum in a manifold stratified
pair has a neighborhood which is the teardrop of a manifold
approximate fibration.
The second theorem is just a more complete
statement.

\proclaim{Theorem 2.1 (Teardrop Neighborhood Existence)}
Let\/ $(X,B)$ be a pair such that $X\setminus B$
is a manifold of dimension $n$. Then\/  $(X,B)$ is a manifold stratified
pair if and  only if 
$B$ has a neighborhood in $X$ which is the teardrop of a manifold approximate
fibration.
\endproclaim

There are two equivalent ways to understand what it means for $B$ to have
a neighborhood in $X$ which is the teardrop of a manifold approximate fibration
as in Theorem 2.1:
\widestnumber\item{(ii)}
\roster
\item"(i)" There exist  a neighborhood $U$ of $B$ in $X$ 
and a manifold approximate fibration $p:V\to B\times\br$ such that
$(U,B)$ is homeomorphic 
to $(V\cup_pB,B)$ \rel $B$.

\item"(ii)" There exists an open neighborhood $U$ of $B$ in $X$
and a proper map $f:U\to B\times (-\infty,+\infty]$ such that
$f^{-1}(B\times\{+\infty\}) = B$,
$f|:B\to B\times\{+\infty\}$ is the identity, and
$f|:U\setminus B\to B\times\br$ is a manifold approximate fibration.
\endroster
That these are equivalent follows from the material in \S3 (see especially
Proposition 3.7). 
Theorem 2.1 follows directly from the following theorem.

\proclaim{Theorem 2.2 (Neighborhood Germ Classification)} 
The teardrop construction defines a bijection from the
set of controlled homeomorphism classes of manifold
approximate fibrations over $B\times \br$ \rom(with
total space of dimension $n$\rom) to the set of germs of
stratified neighborhoods of $B$ \rom(with top stratum of
dimension $n$\rom). \endproclaim

In fact, Theorem 2.2 is just the consequence at the 
$\pi_0$ level of a more general Higher Classification Theorem which asserts
that two simplicial sets are homotopy equivalent (Theorem 2.3 below).
However, a proof of Teardrop Neighborhood Existence (Theorem 2.1) 
is offered in
\S7 which avoids some of the parametric considerations needed
for Theorem 2.3.
Before we can define the simplicial sets appearing in Theorem 2.3 we need
sliced versions of some of the definitions. 

Let $\Delta $ be a space which will play the role of a
parameter space.  Let $(X,A\times \Delta )$ be a pair
of spaces and let $\pi : X\to\Delta $ be a map such
that $\pi|: A\times \Delta\to\Delta $ is the projection.
Then $A\times \Delta $ is said to be {\it sliced forward
tame\/} in $X$ (with respect to $\pi$) if there exists a
neighborhood $N$ of $A\times \Delta $ in $X$ and a
nearly strict deformation $H$ of $N$ into $A\times
\Delta$ such that $H$ is fibre preserving over
$\Delta $ (i\.e\., $\pi H_t = \pi $ for all $t\in I)$.
The {\it sliced homotopy link\/} of $A\times \Delta $ in
$X$ (with respect to $\pi $) is $\ho_\pi (X,
A\times \Delta ) = \{\omega\in\mapst(([0,1],
\{ 0\}), (X, A\times \Delta )) ~|~ \pi\omega(t) = \pi\omega(0)$
for all $t\in I\}$.  Note that evaluation at 0 still
gives a map $q:\ho_\pi  (X,A\times \Delta) \to
A\times \Delta $.

Let $n\ge 0$ be a fixed integer and let $B$ be a manifold
(without
boundary). In \S5 the simplicial set $\SN^n(B)$ of
stratified neighborhoods of $B$ is defined.  Roughly,
its $k$-simplices are $k$-parameter families of
manifold stratified spaces containing $B\times \Delta ^k$ as
the lower stratum using the notions of sliced forward
tameness and the sliced homotopy link.  On the other
hand, the simplicial set $\maf^n(B\times \br )$ of
manifold approximate fibrations over $B\times \br
$ was defined in \cite{29} (see also \S5).  This set has
$k$-simplices consisting of $k$-parameter families of
manifold approximate fibrations over $B\times \br$.

Note that if $p: M\to B\times \br \times \Delta^k$
is a map, then the teardrop construction yields a pair
$(M\cup_p B\times \Delta ^k, B\times \Delta ^k)$.
Define $\Psi (p) = (M\cup _p B\times \Delta ^k, B\times\Delta
^k )$.  The following result is the simplicial set
version of Theorem 2.2.

\proclaim{Theorem 2.3 (Higher Classification)} 
If $B$ is a closed manifold and $n\ge 5$,
then the teardrop construction defines a homotopy
equivalence $\Psi : \maf^n(B\times \br )\to\SN^n
(B)$.  \endproclaim

To see why Theorem 2.2 follows from Theorem  2.3,
recall  that $\pi_0 \maf^n (B\times
\br)$ is the set of controlled homeomorphism
classes of manifold approximate fibrations over
$B\times \br$ (see \cite{29}).  And it is not difficult to see
that $\pi_0 \SN^n (B)$ is the set of germs of
stratified neighborhoods of $B$ (see Corollary 5.6).

Fibre bundles have well-defined fibres up to homeomorphism. Analogously,
manifold approximate fibrations have well-defined fibre germs up to
controlled homeomorphism (see \cite{29}).
Recall that if $p:M\to B$ is a manifold approximate fibration with
$B$ connected, $\dim B=i$ and $\dim M=n\geq 5$, then the {\it fibre
germ\/} of $p$ is the manifold approximate fibration
$q=p|:V=p^{-1}(\br^i)\to\br^i$ where $\br^i\hookrightarrow B$ is an
open embedding (which is orientation preserving if $B$ is oriented).
The theorems above involve manifold approximate fibrations
$p:M\to B\times\br$ and these have fibre germs of the form
$q:V\to\br^{i+1}$.
The teardrop construction yields a manifold stratified pair
$(V\cup_q\br^i,\br^i)\subseteq (M\cup_pB,B)$.
The local holink of $B$ in $M\cup_pB$ is homotopy equivalent to
$V$. For locally conelike stratified pairs $(X,B)$ (see \cite{55})
a neighborhood of $B$ in $X$ is given by the teardrop of a manifold
approximate fibration
$p:M\to B\times\br$ with {\it trivial fibre germ}; that is, the 
projection $F\times\br^{i+1}\to\br^{i+1}$ for some closed manifold $F$.

Let $\maf(B\times\br)_q$ be the simplicial subset of $\maf^n(B\times\br)$
consisting of manifold approximate fibrations with fibre germ
$q:V\to\br^{i+1}$. For trivial fibre germ, we write this simplicial set
as $\maf(B\times\br)_{F\times\br^{i+1}}$.
According to \cite{29}, \cite{30}, $\maf(B\times\br)_q$ is homotopy equivalent to
a simplicial set of lifts
of 
$B\to\btop_{i+1}$ up to
$\btop^{\level}(q)$
where $B\to\btop_{i+1}$ is the composition of the classifying map
$B\to\btop_i$ for the tangent bundle of $B$ with the map
$\btop_i\to\btop_{i+1}$ induced by euclidean stabilization.
The fibre of
$\btop^{\level}(q)\to\btop_{i+1}$ is $\btop^c(q)$, the 
classifying space of controlled homeomorphisms on $q:V\to\br^{i+1}$.
According to \cite{31} $\btop^c(q)\simeq\btop^b(q)$, the classifying
space of bounded homeomorphisms. In the case of trivial fibre germ
$F\times\br^{i+1}\to\br^{i+1}$, this is written as $\btop^b(F\times
\br^{i+1})$.
For relevant information about the homotopy type of $\btop^b(F\times
\br^{i+1})$ see \cite{63}.
For example, if $B\times\br$ is parallelizable, then
$$\maf(B\times\br)_{F\times\br^{i+1}}\simeq 
\map(B,\btop^b(F\times\br^{i+1}))$$
and this classifies neighborhood germs in the locally conelike case.

These classification results together with \cite{63} can be used to give
an alternative proof of Weinberger's surgery theoretic stable classification
theorem \cite{62} in the case of two strata. In fact, this alternative
proof is outlined in \cite{62\rm, 10.3\.A} and discussed here in Remark 8.17(i).
 
In \S8 the classification theorem for manifold approximate fibrations
is combined with the classification of neighborhood germs to construct
examples of manifold stratified pairs in which the lower strata do not
have a neighborhood given by the mapping cylinder of a fibre bundle, 
or even a block bundle. Moreover, the examples do not improve in 
this regard under euclidean stabilization. These examples are locally
conelike and the lower strata do have neighborhoods which are mapping
cylinders of manifold approximate fibrations.

In addition, Theorem 2.2 provides the link between the results 
on approximate fibrations proven in \cite{32} and the tangentiality results 
of \cite{16}, \cite{17}.

Teardrop neighborhoods can also be used in conjunction with the geometric
theory of manifold approximate fibrations \cite{22}, \cite{24} to study the
geometric topology of manifold stratified pairs. We include two examples
here, both of which involve extending a structure on the lower stratum to
a neighborhood of the stratum. This is a very important use of
manifold approximate fibrations which is similar to the way
fibre bundles are used in inductive proofs for smoothly stratified
spaces. The following isotopy extension theorem is established in \S9.

\proclaim{Corollary 2.4 (Parametrized Isotopy Extension)} 
If $(X,B)$ is a manifold stratified pair, $\dim X \ge 5$,
$B$ is a closed manifold 
and $h: B\times \Delta ^k\to B\times\Delta ^k$ is a
$k$-parameter isotopy \rom(i\.e\., $h$ is a homeomorphism,
fibre preserving over $\Delta ^k$, and $h|B\times
\{ 0 \} = \id_{B\times\{ 0\}}$\rom), then there exists a $k$-parameter isotopy
$\tilde h : X\times \Delta ^k \to X\times \Delta ^k$
extending $h$ such that $\tilde h$ is the identity on
the complement of an arbitrarily small neighborhood of $B$.
\endproclaim

In the case that $B$ is a locally flat submanifold of
$X$, this theorem is due to Edwards and Kirby \cite{11}.
For locally conelike stratified spaces with an
arbitrary number of strata, it is due to Siebenmann
\cite{55}.  Finally, Quinn \cite{48} proved this theorem for
manifold stratified spaces in general (with an arbitrary number
of strata), but only in the case $k = 1$.

Also in \S9 we prove an $h$--cobordism extension theorem which can be
used to prove a realization theorem for stratified Whitehead torsions (see
Remark 9.4(i)).

A {\it fibre preserving map \rom(f\.p\.\rom)\/} is a map which preserves the fibres
of  maps to a given parameter space. The parameter space will usually
be a $k$-simplex or an arbitrary space denoted $K$. Specifically,
if $\rho:X\to K$ and  $\sigma:Y\to K$ are maps, then a map $f:X\to Y$
is f\.p\. (or f\.p\. over $K$) if $\sigma f=\rho$.

There is a notion of reverse tameness which, in the presence of
forward tameness, is often equivalent to
the finite domination of local holinks condition discussed above.
See \cite{48\rm, 2.15} and \cite{28\rm, 9.15, 9.17, 9.18} paying special attention to
the point-set topological conditions appearing in \cite{28}.
Moreover, when strata are manifolds, the notions of forward tameness and
reverse tameness are often equivalent (by Poincar\'e duality).
See \cite{48\rm, 2.14} and \cite{28\rm, 10.13, 10.14} 
paying special attention to the $\pi_1$
conditions appearing in \cite{28}.

Hughes and Ranicki's book \cite{28} 
contains many of the the results of this paper in the special
case of stratified pairs with lower stratum a single point.
The reader is advised to consult that work for background, examples and
historical remarks.
The paper \cite{27} contains generalizations to manifold stratified
spaces with more than two strata. The proofs in \cite{27} are often by 
induction on the number of strata and rely on the present 
paper for the beginning
of the induction. More applications to the geometric topology of
manifold stratified spaces are contained in \cite{27}. See also \cite{26}.

\head 3. The topology of the teardrop
\endhead
Let $p: X\to Y\times \br$ be a map.  The {\it teardrop\/}
of $p$, denoted by $X\cup_p Y$, is defined to be the
space with underlying set the disjoint union $X \amalg
Y$ and topology given as follows.  First, let $c:
X\cup_p Y\to Y\times (-\infty , +\infty ]$ be defined
by
$$
c(x) = \cases
     p(x),  & \text{if  $x\in X$}\\
     (x,+\infty ), & \text{if $x\in Y$}.
\endcases$$
Then the topology on $X\cup_p Y$ is the minimal
topology such that
\widestnumber\item{(ii)}
\roster
\item"(i)" $X\subseteq X\cup_p Y$ is an open embedding, and

\item"(ii)" $c$ is continuous.
\endroster
The mapping $c$ is called the {\it collapse\/} mapping for
$X\cup_p Y$.

Note that a basis for this topology is given by 
$$\{ c^{-1} (U)~|~ U \text{is open in }\
     Y\times (-\infty ,+\infty ] \} \cup
     \{ U~|~ U\ \text{is open in} \  X\} .$$

There are two minor variations on this construction
which we will use.  The first occurs when $U$ is an
open subset of $X$ and $p$ is only defined on $U$,
$p: U\to Y\times \br $.  Then we let $X\cup_p Y =
X \cup (U\cup _p Y)$.  The second variation occurs
when the range of $p$ is restricted, usually to
$Y\times [0, +\infty )$.  We can still form $X\cup_p
Y$ and the collapse map $c: X\cup_p Y\to Y\times [0,
+\infty ]$.

Special cases and variations of the teardrop
construction have appeared frequently in the
literature and we now discuss some examples.

\subhead 3.1. Mapping cylinders \endsubhead
If $q: X\to Y$ is a map, let $p: X\times (0,1)\to
Y \times (0,1)$ denote $q\times\id$. Then we define
the {\it open mapping cylinder\/} of $q$ to be the teardrop
$$\ocyl(q) = (X\times (0,1)) \cup_p Y,$$
where we replace $\br $ with $(0,1)$.
The {\it mapping cylinder\/} is
$$\cyl(q) = (X\times [0,1)) \cup_p Y.$$
Note that this is not the usual quotient topology on
the mapping cylinder (except in special cases), but is
more useful geometrically (see \cite{6}, \cite{45}, \cite{48}).
The {\it open cone\/} $\oc(X)$ of a space $X$ is just the open
mapping cylinder (with the teardrop topology) of the constant map
$X\to\{ v\}$ with $v$ the vertex of the cone.

It follows from this example that the teardrop $X\cup_pY$ of a map
$p:X\to Y\times [0,1)$ is a mapping cylinder neighborhood of $Y$ if
there exist a space $Z$, a map $q:Z\to Y$, and a homeomorphism
$h:Z\times[0,1)\to X$ such that $ph = q\times\id_{[0,1)}$.

\subhead 3.2. Joins\endsubhead  
The join of two spaces $X\ast
Y$ can be viewed as a teardrop as follows.  Let $p:
X\times (0,1)\times Y\to Y \times (0,1)$ be defined
by $p(x,t,y) = (y,t)$.  
Identify $X\times (0,1)$ with $\oc
(X)\setminus\{ v\}$.  Then $X\ast Y = (\oc(X)\times Y)\cup
_p Y$.  Again, this is not the quotient topology, but
it is a topology which is often used.

\subhead 3.3. Hadamard's teardrop\endsubhead  
Let $H$ be an 
Hadamard manifold of dimension $n$ (i\.e\., $H$ is a
complete, simply connected Riemannian manifold of
nonpositive curvature) with distance function $d$
induced by the metric.  Fix a point $x_0 \in H$ and
let $S$ denote the unit tangent sphere of $H$ at $x_0$.
For each $x\not= x_0$ in $H$, let $\gamma _x :
[0,+\infty ) \to H$ be the unique unit speed geodesic
such that $\gamma_x (0) = x_0$ and $\gamma_x (d(x_0,x))
= x$.  Define $p: H\setminus\{ x_0\} \to S\times (0, +
\infty )$ by
$$p(x) = (\gamma'_x (0), d(x_0, x)).$$
(It follows from standard facts that $\gamma'_x
(0)$ depends continuously on $x$.)  It is easy to see
that the teardrop $H\cup_p S$ is homeomorphic to the
Eberlein-O'Neill compactification $\ov H = H\cup
H(\infty )$ with the cone topology \cite{10} (in
particular, $H\cup_pS$ is an $n$-cell).  To see this,
let $f: [0,1] \to [0,+\infty ]$ be a homeomorphism,
let $B$ be the unit tangent ball of $H$ at $x_0$ and
let $\psi : B\to H\cup_p S$ be defined by
$$\psi (v) = \cases
         \exp (f(\Vert v\Vert\cdot v) & \text{if $x\notin S$}\\
          v   & \text{if $x\in S$}.
\endcases$$
Then $\psi $ is a homeomorphism (using the continuity
criterion below) and together with \cite{10\rm, Prop\. 2.10}
can be used to get a homeomorphism with $\ov H$.

Another useful construction is as follows.  If $q: M
\to H$ is a map, then the composition $pq:M\setminus
q^{-1} (x_0)\to S\times (0,+\infty )$ yields a
teardrop $M\cup_{pq} S$.  If $q$ is proper, this
amounts to compactifying $M$ by adding the sphere
$S\approx H(\infty )$ at infinity.  This special case
of the teardrop was used in \cite{31} for studying
manifold approximate fibrations over $H$.

\subhead Point-set topology \endsubhead
A pleasant feature of the teardrop topology is that it
is easy to decide when a function into a teardrop is
continuous.  In fact, the proof of the following lemma
follows immediately from the description of the basis above.

\proclaim{Lemma 3.4 (Continuity Criteria)} 
Let $f: Z\to X\cup _p Y$ be a function.  Then $f$ is
continuous if and only if 
\widestnumber\item{(ii)}
\roster
\item"(i)" $f|: f^{-1} (X)\to X$ is continuous, and

\item"(ii)" the composition 
$Z@>f>> X\cup_p Y @>c>> Y\times (-\infty , +\infty ]$
is continuous.\qed
\endroster
\endproclaim

If $(X,Y)$ is a pair of spaces, we now address the
question of the existence of a map $p: X\setminus Y\to Y\times
\br $ such that the identity from $X$ to
$(X\setminus Y)\cup _p Y$ is a homeomorphism.  If this is the
case, then $(X,Y)$ is said to be {\it the teardrop of\/}
$p$. The answers are in Corollaries 3.11 and 3.12.

If $f: X\to Y$ is a map and $A\subseteq Y$, then $f$ is
said to be a {\it closed mapping over\/} $A$ if for each
$y\in A$ and closed subset $K$ of $X$ such that $K\cap
f^{-1} (y) = \emptyset $, it follows that $y\notin \cl(f(K))$
(the closure of $f(K))$.

\remark{Remark \rom{3.5}}
\widestnumber\item{(iii)}
\roster\item"(i)" $f: X\to Y$ is a closed mapping if and only
if $f$ is a closed mapping over $Y$.

\item"(ii)" If $A\subseteq Y$ and $f: X\to Y$ is a
closed mapping over $A$, then $f$ is a closed mapping
over any $B\subseteq A$.

\item"(iii)" If $A$ is closed in $Y$ and $f: X\to Y$
is a closed mapping over $A$, then $f|:
f^{-1}(A)\to A$ is a closed mapping (but not conversely).
\endroster
\endremark

\proclaim{Lemma 3.6} 
If $p: X\to Y\times \br $
is a map, then the collapse $c: X\cup_p Y\to Y\times
(-\infty , + \infty ]$ is a closed mapping over
$Y\times \{ +\infty \}$. \endproclaim

\demo{Proof}  Let $y\in Y$ and let $K$ be a closed
subset of $X\cup_p Y$ such that $y\notin K$ (note $y
= c^{-1} (y,+ \infty ))$.  Then $y\in U = (X\cup_p Y)\setminus K$ 
and $U$ is open.  By the definition of the
teardrop topology, there is an open subset $V$ of $(y,
+ \infty )$ in $Y\times (-\infty , + \infty ]$ such
that $y\in c^{-1} (V)\subseteq U$.  Then $c(K)\cap V =
\emptyset $, so $(y,+ \infty ) \notin \cl(c(K))$.
\qed
\enddemo

\proclaim{Proposition 3.7}
Let $(X,Y)$ be a pair of
spaces for which there is a mapping $f: X\to Y\times
(-\infty , + \infty ]$ such that $f(y) = (y,+\infty )$
for each $y\in Y$ and $f(X\setminus Y)\subseteq Y\times \br
$.  Let 
$$p = f|: X\setminus Y\to Y\times \br .$$  
Then
$(X,Y)$ is the teardrop of $p$ if and only if $f$ is a
closed mapping over $Y\times \{ + \infty \}$. \endproclaim

\demo{Proof}  First note that $f$ is the collapse
$c$ for the teardrop $(X\setminus Y)\cup_p Y$.  It follows that
the identity $X\to (X\setminus Y)\cup_p Y$ is always continuous.
To prove the proposition, assume that the identity is
a homeomorphism.  By Lemma 3.6, $c$ is a closed
mapping over $Y \times \{ +\infty \}$.  Since $f= c$,
so is $f$.

\noindent
Conversely, assume $f$ is a closed mapping over
$Y\times \{ + \infty \}$.  Given an open subset $U$ of
$X$, we will show that $U$ is open in $(X\setminus Y) \cup_p
Y$.  For this, it suffices to consider $y\in U\cap Y$
and show that $U$ is a neighborhood of $y$ in $(X\setminus Y)
\cup _p Y$.  To this end let $K = X\setminus U$ and observe
that since $f^{-1} (y,+\infty ) = y \notin K$, it
follows that $(y,+\infty )\notin \cl(f(K))$.  Thus,
there is an open subset $V$ of $Y\times (-\infty ,
+\infty ]$ such that $(y,+\infty ) \in V$ and $V\cap
f(K) = \emptyset $.  Then $c^{-1} (V)$ is open in
$(X\setminus Y)\cup _p Y$ and $y\in c^{-1} (V)\subseteq U$.
\qed
\enddemo

\proclaim{Corollary 3.8} 
A pair $(X,Y)$ is a teardrop
if and only  if there is a map $f: X\to Y\times
(-\infty , + \infty ]$ which is closed over $Y\times
\{ +\infty \}$ such that $f(y) = (y,+\infty )$ for
each $y\in Y$ and $f(x) \in Y\times \br $ for
each $x\in X\setminus Y$.
\qed
\endproclaim

\proclaim{Proposition 3.9} Let $(X,Y)$ be a pair of
spaces such that $X$ is Hausdorff and $Y$ is locally
compact.  Suppose there exist a proper retraction $r:
X\to Y$ and a map $\phi : X\to (-\infty , +\infty ]$
such that $\phi^{-1} (+\infty ) = Y$.  Then $f =
r\times \phi : X\to Y\times (-\infty , +\infty ]$ is
a closed mapping over $Y\times \{ +\infty \}$.
Consequently, $(X,Y)$ is a teardrop.  \endproclaim

\demo{Proof}  Let $y\in Y$ and let $K$ be a closed
subset of $X$ such that $y\notin K$.  We need to show
that $(y,+\infty ) \notin \cl(f(K))$.  To this end,
let $U$ be open in $X$ such that $y\in U$ and $U\cap K
= \phi $.  Choose an open subset $V$ of $Y$ such that
$y\in V, \cl(V)\subseteq U\cap Y$, and $\cl(V)$ is compact.
Let $K_1 = r^{-1} (\cl(V))\cap K$ and $K_2 =
K\setminus r^{-1}(V)$.  Then $K_1$ is compact and $K = K_1\cup
K_2$.  Since $f(K_1)$ is compact and $(y,+\infty )\notin
f(K_1)$, it suffices to show that $(y,+\infty )\notin
\cl(f(K_2))$.  But $(y,+\infty ) \in V\times (-\infty
, +\infty ]$ and $f(K_2) \cap V\times (-\infty , +
\infty ] = \emptyset$.
That $(X,Y)$ is a teardrop follows form Proposition 3.7.
\qed
\enddemo

Note that such a map $\phi$ in the hypothesis of Proposition
3.9 would exist whenever $X$ is normal and $Y$ is a closed
$G_\delta $-subset.

\proclaim{Theorem 3.10} Let $Y$ be a closed subset of
the metrizable space $X$.  Then $(X,Y)$ is a teardrop
if and only if there exists a metric $d$ for $X$ and a
retraction $r: X\to Y$ such that whenever $\{ x_n\}$
is a sequence in $X$ with $x_n \to \infty $ \rom(i\.e\., $\{ x_n\}$ has
no convergent subsequence\rom) and $d(x_n, Y)\to
0$, it follows that $r(x_n)\to \infty $.  \endproclaim

\demo{Proof}  Suppose first the $(X,Y)$ is the
teardrop of $p: X\setminus Y\to Y\times \br $ and let
$c: X\to Y\times (-\infty , +\infty ]$ be the collapse.
Define $\rho : X\to [0, + \infty )$ to be the composition
$$X @>c>> Y \times (-\infty , + \infty ]
     @>\proj>> (-\infty , + \infty ]
     @>h>> [0, +\infty )$$
where $h$ is a homeomorphism.  Let $D$ be any metric
on $X$ and define $d$ by
$$d(x,x') = D (x,x') + \vert \rho (x)-\rho (x')\vert .$$
It is easy to see that $d$ is indeed a metric and
yields the same topology on $X$ as $D$.

\noindent
Define $r: X\to Y$ to be the composition
$$X @>c>> Y \times (-\infty , + \infty ]
     @>\proj>> Y.$$
To see that $r$ has the desired property, let $\{
x_n\}$ be a sequence in $X$ such that $x_n \to \infty
$ and $d(x_n, Y)\to 0$.  Given $y\in Y$ we will show
that there is no subsequence $\{ x_{n_k}\}$ with $r(x_{n_k})
\to y$.  To this end let
$$K = \bigcup_{n=1}^\infty\{ x_n\} \setminus \{ y\} .$$
Then $K$ is a closed subset of $X$ and $y\notin K$.
Since $c$ closed over $Y \times \{ + \infty \}$ by
Lemma 3.6, it follows that $(y, + \infty ) \notin\cl
(c(K))$.  Thus, if $\{ x_{n_k}\}$ is a subsequence,
$\{ c(x_{n_k})\}$ does not converge to $(y,+\infty )$.
Since $d(x_n, Y)\to 0,\ \rho (x_n) \to 0$.  This
implies $c(x_n)\to Y\times \{ +\infty \}$.  If $r(x_{n_k})
\to y$, then we would have $c(x_{n_k})\to (y,+\infty
)$, a contradiction.

\noindent
Conversely, assume $r$ and $d$ are given as above.  Define
$\phi : X\to (-\infty , + \infty ]$ by 
$$\phi (x) = \cases
      {1\over d(x,Y)},  & \text{if $x\in X\setminus Y$}\\
         + \infty  & \text{if  $x\in Y$}.
\endcases$$
Let $f = r\times \phi : X\to Y\times (-\infty , +\infty ]$.
By Corollary 3.8, it suffices to show that $f$ is
closed over $Y\times \{ + \infty \}$.  To this end let
$K$ be closed in $X$ and $y\in Y\setminus K$.  Suppose $(y,
+\infty )\in \cl(f(K))$.  Then there exists a sequence
$\{ x_n\}$ in $K$ such that $f(x_n)\to (y, +\infty
)$.  Then $r(x_n)\to y$ and $\phi (x_n)\to +\infty
$.  Thus, $d(x_n, Y)\to 0$.  If $\{ x_n\}$ has a
convergent subsequence $\{ x_{n_k}\}$, then
$x_{n_k}\to y _0\in Y\cap K$.  Then $r(x_{n_k}) \to y_0$
so $y = {y_0}$, a contradiction since $y\notin K$.
Thus, we must have $x_n\to \infty $.  So $r(x_n)
\to \infty $, again a contradiction.
\qed
\enddemo

\proclaim{Corollary 3.11}
If $Y$ is a compact subset of the
metric space $X$, then $(X,Y)$ is a teardrop if and
only if there exists a retraction $r: X\to Y$.
\qed
\endproclaim

\proclaim{Corollary 3.12} Let $Y$ be a closed subset
of the locally compact
metric space $X$.
Then $(X,Y)$ is a teardrop if and only if there exists
a retraction $r : X\to Y$.
\endproclaim

\demo{Proof}  If $(X,Y)$ is a teardrop,  let
$r$ be given by Theorem 3.10.  
Conversely, if $r:X\to Y$ is a retraction then by Proposition 3.9, it
suffices to show that $Y$ has a closed neighborhood $N$ in $X$
such that $r|:N\to Y$ is proper. To this end, for each $y\in Y$, let
$N_y$ be a compact neighborhood of $y$ in $X$ and let
$$N = \bigcup\{r^{-1}(N_y\cap Y)\cap N_y ~|~ y\in Y\}.\qed $$
\enddemo

We now observe that there are versions of the
preceding results which are valid near $Y$.  To make
this precise, let $(X,Y)$ be a pair of spaces.  An
open neighborhood $U$ of $Y$ in $X$ is said to be a 
{\it teardrop neighborhood} if the pair $(U,Y)$ is a
teardrop; that is, there is a map
$$p : U \setminus Y \to Y\times \br $$
such that the identity from $X$ to $(X\setminus Y)\cup _pY$
is a homeomorphism.  The following results follow
immediately  from Corollaries 3.11 and 3.12.

\proclaim{Corollary 3.13}  If $Y$ is a compact subset
of the metric space $X$, then $Y$ has a teardrop
neighborhood in $X$ if and only if $Y$ is a
neighborhood retract of $X$. 
\qed
\endproclaim

\proclaim{Corollary 3.14} Let $Y$ be a closed subset
of the locally compact
metric space $X$.
Then $Y$ has a teardrop neighborhood in $X$ if and
only if $Y$ is a 
neighborhood retract of $X$.
\qed
\endproclaim

Next, we prove a lemma which will
be useful in \S4.

\proclaim{Lemma 3.15}
If $X$ and $Y$ are metric spaces
and $p: X\to Y\times \br$ is a map, then the
teardrop $X\cup _p Y$ is metrizable.  
\endproclaim

\demo{Proof}  Let $d_X $ and $d_Y$ be metrics for
$X$ and $Y$, respectively.  Define a function $\rho :
(X\amalg Y)\times (X\amalg Y) \to [0,+\infty )$ by 
$$\rho (a,b) = \cases
      d_X (a,b), & \text{if  $a,b\in X$}\\
      d_Y (a,b), & \text{if  $a,b\in Y$}\\
      0,         & \text{otherwise}.
\endcases$$
Define a metric $d$ on $Y\times (-\infty , + \infty ]$
by
$$d((y_1, t_1), (y_2, t_2)) = \max \{ d_Y (y_1, y_2),
     \ \vert e^{-t_1} -e^{-t_2}\vert \}$$
where $e^{-\infty } = 0$.  Note that $d$ generates the
standard topology.
Define the metric $D$ on $X\cup _p Y$ by 
$$D(a,b) = \rho (a,b) + d (c(a), c(b))$$
where $c: X\cup_p Y\to Y\times (-\infty , + \infty ]$
is the usual collapse.  One checks that $D$ generates
the teardrop topology.
\qed
\enddemo

\subhead Related constructions \endsubhead
Whyburn appears to be the first to have considered a
construction similar to the teardrop (see \cite{66},
\cite{67}), and we now explore the relationship between
the two constructions.  Let $X$ and $Y$ denote
(disjoint) Hausdorff spaces and let $f: X\to Y$ be a
map.  Whyburn defines the {\it unified space\/} $X\oplus
Y$ to be the topological space whose underlying set is
$X \amalg Y$  and with topology ${\cal W}$ given by
$V\in {\cal W}$ if and only if 
\widestnumber\item{(ii)}
\roster
\item"(i)" $V\cap X,\  V\cap Y$ are open in $X,Y$,
respectively, and 

\item"(ii)" for every compact $K\subseteq V\cap Y,\
     f^{-1}(K)\cap (X\setminus V)$ is compact.
\endroster
Whyburn proved that the function $r: X\oplus Y\to Y$
defined by $r(x) = f(x) $ if $x\in X$, and $r(y) = y$
if $y\in Y$, is a (continuous) proper retraction.

For the next two propositions suppose $p: X\to
Y\times [0, + \infty )$ is a map and $f: X\to Y$ is
the composition
$$X @>p>> Y\times [0,+\infty ) 
     @>\proj>> Y.$$

\proclaim{Proposition 3.16}
If $X$ and $Y$ are
Hausdorff, then the identity from $X\oplus Y$ to
$X\cup_p Y$ is continuous if and only if $p$ is proper.
\endproclaim

\demo{Proof}  Assume first that the identity is
continuous, and note that there is a commutative
diagram of maps:
$$\CD X\oplus Y @>r>> Y\\
@V{\id}VV  @AA{\proj}A\\
X\cup_p Y @>c>> Y\times [0+\infty ]
\endCD$$               
Recall that Whyburn showed $r$ is proper.  It follows
easily  that $c$ is proper.  Hence, $p = c|: X\to
Y\times [0, +\infty )$ is proper.

Conversely, assume that $p$ is proper.  By Lemma 3.4,
it suffices to show that the function $q: X\oplus Y
\to Y\times [0,+\infty ]$ defined by
$$ q(x) = \cases
            p(x), & \text{if  $x\in X$}\\
            (x,+\infty )  & \text{if  $x\in Y$}
\endcases$$
is continuous.  For this we need to show
$q^{-1}(U\times (t,+\infty ])$ is open where $U$ is
open in $Y$ and $t\in [0, +\infty )$.  Let $V = q^{-1}
(U\times (t, + \infty ]) = p^{-1} (U\times (t, +\infty
)) \cup U$ and let $K$ be a compact subset of $V\cap Y
= U$.  Then
$$f^{-1} (K)\cap (X\setminus V) = p^{-1} (K\times [0, +\infty ))
     \cap (X\setminus V) = p^{-1}(X\times [0,t])$$
which is compact since $p$ is proper.  Hence, $V$ is open.
\qed
\enddemo

\proclaim{Proposition 3.17} If $X$ and $Y$ are
Hausdorff, $Y$ is locally compact, and $p$  is proper,
then the identity from $X\cup_p Y$ to $X\oplus Y$ is continuous.
\endproclaim

\demo{Proof}  Let $U$ be open in $X\oplus Y$.  To
show $U$ is open in $X\cup_p Y$, it suffices to
consider $y\in U\cap Y$ and show $U$ is a neighborhood
of $y$ in $X\cup_p Y$.  Let $y\in V\subseteq \cl(V)
\subseteq U\cap Y$ where $V$ is open and $\cl(V)$ is
compact.  Then $c^{-1} (V\times [0, +\infty ])$ is
open in $X\cup _p Y$ and $f^{-1} (\cl(V))\cap (X\setminus U)$ is
a compact subset of $X$.  Let
$$W = c^{-1} (V\times [0, +\infty ]) \setminus [ f^{-1} (\cl(V))
     \cap (X\setminus U)].$$
Then $y\in W$, $W$ is open in $X\cup _pY$ and $W\subseteq U$.
\qed
\enddemo

\proclaim{Corollary 3.18} If $X$ and $Y$ are
Hausdorff and $Y$ is locally compact, then the
identity from $X\oplus Y$ to $X\cup _p Y$ is a
homeomorphism if and only if $p$ is proper. 
\qed
\endproclaim

Many authors  (\cite{14}, \cite{15}, \cite{36}, \cite{53})
have used a construction closely related
to Whyburn's unified space and we now briefly discuss
their construction.  Suppose $X,Y$ are disjoint,
Hausdorff, spaces, $X$ is locally compact and
non-compact, $Y$ is compact, $N$ is a neighborhood of
infinity in $X$, and $f: N\to Y$ is a map.  Then
Ferry and Pedersen \cite{15} define a space $X\amalg _f Y$ whose
underlying set is $X\amalg Y$ and with topology
generated by the basis
\block
     $\{ V ~|~ V$ is open in  $X\} \cup
     \{ (f^{-1}(V)\cap V')\cup V ~|~ V$ 
     is open in $Y$ and $V'$ is an open neighborhood
     of infinity in $X\}$.
\endblock
It is easy to see that the identity from $X\cup
(N\oplus Y)$ to $X\amalg_fY$ is a homeomorphism.
For an alternate treatment of related constructions,
one should consult James \cite{34\rm, \S8}.

\subhead Controlled maps \endsubhead
Finally, we use the teardrop topology to
clarify the notion of a controlled map given in
\cite{29\rm, \S12}.  For notation, if $\alpha$ is any
map we will let $\M(\alpha )$ denote the mapping
cylinder of $\alpha$ with the standard quotient
topology.  On the other hand, $\cyl(\alpha)$ will
denote the mapping cylinder with the teardrop
topology as in 3.1.
Suppose $f_t : X_1\to X_2 ,\ 0 \le t< 1$, is a
family of maps such that the induced map $f: X_1
\times [0,1)\to X_2$ is continuous.  Let $p: X_1 \to
Y$ and $q: X_2\to Y$ be given maps.

\proclaim{Proposition 3.19} The following are equivalent\rom:
\widestnumber\item{(iii)}
\roster
\item"(i)" $f_t$ is a controlled map from $p$ to $q$\rom; i\.e\.,
$\hat f :X_1 \times [0,1]\to Y$ given by
$$\hat f(x,t) = \cases qf_t(x), & \text{if $t< 1$}\\ 
          p(x), & \text{if $t = 1$}
\endcases$$
is continuous.

\item"(ii)" $f_\ast : X_1\times [0,1] \to \cyl (q)$
     given by
$$ f_\ast (x,t) = \cases  (f_t(x), t), & \text{if $t< 1$}\\ 
               p(x), & \text{if $t = 1$}
\endcases$$
is continuous.

\item"(iii)" $\tilde f: \M(p) \to\cyl (q) $ given by
$$ \cases
     \tilde f([x,t]) =  (f_t(x),t), & \text{if $t< 1$}\\ 
     \tilde f([y]) = y, & \text{if $y \in Y$}
\endcases$$
is continuous.  \endroster
\endproclaim

\demo{Proof}  (i) {\it implies\/} (ii): Since $\hat f$ is continuous,
     so is $cf_\ast : X_1\times [0,1]\to Y\times [0,1]$.
Lemma 3.4 then implies $f_\ast $ is continuous.

\noindent
(ii) {\it implies\/} (iii): Let $\pi : (X_1 \times [0,1])
\amalg Y\to \M(p) $ be the quotient map.  Then $\tilde
f$ is continuous if $\pi \tilde f$ is.  But $\pi
\tilde f|X_1\times [0,1] = f_\ast $ and $\pi
\tilde f|Y$ is the inclusion.

\noindent
(iii) {\it implies\/} (i): Note that $\hat f$ is the composition
$$X_1 \times [0,1] @>{\pi}>> \M(p)
@>{\tilde f}>> \cyl (q) @>c>> Y\times
[0,1] @>\proj>> Y.\qed$$
\enddemo

\head 4. The teardrop of an approximate fibration
\endhead

In this section we study the teardrop of an approximate 
fibration  $p: X\to Y\times
\br $ and establish two important properties.
First, if $X$ and $Y$ are metric spaces, then the
teardrop $(X\cup _p Y, Y)$ is a homotopically
stratified pair (Theorem 4.7).  Second,
if $p$ is a manifold approximate fibration, then $(X\cup_pY,Y)$ is
a manifold stratified pair (Corollary 4.11). This second result is
part of Theorem 2.1 and does not require the assumption that the 
dimension be greater than $4$.
The main technical tool is Theorem 4.2 which characterizes a
homotopically stratified pair in terms of a certain lifting property.
There are two other useful results.
One (Proposition 4.4) shows that the property of being a homotopically
stratified pair depends only on a neighborhood of the lower stratum.
The other (Proposition 4.8) characterizes (up to
fibre homotopy equivalence) the homotopy link as the the Hurewicz 
fibration associated to the induced map $X\to Y$.

We begin with the definition of the lifting property
which characterizes homotopically stratified pairs.
Let $(X,Y)$ be a pair such that $Y$ is a neighborhood
retract of $X$, and fix an open neighborhood and a
retraction $r : U\to Y$.  Consider the following spaces:
$$\align
W_1 (r) &= \{ (x,\omega)\in Y\times \map(I,Y) ~|~ x = \omega (1) \},\\
W_2 (r) &= \{(x,\omega)\in(U\setminus Y)\times\map(I,Y) ~|~ r(x) = 
\omega (1) \} , \quad \text{and}\\
W(r) &= W_1 (r) \cup W_2 (r) = \{ (x,\omega)\in U
     \times \map (I,Y) ~|~ r(x) = \omega(1)\} .
\endalign$$
Mapping spaces are always given the compact-open topology.
Note that the map $w(r):W(r)\to Y$ defined by $w(r)(x,\omega) = \omega(0)$
is the associated Hurewicz 
fibration of $r$,
and $w(r)|:W_2(r)\to Y$ is the associated Hurewicz fibration of $r|:U\setminus Y
\to Y$.

\definition{Definition 4.1} 
The pair $(X,Y)$ has the $W(r)$-{\it lifting property\/}
if there exists a map
$$\alpha : W(r)\to \map (I,X)$$
such that
\roster
\item $\alpha (x,\omega)(0) = \omega(0) $ for all $(x,\omega)\in W(r)$,
\item $\alpha (x,\omega)(1) = x $ for all $(x,\omega)\in W(r)$,
\item if $(x,\omega)\in W_1(r)$, then $\alpha (x,\omega)=\omega$, and
\item if $(x,\omega)\in W_2 (r)$, then $\alpha (x,\omega)\in \map
_s ((I,0) , (X,Y)) = \ho(X,Y)$.
\endroster
\enddefinition

\proclaim{Theorem 4.2}  If $X$ is a metric space and $Y
\subseteq X$, then the following are equivalent\rom:
\widestnumber\item{(iii)}
\roster
\item"(i)" $(X,Y)$ is homotopically stratified,

\item"(ii)" $Y$ is a neighborhood retract of $X$ and for every
sufficiently small neighborhood $U$ of $Y$ and retraction 
$r: U\to Y,\ (X,Y)$ has the $W(r)$-lifting property,

\item"(iii)" there exist a neighborhood $U$ of $Y$ and a
retraction $r: U\to Y$ such that $(X,Y)$ has the 
$W(r)$-lifting property. 
\endroster
\endproclaim

\demo{Proof}  (i) {\it implies\/} (ii): Since $(X,Y)$
is homotopically stratified, hence forward tame, there
exists a neighborhood $N$ of $Y$ and a nearly strict
deformation %
$$H: (N\times I, Y\times I\cup N\times \{ 0\})\to (X,Y).$$
In particular, $Y$ is a neighborhood retract of $X$.
Let $U$ be any neighborhood of $Y$ such that $U\subseteq
N$ and let $r: U\to Y$ be any retraction.  We will
show that $(X,Y)$ has the $W(r)$-lifting property.
Define a map $\beta : W(r) \to \map (I,Y)$ by the formula
$$\beta (x,\omega) (t) =
           \cases  rH(x,2t), & \text{if $0 \le t\le {1\over 2}$}\\
           \omega(2-2t), & \text{if ${1\over 2} \le t\le 1$}.
\endcases$$
Define $f: W_2 (r) \to \ho(X,Y)$ by $f(x,\omega)(t)
     = H(x,t) $ for $t\in I$, and define
$$F: W_2 (r) \times I\to Y$$
by $F(x,\omega,t) = \beta (x,\omega) (t) $ for $t\in I$.
Note that we have a lifting problem:
$$\CD W_2(r) @>f>> \ho(X,Y)\\
      @V{\times 0}VV @VVqV\\
      W_2 (r)\times I  @>F>> Y.
\endCD$$
(Recall that $q$ is evaluation at 0).  Since part of our
hypothesis is that $q$ is a fibration, we have a
solution $\tilde F$.  We will use $\tilde F$ to define
$\alpha $, but to make sure that a certain extension
to $W(r)$ is continuous on $W_1 (r)$, we first need a lemma whose proof is
postponed until later in this section.

\proclaim{Lemma 4.3}      There exists a map $\gamma : W_2 (r) \times I
\to [0,1]$ such that
\roster
\item $\gamma (x,\omega,0) = 1 $ for all $(x,\omega) \in W_2 (r)$,
\item $\diam  \{ \tilde F(x,\omega,t)(s) ~|~ 0\le s\le \gamma
          (x,\omega,t)\} \le 2 \diam \{ \tilde F(x,\omega,0)
     (s) ~|~ s\in I\}$ for all $(x,\omega,t)\in W_2 (r)
     \times I$,
\item $\gamma (x,\omega,t) = 0$  if and only if $t = 1$,
     for all $(x,\omega) \in W_2 (r)$. 
\endroster
\endproclaim

Assuming the lemma 
we complete the proof that (i) {\it implies\/}
(ii) in Theorem 4.2.  Define
$$\alpha : W_2 (r) \to \ho(X,Y)
~~\text{by}~~
\alpha (x,\omega) (t) = \tilde F (x,\omega,1-t)
     (\gamma (x,\omega,1-t)). $$
Then $\alpha $ extends to a map $\alpha : W(r)\to
     \map (I,X)$ by setting
$\alpha (x,\omega) = \omega$
for $(x,\omega)\in W_1(r)$.  It is straightforward to
verify that $\alpha $ is continuous and satisfies the
condition of the 
$W(r)$-lifting property.

\noindent
(ii) {\it implies\/} (iii) is obvious.

\noindent
(iii) {\it implies\/} (i):
Let $\alpha : W(r)\to \map (I,X)$ satisfy the
definition of the $W(r)$-lifting property where $r:
U\to Y$ is some retraction of a neighborhood of $Y$.
For each $x\in U$ let $\omega _x$ denote the constant path
at $r(x)$.  Define $H: U\times I \to X$ by 
$$H(x,t) = \alpha (x,\omega_x)(t).$$
Then $H$ is a nearly strict deformation of $U$ into
$Y$, so $Y$ is forward tame in $X$.
To see that $q:\ho(X,Y)\to Y$ is a fibration,
consider a lifting problem
$$\CD Z @>f>> \ho(X,Y)\\
     @V{\times 0}VV @VqVV\\
     Z\times I @>F>> Y.
\endCD$$
We may assume that $Z$ is metric.  Using a partition of
unity one can construct a map
$\epsilon : Z\to (0,1]$
such that for every $z\in Z$ and $0\le t\le \epsilon (z)$,
we have $f(z)(t)\in U$.
Define a map $\omega : Z\times I\to \map (I,Y)$ by
$$\omega (z,t)(s) =\cases
      F(z,t-2ts), & \text{if $0\le s\le 1/2$}\\
      r(f(z)(\epsilon (z)(2ts-t))), & \text{if $1/2 \le s\le 1$}.
\endcases$$
Note that $\omega (z,0)(s) = F(z,0) = f(z)(0)$ for all $z\in
Z$ and $s\in I$.
Now define $$\delta : Z\times I\to \map (I,X) ~~\text{by}~~
\delta (z,t)   =
     \alpha (f(z)(\epsilon (z)t),
          \omega (z,t))$$
and note that
\roster
\item $\delta (z,0)(s) = F(z,0),$
\item $\delta (z,t) (1) = f(z)(\epsilon(z)t),$
\item $\delta (z,t)(0) = F(z,t).$
\endroster
Finally, define a solution $\tilde F: Z\times I\to \ho
(X,Y)$ of the lifting problem by
$$\tilde F(z,t) (s) =\cases
     \delta (z,t)(s/\epsilon (z)t), &\text{if $0 \le s < \epsilon (z)t$}\\
      f(z)(s), &\text{if $\epsilon (z)t\le s\le 1$}.\qed
\endcases$$
\enddemo

\demo{Proof of Lemma \rom{4.3}}  First note that $\{ \tilde F(x,
\omega,0)(s)\mid s\in I\} = \{ H(x,s)\mid s\in I\}$ for each
$(x,\omega)\in W_2 (r)$.  Now for $x\in U\setminus Y$, let $c(x) = \diam
\{ H(x,s) ~|~ s\in I\}$.  Note that $0< c(x)$.
For each $(x,\omega,t)\in W_2 (r)\times I$, let 
$$\delta (x,\omega,t) 
=  \lub \{ s\in I\vert  \diam \{ \tilde F(x,\omega,t)(s')
\vert 0\le s'\le s\} \le c (x)\}.$$  
Note that $0<\delta (x,\omega,t) \le 1$.
For each $(x,\omega,t)\in W_2 (r)\times I$, let $V(x,\omega,t)$ be
a neighborhood of $(x,\omega,t)$ such that whenever $(x',\omega',t')\in 
V(x,\omega,t)$, then 
\roster
\item $\diam \{ \tilde F(x', \omega', t')(s)
     ~|~ 0 \le s\le \delta (x,\omega,t)\} < 3c(x)/2$, and
\item $c(x) \le 4 c(x')/3$.
\endroster
Let $\{ V_\alpha\}$ be a locally finite refinement of
$\{ V(x,\omega,t)\}$ and let $\{ \phi _\alpha \}$ be a
partition of unity subordinate to $\{ V_\alpha \}$.
For each $\alpha $ choose $(x,\omega,t)$ such that
$V_\alpha \subseteq V(x,\omega,t)$ and set $\delta _\alpha =
\delta (x,\omega,t)$.  Define $\hat\gamma :W_2 (r)\times
I\to I $ by $\hat\gamma = \sum \delta _\alpha \phi
_\alpha $.  Clearly $\hat\gamma$ satisfies
item (2) of the lemma, but we need to modify $\hat
\gamma$ to achieve the other conditions. 
Using the paracompactness of $W_2 (r)$, choose a
neighborhood $V$ of $W_2 (r)\times \{ 0\}$ in $W_2
(r)\times I$ such that if $(x,\omega,t)\in V$, then
$$\diam \{ \tilde F (x,\omega,t)(s) ~|~ s\in I\} \le 2c(x).$$
Let $\psi : W_2 (r)\times I\to I$ be a map such that
$\psi = 1$ on $W_2 (r)\times \{ 0\},\  \psi = 0 $ off
of $V$, and $\psi > 0$ on $V$.  Finally set
$$\gamma (x,\omega,t) = (1-t) [(1-\psi (x,\omega,t))\hat\gamma
     (x,\omega,t) + \psi (x,\omega,t)]. \qed$$
\enddemo

\proclaim{Proposition 4.4}
If $X$ is a metric space and $Y\subseteq X$, then the following are 
equivalent\rom:
\widestnumber\item{(iii)}
\roster
\item"(i)" $(X,Y)$ is a homotopically stratified pair,
\item"(ii)" for every neighborhood $U$ of $Y$ in $X$, $(U,Y)$ is a homotopically
stratified pair,
\item"(iii)" there exists a neighborhood $U$ of $Y$ in $X$ such that
$(U,Y)$ is a homotopically stratified pair.
\endroster
\endproclaim

\demo{Proof} (i) {\it implies\/} (ii): Let $U$ be a neighborhood of $Y$ in
$X$. Forward tameness implies there exist a neighborhood $N$ of $Y$
in $X$ such that $N\subseteq U$ and a nearly strict deformation of
$N$ to $Y$ in $U$ which gives a retraction $r:N\to Y$.
The proof of Theorem 4.2 (i) implies (ii) shows that if $N$ is a sufficiently
small neighborhood of $Y$ in $U$, then $(U,Y)$ has the $W(r)$--lifting
property so that Theorem 4.2 may be invoked.

\noindent
(ii) {\it implies\/} (iii) is obvious.

\noindent
(iii) {\it implies\/} (i): By Theorem 4.2 we know that $(U,Y)$ has the
$W(r)$--lifting property for some $r$. It follows that $(X,Y)$ has the 
$W(r)$--lifting property and Theorem 4.2 may be invoked once again.
\qed
\enddemo

We now recall the definition of approximate fibrations as given in
\cite{29}. See \cite{29\rm, \S12} for an explanation of how this definition
relates to others in the literature.

\definition{Definition 4.5} 
A map $p:E\to B$ is an {\it approximate fibration\/} if for every
commuting diagram
$$\CD Z @>f>> E\\ 
      @V{\times 0}VV @VVpV\\
      Z\times [0,1] @>F>> B
\endCD$$
there is a controlled map  $\tilde F:Z\times [0,1]\times [0,1)\to E$
from $F$ to $p$ such that $\tilde F(x,0,u)=f(x)$ for all $(x,u)\in
Z\times [0,1)$.
To say $\tilde F$ is a {\it controlled map\/} from $F$ to $p$ means
the function
$G: Z\times [0,1]\times [0,1]
     \to B$ defined by
$$G(z,t,u) = \cases  p\tilde F(z,t,u), & \text{if $u< 1$}\\
                  F(z,t), & \text{if $u=1$}
\endcases$$
is continuous.
\enddefinition

\proclaim{Lemma 4.6  (Open ended homotopies)}
Suppose that $p : E\to B$ is an approximate fibration
and that the following lifting problem is given\rom:
$$\CD  Z  @>f>> E\\
       @V{\times 0}VV @VVpV\\
         Z\times [0,1)  @>F>> B.
\endCD$$
Then there exists a controlled lift $\tilde F$, i\.e\.,
a map $\tilde F: Z\times [0,1)\times [0,1)\to E$
such that
\widestnumber\item{(ii)}
\roster
\item"(i)" $\tilde F(z,0,u) = f(z)$ for all $u\in [0,1)$,  and

\item"(ii)" the function $G: Z\times [0,1)\times [0,1]
     \to B$ defined by
$$G(z,t,u) = \cases p\tilde F(z,t,u), & \text{if $ u < 1$}\\
                F(z,t), & \text{if $u=1$}
\endcases$$
is continuous.
\endroster\endproclaim

\demo{Proof}
Let $\pi:{\cal E}\to B$ be the Hurewicz fibration associated to
$p:E\to B$ and let $i:E\to{\cal E}$ be the inclusion.
According to \cite{29\rm, 12.5} there is a controlled map
$R:{\cal E}\times [0,1)\to E$ from $\pi$ to $p$ and a controlled homotopy
$H:E\times [0,1]\times [0,1)\to E$ from $\id_e$ to $Ri$. This means that
the function
$\ov{R}:{\cal E}\times [0,1]\to B$ defined by
$$\ov{R}(x,t) = \cases pR(x,t), & \text{if $t<1$}\\
                             \pi(x),  & \text{if $t=1$}
\endcases$$
is continuous, that $H$ satisfies
$H(x,0,t) = x$ and $H(x,1,t) = R(i(x),t)$ for all $(X,t)\in E\times [0,1)$,
and that the function $\ov{H}:E\times [0,1]\times [0,1]\to B$
defined by
$$\ov{H}(x,s,t) = \cases pH(x,s,t), & \text{if $t<1$}\\
                             p(x),  & \text{if $t=1$}
\endcases$$
is continuous.
Given a lifting problem of the form
$$\CD  Z  @>f>> E\\
       @V{\times 0}VV @VVpV\\
         Z\times [0,1)  @>F>> B.
\endCD$$
there is an induced problem
$$\CD  Z  @>if>> {\cal E}\\
       @V{\times 0}VV @VV{\pi}V\\
         Z\times [0,1)  @>F>> B.
\endCD$$
Since $\pi$ is a fibration, this second problem has an exact solution
$\hat F:Z\times [0,1)\to {\cal E}$.
Define $F':Z\times [0,1)\times [0,1)\to E$ by
$F'(z,s,t) = R(\hat F(z,s),t)$.
Then a controlled solution $\widetilde{F}:Z\times [0,1)\times [0,1)\to E$
to the first problem can be defined by
$$\widetilde{F}(z,s,t) = \cases H(f(z),{s\over 1-t},t), & 
\text{if $0\leq s\leq 1-t$}\\
           F'(z,{s-1-t\over t},t), & \text{if $1-t\leq s\leq 1$.}
\endcases$$
One checks that the function $G$ defined in the statement is continuous.
\qed
\enddemo

\proclaim{Theorem 4.7} If $X$ and $Y$ are metric
spaces and $p: X\to Y\times \br $ is an
approximate fibration, then the teardrop $(X\cup_p Y,
Y)$ is a homotopically stratified pair. 
\endproclaim

\demo{Proof}  There exists a retraction $r:
X\cup_p Y\to Y$ given by the composition
$$X\cup_p Y @>c>> Y\times (-\infty , +\infty ]
     @>{\proj}>> Y.$$
Since $X\cup_p Y$ is metric by Lemma 3.15, it suffices
by Theorem 4.2 to show that $(X\cup_p Y, Y)$ has the $W(r)$--lifting property.
We will first define $\alpha $ on $W_2(r)$ and then
extend it to all of $W(r)$.  To this end define
$$F: W_2(r)\times [0,1) \to Y\times \br 
\qquad\text{by}\qquad
F(x,\omega ,t) = (\omega (1-t), {s\over 1-t} )
$$
where $s$ is defined by $p(x) = (r(x), s) \in Y\times
\br $.  Define $f: W_2 (r)\to X$ by $f(x,\omega ) =
x$.  Then we have a lifting problem
$$\CD W_2 (r) @>f>> X\\
        @V{\times 0}VV  @VVpV\\
      W_2 (r)\times [0,1) @>F>> Y\times \br
\endCD$$
to which we can apply Lemma 4.6 and get a controlled lift
$$\tilde F: W_2 (r)\times [0,1)\times [0,1)\to X.$$
Let $G: W_2 (r)\times [0,1)\times [0,1]\to Y\times
\br $ be the map defined in Lemma 4.6.  Using the
paracompactness of $W_2 (r)\times [0,1)$, there exists
a map $\gamma : W_2 (r)\times [0,1)\to [0,1)$ such
that if $(x,\omega) \in W_2 (r)$ and $1- {1/i}\le t\le 1-
1/(i+1)$, then
$$\diam G(\{ x,\omega ,t\} \times [\gamma (x,\omega ,t),1])
     < {1/i}.$$
Then define $\hat F :W_2 (r)\times [0,1]\to X
\cup_p Y$ by
$$\hat F(x,\omega ,t)=
     \cases \tilde F(x,\omega ,\gamma (x,\omega ,t)),
               & \text{if  $0\le t< 1$}\\
          \omega (0),
               & \text{if $t = 1$}.
\endcases$$
And define $\alpha : W_2 (r)\to\ho(X\cup_pY,
Y)$ by 
$$\alpha (x,\omega )(t) = \hat F(x,\omega ,1-t).$$
Then $\alpha $ extends continuously to $\alpha : W(r)
\to\map(I,X\cup_p Y)$ by setting $\alpha (x,\omega ) = \omega $
for $(x,\omega )\in W_1 (r)$.
\qed
\enddemo

\proclaim{Proposition 4.8}
If $X$ and $Y$ are metric spaces and $p:X\to Y\times\br$ is an approximate
fibration, then $q:\ho(X\cup_pY,Y)\to Y$ is fibre homotopy equivalent
to the Hurewicz fibration associated to the composition
$$X @>p>> Y\times\br @>{\proj}>> Y.$$
\endproclaim

\demo{Proof} Let $r:X\cup_pY\to Y$ be the retraction
$X\cup_p Y @>c>> Y\times (-\infty , +\infty ] @>{\proj}>> Y$.
Let $\pi = w(r)|:W_2(r)\to Y$ which is the Hurewicz fibration associated to
$r|:X\to Y$. We must show that $\pi$ is fibre homotopy equivalent to
$q:\ho(X\cup_pY,Y)\to Y$.
It follows from the proof of Theorem 4.7 that 
$(X\cup_pY,Y)$ has the $W(r)$--lifting property.
Let $\alpha:W(r)\to\map(I,X\cup_pY)$ be a map as in Definition 4.1.
Define $f:W_2(r)\to\ho(X\cup_pY,Y)$  to be the restriction of $\alpha$
so that $f(x,\omega) = \alpha(x,\omega)$. We will show that $f$ is a
fibre homotopy equivalence with fibre homotopy inverse
$g:\ho(X\cup_pY,Y)\to W_2(r)$ defined by $g(\omega) = (\omega(1),r\omega)$.
We will define a fibre homotopy $G:gf\simeq\id_{W_2(r)}$ as follows.
If $\omega\in\map(I,Y)$ and $s\in I$, define
$\omega_s^+:I\to Y$ by
$\omega_s^+(t) = \omega((1-s)t + s)$. Define a homotopy
$E:W_2(r)\times I\to \map(I,Y)$ by
$$E(x,\omega,s)(t) =\cases \omega(t), &  \text{if $0\leq t\leq s$}\\
       r\alpha(x,\omega_s^+)({t-s\over 1-s}), &  \text{if $s\leq t < 1$}\\
       r(x), & \text{if $t= 1$.}
\endcases$$
Then let $G((x,\omega),s) = (x, E(x,\omega,s))$.
We will now define a fibre homotopy $F:\id_{{\ho}(X\cup_pY,Y)}\simeq fg$ as
follows. If $\omega\in\ho(X\cup_pY,Y)$ and $s\in I$, define
$\omega_s:I\to X\cup_pY$ by
$\omega_s(t) = \omega(ts)$. Then define $F$ by
$$F(\omega,s)(t) = \cases \omega(0), & \text{if $t=0$}\\
                   \alpha(\omega(s),r\omega_s)({t\over s}), & \text{if 
                               $0<t\leq s$}\\
                   \omega(t), &  \text{if $s\leq t\leq 1$.}
\qed\endcases$$
\enddemo

\proclaim{Lemma 4.9 (Folklore)} If $p:X\to Y$ is a proper
approximate fibration between ANRs \rom(locally compact, separable metric\rom),
then the homotopy fibre of $p$ is finitely dominated.
\endproclaim

\demo{Proof} 
Fix a basepoint $y_0\in Y$. The homotopy fibre of $p$ is
$$W=\{ (x,\omega)\in X\times Y^I ~|~ \omega(0)=p(x), \omega(1)=y_0\}.$$
Let $U$ be an open neighborhood of $y_0$ which contracts to $y_0$ in $Y$;
that is, there exists a homotopy $H:U\times I\to Y$ such that $H_0=
\inclusion:U\to Y$, $H_1(U) =\{ y_0\}$ and $H_t(y_0)=y_0$ for all
$t\in I$.
Let $V$ be a compact neighborhood of $y_0$ such that $H(V\times I)\subseteq U$.
It is well-known that for every open cover $\cal U$ of $X$ there is a locally
finite simplicial complex which $\cal U$-dominates $X$ (see e.g.
\cite{42}). This fact together with the compactness of
$p^{-1}(V)$ implies that there exist a locally finite simplicial complex
$L$, maps $f:L\to X$, $g:X\to L$, and a homotopy
$J:\id_X\simeq fg$ such that $J(p^{-1}(V)\times I)\subseteq p^{-1}(U)$.
Note that $g(p^{-1}(V))\subseteq f^{-1}(p^{-1}(U))$ and use the compactness of
$p^{-1}(V)$ again to find a finite subcomplex $K$ of $L$ 
(in some fine triangulation)
such that $g(p^{-1}(V))\subseteq K$ and $ f(K)\subseteq p^{-1}(U)$.
We will show that $K$ dominates $W$.
Consider the lifting problem
$$\CD
W @>{g}>> X\\
@V{\times 0}VV @VV{p}V\\
W\times I @>{G}>> Y\\
\endCD$$
where $G((x,\omega),t)=\omega(t)$ and $g(x,\omega)=x$.
Since $p$ is an approximate fibration there is an approximate solution
$\widetilde{G}:W\times I\to X$.
Assume that $p\widetilde G$ is so close to $G$ that the image of 
$\widetilde{G}_1$ is in $p^{-1}(V)$ and that there is a homotopy
$F:p\widetilde{G}\simeq G$ rel $W\times\{0\}$. Using the homotopy extension
theorem we can insist that $F|W\times\{ 1\}\times I$ is given by
$F((x,\omega), 1, s) = H(p\widetilde{G}_1(x,\omega),s)$.
It follows that there is a homotopy $A:W\times I\times I\to Y$ such that
\roster
\item $A((x,\omega),0,s) = \omega(0)$,
\item $A((x,\omega)1,s)=H(pfg\widetilde{G}_1(x,\omega),s)$,
\item $A((x,\omega),t,1)=\omega(t)$,
\item $A((x,\omega)t,0)=
\cases p\widetilde{G}((x,\omega),2t), &0\leq t\leq {1\over 2}\\
pJ(\widetilde{G}_1(x,\omega),2t-1), &{1\over 2}\leq t\leq 1.\\
\endcases$
\endroster
Define $d:K\to W$ and $u:W\to K$ by
$d(x) = (f(x), H(pf(x),\cdot))$ and 
$u(x,\omega) = g(\widetilde{G}_1(x,\omega))$.
The homotopy $A$ can be used to construct a homotopy $du\simeq\id_W$.
\qed
\enddemo

\proclaim{Corollary 4.10} If $X$ and $Y$ are ANRs 
\rom(locally compact, separable
metric\rom) and $p:X\to Y\times\br$ is a proper approximate fibration, then
the \rom(homotopy\rom) fibre of $q:\ho(X\cup_pY,Y)\to Y$ is finitely dominated.
Moreover, $(X\cup_pY,Y)$ is a homotopically stratified locally compact,
separable metric pair with
finitely dominated local holinks. 
\endproclaim

\demo{Proof} It follows from Lemma 3.15 that $X\cup_pY$ is metrizable. Since
$X$ and $Y$ are separable, so is $X\cup_pY$. Since $p$ is proper, it 
follows easily that the teardrop collapse $c:X\cup_pY\to Y\times 
(-\infty,+\infty]$ is also proper (cf. the proof of Proposition 3.16).
In particular, $X\cup_pY$ is locally compact. By Theorem 4.7,
$(X\cup_pY,Y)$ is homotopically stratified. It follows from 
Proposition 4.8 that the homotopy fibre of 
$\ho(X\cup_pY,Y)\to Y$ is homotopy equivalent to the homotopy fibre of
$p$ which is finitely dominated by Lemma 4.9.
Thus, $(X\cup_pY,Y)$ has finitely dominated local holinks.
\qed
\enddemo

\proclaim{Corollary 4.11} If $B$ is a closed manifold and
$p:M\to B\times\br$ is a manifold approximate fibration, then the
teardrop $(M\cup_pB,B)$ is a manifold stratified pair.
\endproclaim

\demo{Proof} This follows immediately from Corollary 4.10.
\qed
\enddemo

\head 5. Spaces of stratified neighborhoods and 
manifold approximate fibrations
\endhead

This section contains the details of the definitions of
the simplicial set $\maf^n (B\times \br ) $ of manifold approximate fibrations
and the simplicial set 
$\SN^n(B)$ of stratified neighborhoods. Facts are established which are
needed to
define the simplicial map $\Psi:\maf^n(B\times\br )\to\SN^n(B)$.

\definition{Definition 5.1}
Suppose $A\times K$ is a closed subset of $X$ and $\pi: X\to K$ is a map
such that
$\pi|:A\times K\to K$ is the projection.

\noindent
(1) The pair $(X,A\times K)$ is a {\it sliced homotopically
stratified pair\/ \rom(with respect to\/ $\pi$\rom)} if
\widestnumber\item{(iii)}
\roster
\item"(i)" $A\times K$ is sliced forward tame in $X$ with respect
to $\pi$,
\item"(ii)" the evaluation $q:\ho_\pi(X,A\times K)\to A\times K$
is a fibration,
\item"(iii)" (Local triviality near $A\times K$) there exist an open
neighborhood $W$ of $A\times K$ in $X$ and a space $U$ containing
$A$ such that for each $t\in K$ there exist an open neighborhood
$V$ of $t$ in $K$ and 
a f\.p\. open embedding $h:U\times V\to X$ such that
$h|:A\times V\to X$ is the inclusion 
and $h(U\times V) = W\cap \pi^{-1}(V)$. That is, $\pi|:W\to K$ is
a fibre bundle projection containing $A\times K\to K$ as a subbundle.
In this case $W$ is said
to be a {\it locally trivial neighborhood of\/} $A\times K$ {\it in\/} $X$.
If $V=K$, then $W$ is said 
to be a {\it trivial neighborhood of\/} $A\times K$ {\it in\/} $X$.
\endroster
(2) The pair $(X,A\times K)$ has {\it finitely dominated 
local holinks \rom(with respect to\/ $\pi$\rom)}
if the fibre of $q:\ho_\pi(X,A\times K)\to A\times K$
is finitely dominated.

\noindent
(3) The pair $(X,A\times K)$ is a {\it sliced manifold 
stratified pair \rom(with respect to\/ $\pi$\rom)} if it is a sliced homotopically stratified pair
with finitely dominated local holinks, $X$ is a locally compact
separable metric space, $A$ is a manifold, and for each $t\in K$
$\pi^{-1}(t)\setminus A\times\{ t\}$ is a manifold.
\enddefinition

Note that if $K$ is contractible, then the local triviality condition
near $A\times K$ implies that $A\times K$ has a trivial neighborhood in
$X$.

\proclaim{Proposition 5.2}
Suppose $A\times K$ is a closed subset of a metric space
$X$ and $\pi: X\to K$ is a map such that
$\pi|:A\times K\to K$ is the projection.

\noindent
{\rm (i)} If $N$ is a neighborhood of $A\times K$ in $X$, then
the inclusion $\ho_\pi(N,A\times K) \to
\ho_\pi(X,A\times K)$ is a fibre homotopy equivalence from
$q:\ho_\pi(N,A\times K)\to A\times K$ to
$q:\ho_\pi(X,A\times K)\to A\times K$.

\noindent
{\rm (ii)} If $N$ is a neighborhood of $A\times K$ in $X$, then
$q:\ho_\pi(X,A\times K) \to A\times K$ is a fibration if and only if
$q:\ho_\pi(N,A\times K)\to A\times K$ is.

\noindent
{\rm (iii)} If $K$ is compact, the following are equivalent\rom:
\roster
\item"(a)" $(X,A\times K)$ is a sliced homotopically stratified pair,
\item"(b)" for every neighborhood $N$ of $A\times K$ in $X$,
$(N,A\times K)$ is a homotopically stratified pair,
\item"(c)" there exists a neighborhood $N$ of $A\times K$ in $X$
such that $(N,A\times K)$ is a homotopically stratified pair.
\endroster
{\rm (iv)} If $N$ is a neighborhood of $A\times K$ in $X$, then
$(X,A\times K)$ has finitely dominated local holinks
if and only if $(N,A\times K)$ does.

\noindent
{\rm (v)} If $K$ is compact and $N$ is open an 
open neighborhood of $A\times K$ 
in $X$ and $(X,A\times K)$ is a sliced manifold
stratified pair, then so is  $(N,A\times K)$.
\endproclaim

\demo{Proof}(i)
(cf\. \cite{28\rm, 1.12}) For each $\omega\in\ho_\pi(X,A\times K)$ choose a
number $t_\omega\in (0,1]$ such that 
$\omega([0,t_\omega])\subseteq\inr(N)$.
Let $U(\omega)$ be an open neighborhood of $\omega$ in
$\ho_\pi(X,A\times K)$ such that
$\alpha([0,t_\omega])\subseteq{\inr}(N)$ for all $\alpha\in U(\omega)$.
Since $\ho_\pi(X,A\times K)$ is a metric space, there is a locally finite
refinement $\{U_i\}$ for the cover
$\{U(\omega)~|~\omega\in\ho_\pi(X,A\times K)\}$ of
$\ho_\pi(X,A\times K)$ and a partition of unity $\{\phi_i\}$ subordinate
to $\{U_i\}$.
For each $i$ choose
$\omega_i\in\ho_\pi(X,A\times K)$ such that $U_i\subseteq U(\omega_i)$
and let $t_i=t_{\omega_i}$.
For each $\omega\in\ho_\pi(X,A\times K)$ let 
$m_\omega =\max\{ t_i ~|~ \phi_i(\omega)\not= 0\}$.
Note that $\omega([0,m_\omega])\subseteq{\inr}(N)$ and
$\sum_i \phi_i(\omega)t_i\leq m_\omega$ for all $\omega$.
Define a homotopy
$R:\ho_\pi(X,A\times K)\times I\to \ho_\pi(X,A\times K)$ by
$$R(\omega,t)(s) = \cases
 \omega(s), & \text{if $0\leq s\leq \sum_i \phi_i(\omega)t_i$}\\
       \omega((1-t)s + t\sum_i\phi_i(\omega)t_i), & \text{if 
$\sum_i(\omega)t_i\leq   s\leq 1$.}
\endcases$$
Then $R$ is a fibre deformation with $R_0=\id$, 
$$R_1(\ho_\pi(X,A\times K))
\subseteq\ho_\pi(N,A\times K)$$ 
and
$R_t(\ho_\pi(N,A\times K))\subseteq\ho_\pi(N,A\times K)$ for each $t$.
The result follows immediately. Note also that if $\rho:
\ho_\pi(X,A\times K)\to (0,1]$ is defined by $\rho(\omega)=
\sum_i\phi_i(\omega)t_i$, then $\rho$ is continuous and $R_t(\omega)(s) =
\omega(s)$ for all $0\leq t\leq 1$ and $0\leq s\leq\rho(\omega)$.

\noindent
(ii) Let $R$ and $\rho$ be given as in the proof of (i).
Suppose first that  $q:\ho(N,A\times K)\to A\times K$ 
is a fibration. Then a  homotopy lifting problem
$$\CD Z @>f>> \ho_\pi(X,A\times K)\\
     @VVV  @VVqV\\
      Z\times I @>F>> A\times K
\endCD$$
for $\ho_\pi(X, A\times K)\to A\times K$ induces a problem
$$\CD  Z @>{R_1f}>>  \ho_\pi(N,A\times K)\\
     @VVV  @VVqV\\
      Z\times I @>F>> A\times K
\endCD$$
for $\ho_\pi(N,A\times K)\to A\times K$ which has a solution
$G:Z\times I\to\ho_\pi(N,A\times K)$.
For each $\omega\in\ho_\pi(X,A\times K)$ define 
$$\tau_\omega:[0,\rho(\omega)]\times [0,1]\to [0,1]\times[0,\rho(\omega)]
~~\text{ by }~~ 
\tau_\omega(s,t) = (t - {ts\over\rho(\omega)}, s).$$
Then a solution 
$\widetilde{F}:Z\times I\to \ho_\pi(X,A\times K)$
of the original problem can be defined by
$$\widetilde{F}(z,t)(s) = \cases \widehat{G}(z,\tau_{f(z)}(s,t)),
                 & \text{if $0\leq s\leq\rho(f(z))$}\\
               f(z)(s), & \text{if $\rho(f(z))\leq s\leq 1$}
\endcases$$
where $\widehat{G}$ is the adjoint of $G$.

\noindent
Conversely, suppose 
$q:\ho(X,A\times K)\to A\times K$ is a fibration and $N$ is a neighborhood
of $A\times K$ in $X$. To show that 
$\ho_\pi(N,A\times K)\to A\times K$ is a fibration, we may use the 
converse just proven to assume that $N$ is open in $X$.
Let 
$$\CD Z @>f>> \ho_\pi(N,A\times K)\\
      @VVV @VVqV\\
      Z\times I @>F>> A\times K
\endCD$$
be a homotopy lifting problem which
by inclusion is also a problem for 
$$\ho_\pi(X,A\times K)\to A\times K.$$
Thus, there is a solution $G:Z\times I\to\ho_\pi(X,A\times K)$
to this second problem.
Let $U$ be an open neighborhood of $Z\times\{0\}$ in $Z\times I$ such that
$G(U)\subseteq\ho_\pi(N,A\times K)$. Since it suffices to solve
an universal problem, we may assume that $Z$ is a metric space. Thus, there
is a map $\sigma:Z\times I\to I$ such that
$\sigma^{-1}(0) = Z\times \{0\}$ and $\sigma^{-1}(1) = (Z\times I)\setminus
U$. Then $\widetilde{F}:Z\times I\to \ho_\pi(N,A\times K)$ defined by
$\widetilde{F}(z,t) = R(G(z,t),\sigma(z,t))$ is a solution of the 
original problem.

\noindent
(iii) (a) {\it implies\/} (b): If $N$ is a neighborhood of $A\times K$ in
$X$, then $(N,A\times K)$ obviously satisfies the sliced forward tameness
condition. From the fact that $K$ is compact. it follows that 
$(N,A\times K)$ satisfies local triviality near $A\times K$. The holink
fibration condition follows from (ii).

\noindent
(b) {\it implies\/} (c) is obvious.

\noindent
(c) {\it implies\/} (a): 
 The sliced forward tameness and  local triviality
conditions obviously hold for $(X,A\times K)$ if they hold
for $(N,A\times K)$. 
The holink fibration condition follows from (ii).

\noindent
(iv) follows directly from (i).

\noindent
(v) follows (iii) and (iv).
\qed
\enddemo

\proclaim{Lemma 5.3}
Suppose $A\times K$ is a closed subset of a space
$X$ and $\pi: X\to K$ is a map such that
$\pi|:A\times K\to K$ is the projection.
Let $f:K'\to K$ be a map and form the pull--back diagram
$$\CD (X',A\times K') @>>> (X,A\times K)\\
      @V\pi'VV  @VV{\pi}V\\
      K' @>f>>  K.
\endCD$$
\widestnumber\item{(iii)}
\roster
\item"(i)" There is an induced pullback diagram
$$\CD  \ho_{\pi'}(X',A\times K') @>>> \ho_\pi(X,A\times K)\\
      @V{q'}VV  @VVqV\\
      A\times K' @>{\id_A\times f}>> A\times K
\endCD$$
\item"(ii)" If $(X,A\times K)$ is a sliced homotopically stratified pair,
then so is  $(X',A\times K')$.
\item"(iii)" If $(X,A\times K)$ has finitely dominated local holinks,
then so does $(X',A\times K)$.
\item"(iv)" If $(X,A\times K)$ is a sliced manifold
stratified pair, then so is  $(X',A\times K')$.
\endroster
\endproclaim

\demo{Proof}(i) and (ii) are elementary. The other parts follow
immediately.
\qed
\enddemo

For the remainder of this section, $B$ is an $i$-dimensional manifold
without boundary together with a fixed embedding
$B\subseteq \ell_2$ (of small capacity; e\.g\., we could
take $B$ to be inside of a finite dimensional subspace
$\br^L$ of $\ell _2$) and let $n\ge 5$ be a fixed integer.

\definition{Definition 5.4}
The {\it space of stratified neighborhoods\/} of $B$ is
the simplicial set $\SN^n(B)$ whose $k$-simplices are
subsets $X$ of $\ell_2 \times \Delta ^k$ of small
capacity (see \cite{29}) such that 
if $\pi : X\to \Delta ^k$ is the
restriction of the projection $\ell_2 \times \Delta ^k
\to \Delta ^k$,
then $(X,B\times\Delta^k)$ is a sliced manifold stratified pair
with respect to $\pi$   
with  $\dim(\pi ^{-1}(t)) = n$
for each $t\in \Delta ^k$.
\enddefinition

We will denote a typical $k$--simplex of $\SN^n(B)$ by
$\pi:(X,B\times\Delta^k)\to\Delta^k$ or, sometimes,
just by $\pi:X\to\Delta^k$ and consider the embeddings
$B\times\Delta^k\subseteq X$ and $X\subseteq\ell_2\times\Delta^k$
understood.
If $\pi:X\to\Delta^k$ is a $k$--simplex of $\SN^n(B)$, let 
$\partial X=\pi^{-1}(\partial\Delta^k)$ and let 
$\partial\pi = \pi|:\partial X\to\partial\Delta^k$,
Thus $\partial\pi:\partial X\to\partial\Delta^k$ is a union of 
$k+1$ $(k-1)$--simplices of $\SN^n(B)$.

The following result characterizes the homotopy relation in $\SN^n(B)$.
For notation, fix a base vertex of $\SN^n(B)$; that is, a manifold
stratified pair $(Y,B)$ with constant map $Y\to\Delta^0$. For
each $k\geq 0$ the degenerate $k$--simplex on $(Y,B)$  is
the  pair $(Y\times\Delta^k, B\times\Delta^k)$ with projection
$Y\times\Delta^k\to\Delta^k$.

\proclaim{Proposition 5.5} Let $B$ be a closed manifold.
Suppose $\pi:X\to\Delta^k$ and $\pi':X'\to\Delta^k$ are two simplices
of $\SN^n(B)$ such that
$\partial\pi=\partial\pi':\partial X=\partial X'=Y\times\partial
\Delta^k\to\partial\Delta^k$ is the projection.
The following are equivalent\rom:

\noindent
{\rm (i)} $\pi:X\to\Delta^k$ and $\pi':X'\to\Delta^k$ are homotopic
\rel $\partial$,

\noindent
{\rm (ii)} there exists a sliced manifold stratified pair
$(W,B\times\Delta^k\times I)$ with map $\widetilde{\pi}:W\to\Delta^k\times
I$ such that
\roster
\item{} $\widetilde{\pi}|=\pi:\widetilde{\pi}^{-1}(\Delta^k\times\{0\})
=X\to\Delta^k\times\{0\}=\Delta^k$,
\item{} $\widetilde{\pi}|=\pi':\widetilde{\pi}^{-1}(\Delta^k\times\{1\})
=X'\to\Delta^k\times\{1\}=\Delta^k$, and
\item{} $\widetilde{\pi}|=\partial\pi\times\id_I
=\partial\pi'\times\id_I=~ \proj :
\partial X\times I=\partial X'\times I= Y\times\partial\Delta^k\times I
\to\partial\Delta^k\times I$,
\endroster
{\rm (iii)} there exist an open neighborhood $U$ 
of $B\times\Delta^k$ in $X$
and a f\.p\. open embedding $h:U\to X'$ such that
$h|:(B\times\Delta^k)\cup(U\cap\partial X)\to(B\times\Delta^k)\cup(U\cap\partial X')$
is the identity.
\endproclaim

\demo{Proof}
(i) {\it implies\/} (ii):
Let $\widehat{\pi}:\widehat{W}\to \Delta^{k+1}$ be a homotopy 
\rel $\partial$ from
$\pi:X\to\Delta^k$ to $\pi:X'\to\Delta^k$ in $\SN^n(B)$.
Thus, $\widehat{\pi}=\pi$ over $\partial_{k+1}\Delta^{k+1}$,
$\widehat{\pi}=\pi'$ over $\partial_0\Delta^{k+1}$ and
$\widehat{\pi}|=~ \proj:Y\times\partial_i\Delta^{k+1}\to
\partial_i\Delta^{k+1}$ for $0<i<k+1$.
Consider the standard PL map
$\rho:\Delta^k\times I\to\Delta^{k+1}$
such that $\rho^{-1}(\partial\Delta^{k+1}) =\partial(\Delta^k\times I)$
and $\rho$ restricts to homeomorphisms $\Delta^k\times\{0\}\to
\partial_{k+1}\Delta^{k+1}$ and
$\Delta^k\times\{1\}\to
\partial_0\Delta^{k+1}$.
Form the pullback diagram
$$\CD W   @>{\widetilde{\pi}}>> \Delta^k\times I\\
      @VVV  @VV{\rho}V\\
      \widehat{W} @>{\widehat{\pi}}>> \Delta^{k+1}
\endCD$$
It follows from Lemma 5.3(iv)
that $(W,B\times\Delta^k\times I)$
is a sliced manifold stratified pair with map $\widetilde{\pi}$.

\noindent
(ii) {\it implies\/} (iii):
Let $V$ be an open neighborhood of $B\times\Delta^k\times I$ in $W$
such that $\widetilde{\pi}|:V\to \Delta^k\times I$ is a (trivial) fibre
bundle projection containing $B\times\Delta^k\times I\to \Delta^k\times I$
as a subbundle. Choose an open neighborhood $U$ of $B\times\Delta^k\times
\{0\}$ in $V\cap\widetilde{\pi}^{-1}(\Delta^k\times\{0\}) = X$ 
such that 
$$[U\cap\widetilde{\pi}^{-1}(\partial\Delta^k\times\{0\})]\times I
\subseteq V\cap\widetilde{\pi}^{-1}(\partial\Delta^k\times I)
\subseteq \partial X\times I=\partial X'\times I.$$
Let $J=(\Delta^k\times\{0\})\cup(\partial\Delta^k\times I)\subseteq
\Delta^k\times I$ and choose a homeomorphism
$\alpha:J\times I\to \Delta^k\times I$ such that
$\alpha|:J\times\{0\}\to J$ is the identity.
Since $\widetilde{\pi}|:V\to\Delta^k\times I$ is trivial, there exists a
homeomorphism
$g:[\widetilde{\pi}^{-1}(J)\cap V]\times I\to V$ such that
$$\CD  [\widetilde{\pi}^{-1}(J)\cap V]\times I
            @>g>> V\\
@V{\widetilde{\pi}\times\id_I}VV @VV{\widetilde{\pi}|}VV\\
      J\times I @>{\alpha}>>  \Delta^k\times I
\endCD$$
commutes, $g|B\times J\times I$ equals $\id_B\times\alpha:
B\times J\times I\to B\times\Delta^k\times I\subseteq V$, and
$g|:[\widetilde{\pi}^{-1}(J)\cap V]\times\{0\}\to\widetilde{\pi}^{-1}(J)\cap
V$ is the identity.
Define $h:U\to\widetilde{\pi}^{-1}(\Delta^k\times\{1\}) = X'$
by setting $h(x) = g(x,1)$ for all $x\in U$.

\noindent
(iii) {\it implies\/} (i):
Let $N$ be a compact neighborhood of $B\times\Delta^k$ in $X$ such that
$N\subseteq U$. By the small capacity assumption, there exists a f\.p\.
isotopy $H_t:\ell_2\times\Delta^k\to\ell_2\times\Delta^k$, $0\leq t\leq 1$,
such that $H_0=\id_{\ell_2}$,
$H_t|(B\times\Delta^k)\cup\ell_2\times\partial\Delta^k$ is the
identity for each $t\in I$, 
and $H_1|N = h|N:N\to X'\subseteq\ell_2\times\Delta^k$.
Let 
$$W=(\partial X\times I)\cup(X\times\{0\})\cup(X'\times\{1\})\cup
\{(H_t(x),t) ~|~ x\in{\inr}(N), t\in I\}.$$
Proposition 5.2 implies that
$(N\times I,B\times\Delta^k\times I)$ is a sliced homotopically
stratified pair with finitely dominated local holinks, which in turn implies
that $(W,B\times\Delta^k\times I)$ is a sliced manifold stratified pair.
Now $W$ induces a sliced manifold stratified pair
$(\widehat{W}, B\times\Delta^{k+1})$ such that $W$ is the pullback of
$\widehat{W}$ along the map $\rho:\Delta^k\times I\to \Delta^{k+1}$ of (i),
and $(\widehat{W},B\times\Delta^{k+1})$ is the desired homotopy from
$X$ to $X'$ \rel $\partial$.
\qed
\enddemo

The next result follows from Proposition 5.5 by setting $k=0$.

\proclaim{Corollary 5.6} Let $B$ be a closed manifold.
Two vertices $(X,B)$, $(X',B)$ are in the same component of
$\SN^n(B)$ if and only if they are germ equivalent\rom; that is, there
exist an open neighborhood $U$ of $B$ in $X$ and an open embedding
$h:U\to X'$ such that $h|:B\to X'$ is the inclusion.
\qed
\endproclaim

In order for homotopy theory to work well on the space of stratified
neighborhoods, we need the following observation.

\proclaim{Proposition 5.7} $\SN^n(B)$ satisfies the Kan condition.
\endproclaim

\demo{Proof} Suppose there is a collection of $k+1$ $k$--simplices
$(X_j,B\times\partial_j\Delta^{k+1})$ of $\SN^n(B)$,
$j=0,1,\dots,i-1,i+1,\dots,k+1$, which satisfy the compatibility condition
(see \cite{41\rm, p\. 2}). For $X=\cup X_j$ there is a natural map
$X\to B\times w_i\Delta^{k+1}$ where $w_i\Delta^{k+1}$ is the union of
all $k$--dimensional faces of $\Delta^{k+1}$ save $\partial_i\Delta^{k+1}$.
It is elementary to verify that $(X,B\times w_i\Delta^{k+1})$ is a sliced
manifold stratified pair. A possible exception is in the verification of
the holink fibration condition, but that condition follows from 
\cite{29\rm, 16.2}. Pulling back along a retraction  $\Delta^{k+1}\to w_i\Delta^{k+1}$
gives (by Lemma 5.3) a sliced manifold stratified pair
$(\widetilde{X}, B\times\Delta^{k+1})$ which is the required $(k+1)$--simplex
of $\SN^n(B)$.
\qed
\enddemo

Now recall the following definition  from \cite{29}.

\definition{Definition 5.8} The {\it
space of manifold approximate fibrations over\/}
$B\times \br $ is the simplicial set $\maf^n
(B\times \br) $ whose $k$-simplices are subsets
$M$ of $\ell_2 \times B\times \br \times \Delta
^k$ of small capacity such that
\widestnumber\item{(ii)}
\roster
\item"(i)" the restriction of projection $M\to \Delta ^k$ is
a fibre bundle projection with fibres $n$-dimensional
manifolds without boundary.  Let $M_t $ denote the
fibre over $t\in \Delta ^k$.
\item"(ii)" the restriction of projection $p: M\to B \times
\br \times \Delta ^k$ has the property that $p_t
= p|: M_t\to B\times \br \times \{ t\}$ is a
manifold approximate fibration for each $t\in \Delta ^k$.
\endroster
\enddefinition

We will denote a typical $k$--simplex of $\maf^n(B\times\br)$ by
$p: M\to B\times \br \times \Delta ^k$ and consider the 
embeddings
$B\times\Delta^k\subseteq X$ and $X\subseteq\ell_2\times\Delta^k$
understood.

\subhead Definition of\/ $\Psi:\maf^n(B\times\br)\to\SN^n(B)$ \endsubhead
It will be convenient to fix a teardrop of $B$ in
$\ell_2$ which contains all the teardrops constructed
form $\maf^n(B\times \br )$.  To this end let
$$\mu : \ell_2\times B\times \br \to B\times
\br $$
denote projection and let 
$$T(B) = (\ell _2 \times B\times \br )
     \cup_\mu B$$
be the teardrop of $\mu$. It follows from Lemma 4.3
that $T(B)$ is metrizable.  Since $B$ is separable,
$T(B)$ is also separable.  Hence, $T(B)$ embeds in
$\ell_2$ and we fix an embedding $T(B)\subseteq \ell_2$
of small capacity such that $B\subseteq T(B)\subset
\ell_2$ is the original fixed embedding $B\subseteq \ell_2$.

We now define the simplicial map $\Psi : \maf^n
(B\times \br )\to \SN^n (B)$.  Given a
$k$-simplex $M\subseteq \ell_2 \times B\times \br \times
\Delta ^k$ of $\maf^n (B\times \br )$, we get a
commuting diagram 
$$\CD M @>{\subseteq}>> \ell_2 \times B\times\br \times \Delta ^k\\
    @VpVV     @VV{\mu\times\id_{\Delta^k}}V\\
{B\times\br\times\Delta^k} @>=>> {B\times\br\times\Delta^k}
\endCD$$
Thus, $M\cup_p(B\times \Delta ^k)
\subseteq (\ell_2\times B\times \br \times \Delta ^k)
     \cup_{\mu\times \id_{\Delta^k}}
     (B\times \Delta ^k) = T(B)\times \Delta ^k 
     \subseteq \ell_2 \times \Delta ^k$.  It will be shown
below that $(M\cup_p (B\times \Delta ^k), B\times\Delta^k)$
is a $k$-simplex of $\SN^n(B)$, and so we set $\Psi (M) = 
(M\cup_p (B\times \Delta ^k), B\times\Delta^k)$.

\subhead Proof that $\Psi(M)$ is a $k$--simplex of $\SN^n(B)$ \endsubhead
It is clear from the construction that $M\cup_p(B\times\Delta^k)$ is
a subset of $\ell_2\times\Delta^k$ of small capacity.
Since each $p_t : M_t \to B\times
\br \times \{ t\}$ is a manifold approximate
fibration, it follows from Corollary 4.11 that
$(M_t \cup_{p_t} B\times
\{ t\}, B\times \{ t\})$ is a manifold stratified pair
for each $t\in \Delta ^k$.  Therefore,
the sliced forward tameness, holink fibration and finitely dominated local
holinks conditions follow from Claim 5.9 and Lemmas 5.10 and 5.11 below.
To verify the local triviality condition
let ${\cal U}$ be the open cover of $B\times
\br $ consisting of all sets of the form
$$B(x, {1\over \vert y\vert +\vert } ) \times
     (y- {1\over \vert y\vert +1}, y +
     {1\over \vert y\vert + 1} ) $$
where $(x,y) \in B\times \br$ and $B(x,r)$
denotes the ball about $x$ in $B$ of radius $r$.  The
point is that the diameters of members of ${\cal U}$
are small near $B\times \{ + \infty \}$ and there is a
maximum diameter.
By \cite{24} there is a homeomorphism $H: M\times \Delta^k
\to M\times \Delta ^k$ such that $H$ is fibre
preserving over $\Delta^k,\  H_0 = id$, and $pH$ is
${\cal U}\times \Delta ^k$-close to $p_0 \times \id
_{\Delta^k}$.
The local triviality condition follows from the following claim and
the fact that
$(M\cup_{p_0} B)\times\Delta ^k =
(M\times \Delta ^k) \cup_{p_0\times \id
_{\Delta ^k}} (B\times \Delta ^k)$.

\proclaim{Claim 5.9} The map $h: (M\times \Delta ^k) \cup_{
p_0\times \id_{\Delta^k}} (B\times \Delta ^k) \to (M\times \Delta
^k ) \cup_{p} (B\times \Delta ^k )$,
defined by $h\mid : M\times \Delta ^k \to M\times
\Delta ^k $ is $H$ and $h\mid : B\times \Delta ^k\to
B\times \Delta ^k$ is the identity, is a homeomorphism.
\endproclaim

\demo{Proof}  We show that the map
$$g: (M\times \Delta ^k)\cup_{p_0\times \id_{\Delta^k}}(B\times \Delta^k)
@>h>> (M\times \Delta ^k)\cup_p (B\times \Delta^k)
@>c>> B\times (-\infty , +\infty ] \times\Delta ^k$$
is continuous with $c$ the teardrop collapse for $p$.  
For this it suffices to show that if
$(x_n, t_n) \in M\times \Delta ^k,\  (b,t)\in B\times
\Delta ^k$ and $(x_n,t_n)\to (b,t)$ in $(M\times
\Delta ^k) \cup_{p_0\times \id}
(B\times \Delta ^k)$, then $g(x_n, t_n)\to (b, +
\infty , t)$ in $B\times (-\infty ,+\infty ]\times
\Delta ^k$.  Let $c': (M\times \Delta ^k) {\cup}
_{p_0 \times \id} (B\times \Delta ^k) \to B\times
(-\infty , + \infty ] \times \Delta ^k$ be the collapse.
Since $c'$ is continuous, $c'(x_n, t_n)\to (b,
+ \infty , t)$ and so $(p_0 (x_n), t_n) \to (b,
+\infty , t)$.
Given $\epsilon > 0$ there exists an integer $K$ such
that if $U \in {\cal U} $ meets $B\times [K, +\infty
)$, then $\diam U < \epsilon $.  There exists a
positive integer $M$ such that if $n\ge M$, then $p_0
(x_n) \in B\times [K,+\infty )$ and $(p_0 (x_n),t_n)$ 
is $\epsilon$-close to $(b,+\infty , t )$.
Now suppose $n\ge M$ and consider $g(x_n, t_n)$.  Note
that $g(x_n, t_n) = pH(x_n, t_n)$.  There exists $U\in
{\cal U}$ such that $pH(x_n,t_n)$ and $(p_0 (x_n), t_n)$
are both in $U\times \Delta ^k$; i\.e\., $p_{t_n}
H_{t_n}(x_n),\ p_0(x_n)\in {\cal U}$.  Since $p_0 (x_n)
\in B\times [K,+\infty ),\diam U< \epsilon $.  Thus, $pH(
x_n, t_n)$ and $(p_0 (x_n), t_n)$ are
$\epsilon$-close measured in $B\times \br\times \Delta^k$.
Since $(p_0 (x_n), t_n)$ is $\epsilon$-close to $(b,
+\infty , t)$, we have shown that $g(x_n,t_n)$ is 
$\epsilon'$-close to $(b,+\infty , t)$ where
$\epsilon'> 0$ is small if $\epsilon $ is.  Thus,
$g$ is continuous.
This shows $h$ is continuous by Lemma 3.4.  Since $p$
is also ${\cal U}\times \Delta ^k$-close to $(p_0
\times \Delta ^k)H^{-1}$, a similar proof shows that
$h^{-1}$ is continuous.
\qed
\enddemo

We finish this section with the two lemmas mentioned above.

\proclaim{Lemma 5.10} Suppose $B$ is forward tame in $X$. 
\widestnumber\item{(ii)}
\roster
\item"(i)" If $Y$ is any space, then $B\times Y$ is
sliced forward tame in $X\times Y$ with respect to
projection $X\times Y\to Y$.

\item"(ii)" If $\pi : E\to Y$ is a map of spaces and
$h: X\times Y\to E$ is a homeomorphism such that $\pi
h$ is projection, then $h(B\times Y)$ is sliced
forward tame in $E$ with respect to $\pi$.
\endroster
\endproclaim

\demo{Proof} (i) is obvious, and (ii) follows from (i)
by using a sliced nearly strict deformation in
$X\times Y$ conjugated with $h$.
\qed
\enddemo

\proclaim{Lemma 5.11} Suppose $B\subseteq X$ and $\ho (X,B)
\to B$ is a fibration.
\widestnumber\item{(ii)}
\roster
\item"(i)" If $Y$ is any space, then $\ho_{p_2}
     (X\times Y, B\times Y)\to B\times Y$ is a
fibration where $p_2$ is second coordinate projection.

\item"(ii)" If $\pi : E\to Y$ is a map of spaces and $h:
X\times Y\to E$ is a homeomorphism such that $\pi h$
is projection, then $\ho_\pi (E, h(B\times Y))\to
h (B\times Y)$ is a fibration.
\endroster
\endproclaim

\demo{Proof} For (i) note that we have the
following commuting diagram where $\nu (\omega) = (\omega', p_2
\omega (0))$ and $\omega'$ is 
$[0,1] @>{\omega}>> X\times Y @>{\proj}>> X\ $:
$$\CD \ho_{p_2}(X\times Y,B\times Y) @>{\nu}>> \ho(X,B)\times Y\\
@VVV @VVV\\
 B\times Y @>=>>  B\times Y.
\endCD$$
For (ii) note that we have the following commuting
diagram where $\lambda $ is the homeomorphism defined
by $\lambda (\omega ) = h\circ \omega\ $ :
$$\CD  \ho_{p_2}(X\times Y, B\times Y)  @>{\lambda}>> \ho_\pi(E,h(B\times Y))\\
@VVV @VVV\\
 B\times Y @>h>>  h(B\times Y).\qed
\endCD$$
\enddemo

\head 6. Homotopy near the lower stratum \endhead
The main theorems of this paper on Teardrop Neighborhood Existence (2.1)
and Neighborhood Germ Classification (2.2, 2.3) have two aspects in their
proofs: homotopy theoretic and manifold theoretic.
This is already evident in \S4 if one compares Theorem 4.7,
which says that the teardrop of an approximate fibration is a homotopically
stratified pair, with Corollary 4.11 which says that the teardrop of a 
manifold approximate fibration is a manifold stratified pair.
This section contains the homotopy theoretic part of the remaining
aspects of this paper's main existence and classification theorems. The main
result here, Theorem 6.8, produces from a homotopically stratified
pair $(X,A)$ with finitely dominated local holinks,
a $\cal U$--fibration over $A\times (0,+\infty)$ for arbitrarily 
small open covers $\cal U$ of $A\times (0,+\infty)$
(outside the setting of manifolds this is not quite the same notion as 
an approximate fibration).
The proof involves showing that the mapping cylinder of the holink
evaluation is a good homotopy model for a neighborhood germ of $A$ in
$X$. The idea of a good homotopy model is made precise with the notion
of a `strong ${\cal U}$--homotopy equivalence near $A$' in
Definition 6.1.

There are three main steps to the proof of 6.8 corresponding to the 
three main hypotheses: holink evaluation is a fibration, forward tameness
and finitely dominated local holinks. The first step is Proposition 6.3
which shows how being modelled on the mapping cylinder of a fibration
yields $\cal U$--fibrations (we apply this to the holink evaluation
fibration).
The second step, Proposition 6.5, shows that forward tameness is enough to
get started in showing that the mapping cylinder of holink evaluation
is a good model for a neighborhood of $A$ in $X$.
Finally, the third step, Proposition 6.7, adds the finitely dominated local 
holinks condition to produce the strong $\cal U$--homotopy equivalence
near $A$.
Of course, all of this must be done sliced (or fibre preserving) over
$\Delta^k$ in order to obtain the Higher Classification Theorem 2.3.

We begin with the following definition of strong homotopy 
equivalences near $A$.

\definition{Definition 6.1} Suppose $X_1$ and $X_2$ are spaces containing
$A\times\Delta^k$ with maps $\pi_i:X_i\to\Delta^k$ such that
$\pi_i|:A\times\Delta^k\to\Delta^k$ is projection for $i=1,2$.
Suppose $p:X_2\to A\times (-\infty,+\infty]\times\Delta^k$
is a map which is fibre preserving over $\Delta^k$ and
such that $p^{-1}(A\times\{+\infty\}\times\Delta^k) = A\times
\Delta^k$ and $p|:A\times\Delta^k\to A\times\{+\infty\}\times\Delta^k$
is the identity. Suppose $\cal U$ is an open cover of $A\times\br\times
\Delta^k$. A {\it strong f\.p\. $\cal U$--homotopy equivalence near\/} 
$A\times\Delta^k$
$$(f,g,X_1',X_2'): X_1\to X_2$$
is defined by maps
$$f:X_1'\to X_2, \qquad g:X_2'\to X_1$$
such that
\widestnumber\item{(iii)}
\roster
\item"(i)" $X_1'$ a closed neighborhood of $A\times\Delta^k$ in $X_1$
and $X_2' = p^{-1}(A\times [t_2,+\infty]\times\Delta^k)$ 
for some $t_2\in\br$,
\item"(ii)" the maps 
$$\align f:(X_1',A\times\Delta^k) &\to (X_2,A\times\Delta^k),\\
g:(X_2',A\times\Delta^k) &\to (X_1,A\times\Delta^k)
\endalign$$
are fibre preserving over $\Delta^k$, strict and the identity on 
$A\times\Delta^k$,
\endroster
together with  homotopies
$$\align F:gf| &\simeq \text{\rm inclusion}:f^{-1}(X_2')\to X_1,\\
G:fg| &\simeq \text{\rm inclusion}:g^{-1}(X_1')\to X_2
\endalign$$
such that 
\widestnumber\item{(iii)}
\roster
\item"(iii)" $F,G$ are fibre preserving over $\Delta^k$, 
rel $A\times\Delta^k$, and
strict as homotopies between pairs
$(f^{-1}(X_2'), A\times\Delta^k)\times I\to (X_1, A\times\Delta^k)$ 
and
$(g^{-1}(X_1'), A\times\Delta^k)\times I\to (X_2, A\times\Delta^k)$,
\item"(iv)" for every $x\in f^{-1}(X_2')\setminus (A\times\Delta^k)$ 
with $\{ x\}\times I\subseteq F^{-1}(X_1')$ there exists
$U\in {\cal U}$  such that
$pfF(\{ x\}\times I)\subseteq U$,
\item"(v)" for every $x\in g^{-1}(X_1')\setminus (A\times\Delta^k)$ 
there exists
$U\in {\cal U}$ such that
$pG(\{ x\}\times I)\subseteq U$.
\endroster
\enddefinition

Sliced homotopy lifting properties are just the parametric versions of
ordinary lifting properties. These are used to define
sliced $\cal U$--fibrations, sliced approximate fibrations and sliced
manifold approximate fibrations (see \cite{22}). We include the following 
definition for completeness.

\definition{Definition 6.2}
Suppose $p:E\to A\times\Delta$ is a map (with $\Delta$ playing the role
of the parameter space), $V\subseteq A\times\Delta$ and $\cal U$ is an
open cover of $A\times\Delta$. Then $P$ is a
{\it sliced $\cal U$--fibration over\/} $V$ if for every commuting diagram
of maps which are f\.p\. over $\Delta$
$$\CD
Z\times\Delta @>f>> E \\
@V{\times 0}VV  @VVpV\\
Z\times\Delta\times I @>F>> A\times\Delta
\endCD$$
with $\im(F) \subseteq V$, there exists an f\.p\. (over $\Delta$)
map $\tilde F:Z\times\Delta\times I\to E$ such that
$\tilde F_0=f$ and $p\tilde F$ is $\cal U$--close to $F$.
If $V=A\times\Delta$, then $p$ is a {\it sliced $\cal U$--fibration}.
If $p$ is a sliced $\cal U$--fibration for every open cover $\cal U$, then
$p$ is a {\it sliced approximate fibration}. If $E\to\Delta$ is a fibre
bundle projection with manifold fibres (without boundary), $A$ is a 
manifold (without boundary) and $p$ is a proper sliced approximate
fibration, then $p$ is said to be a {\it sliced manifold approximate
fibration}.
\enddefinition

A map $p:E\to A$ is {\it proper over a subspace\/} $V\subseteq A$ if
for every compact subspace $K\subseteq V$, $p^{-1}(K)$ is compact.
We do not insist that proper maps be onto.

The following result shows that it is significant to be strongly
f\.p\. $\cal U$--homotopy equivalent to the mapping cylinder of a fibration
near the base of the mapping cylinder.

\proclaim{Proposition 6.3}
Suppose $q:E\to A\times\Delta^k$ is a fibration and
$Q:\ocyl(q)\to A\times (-\infty,+\infty]\times\Delta^k$
is the teardrop collapse.
Suppose $X$ is a 
locally compact separable metric
space containing $A\times\Delta^k$ with a map
$\pi:X\to\Delta^k$ such that 
$\pi|:A\times\Delta^k\to\Delta^k$ is projection and
$\cal U$ is an open cover of $A\times\br\times\Delta^k$.
Suppose
$(f,g,X_1',X_2'):X\to\ocyl(q)$
is a strong f\.p\. $\cal U$--homotopy equivalence near
$A\times\Delta^k$
and $Qf:X_1'\to A\times (-\infty,+\infty]\times\Delta^k$ is proper.
Then there exists an open neighborhood $V$ of 
$A\times\{+\infty\}\times\Delta^k$
in $A\times (-\infty,+\infty]\times\Delta^k$ such that
$Qf:X_1'\to A\times (-\infty,+\infty]\times\Delta^k$
is a sliced $\st^2({\cal U})$--fibration over
$(A\times \br\times\Delta^k)\cap V$.
\endproclaim

\demo{Proof}
If $X_2' = Q^{-1}(A\times [t_2,+\infty]\times\Delta^k)$ choose an
open neighborhood $V$ of $A\times\{+\infty\}\times\Delta^k$ in
$A\times (-\infty,+\infty]\times\Delta^k$ such that
\widestnumber\item{(iii)}
\roster
\item"(i)" $V\subseteq A\times [t_2,+\infty]\times\Delta^k$,
\item"(ii)" $Q^{-1}(V)\subseteq g^{-1}(X_1')$\ (this is possible since $Q$
is a closed map over $A\times\Delta^k$),
\item"(iii)" $(Qf)^{-1}(V)\subseteq f^{-1}(X_2')$\ (this is possible
since $Qf$ is proper, hence a closed map), and
\item"(iv)" $(Qf)^{-1}(V)\times I\subseteq F^{-1}(X_1')$\ (this is 
possible since $Qf$ is proper and $F$ is the identity on $A\times\Delta^k
\times I$).
\endroster
A sliced homotopy lifting problem
$$\CD
Z\times\Delta^k @>d>> X_1' \\
@V{\times 0}VV  @VV{Qf}V \\
Z\times\Delta^k\times I @>D>> A\times (-\infty,+\infty]\times
\Delta^k
\endCD$$
with ${\im}(D)\subseteq (A\times\br\times\Delta^k)\cap V$
yields another lifting problem
$$\CD
Z\times\Delta^k @>{fd}>> \ocyl(q)\setminus (A\times\Delta^k) \\
@V{\times 0}VV   @VV{Q|}V \\
Z\times\Delta^k\times I @>D>> A\times (-\infty,+\infty]\times
\Delta^k.
\endCD$$
Since $\ocyl(q)\setminus(A\times\Delta^k) = E\times\br$ and
$Q| = q\times\id_\br$ is a fibration, this second problem has an exact
solution
$\widetilde{D}^1:Z\times\Delta^k\times I\to E\times\br$ (so that
$\widetilde{D}^1|Z\times\Delta^k\times\{ 0\} = fd$ and $(q\times\id_\br)
\widetilde{D}^1 = D$).
By choice of $V$, ${\im}(\widetilde{D}^1)\subseteq X_2'$
and ${\im}(g\widetilde{D}^1) \subseteq X_1'$.
Define $\widetilde{D}^2 = g\widetilde{D}^1: Z\times\Delta^k\times I\to
X_1'$ and note that 
$Qf\widetilde{D}^2 = Qfg\widetilde{D}^1$ is $\cal U$--close to
$Q\widetilde{D}^1 = D$.
Except for the fact that $\widetilde{D}^2|Z\times\Delta^k\times\{0\}$ need
not equal $d$, $\widetilde{D}^2$ would be an approximate solution to the
original problem. However, $\widetilde{D}^2|Z\times\Delta^k\times\{0\}
= g\widetilde{D}^1| = gfd$ and $gfd$ is $(Qf)^{-1}({\cal U})$--homotopic to
$d$. Thus a standard argument using paracompactness allows a modification of
$\widetilde{D}^2$ to get a $\st^2({\cal U})$--solution
$\widetilde{D}:Z\times\Delta^k\times I\to X_1'$ (see \cite{28\rm, Prop\. 16.3}).
\qed
\enddemo

\remark{Notation \rom{6.4}}
For the remainder of this section 
suppose $A\times\Delta^k\subseteq X$
and $\pi:X\to\Delta^k$ is a map such that
$\pi|:A\times\Delta^k\to\Delta^k$ is the projection and
$q:\ho_\pi(X,A\times\Delta^k)\to A\times\Delta^k$
is the evaluation.
The open mapping cylinder of $q$ is identified with the teardrop
$$\ocyl(q) = (\ho_\pi(X,A\times\Delta^k)\times\br)\cup_{q\times\id} 
(A\times\Delta^k)$$
where $q\times\id : \ho_\pi(X,A\times\Delta^k)\times\br\to
A\times\br\times\Delta^k$.
Let $Q:\ocyl(q)\to A\times (-\infty,+\infty]\times\Delta^k$ be the
teardrop collapse.
\endremark

The genesis of the ideas in the next two results is in \cite{24\rm, 4.7}
and \cite{48\rm, 2.4}. See especially \cite{28\rm, 9.13, 9.14}.

\proclaim{Proposition 6.5}
Suppose $X$ is a locally compact separable metric space, $A$ is compact and
$A\times\Delta^k$ is sliced forward tame in $X$ with respect to $\pi$.
Then there exist a compact neighborhood $Y$ of $A\times\Delta^k$ in $X$ and
maps
$$f:Y\to\ocyl(q), \qquad
g:\ocyl(q)\to Y$$
together with homotopies
$$F:igf\simeq i:Y\to X, \qquad
G:fg\simeq\id:\ocyl(q)\to\ocyl(q)$$
with $i:Y\to X$ the inclusion such that
\widestnumber\item{(iii)}
\roster
\item"(i)" $f,g,F,G$ are \rel $A\times\Delta^k$,
\item"(ii)" $f,g,F,G$ are f\.p\. over $\Delta^k$,
\item"(iii)" $f,g,F,G$ are strict maps or homotopies between the pairs
$(X,A\times\Delta^k)$ and $(\ocyl(q),A\times\Delta^k)$,
\item"(iv)" for every $N\geq 0$ there exists $M\geq 0$
such that
$$(Qfg)^{-1}(A\times(-\infty,N]\times\Delta^k)\subseteq Q^{-1}(A\times
(-\infty,M]\times\Delta^k),$$
\item"(v)" for every $N\geq 0$ there exists $M\geq 0$ 
such that
$$G(Q^{-1}(A\times [M,+\infty]\times\Delta^k)\times I)\subseteq
Q^{-1}(A\times [N,+\infty]\times\Delta^k).$$
\endroster
\endproclaim

\demo{Proof} (cf\. \cite{28\rm,9.13}.)
Let $d$ be a metric for $X$ and
let $Y$ be a compact neighborhood of $A\times\Delta^k$ in $X$ for which there
exists a nearly strict deformation
$H:(Y\times I, A\times\Delta^k\times I\cup Y\times\{0\})\to
(X,A\times\Delta^k)$ of $Y$ into $A\times\Delta^k$ which is f\.p\. 
over $\Delta^k$. 
It is easy to modify $H$ so that it has the additional property
that if $N=1,2,3,\dots$ and $x\in H(Y\times [0,1/N])$, then
$d(x,A)\leq 1/N$.
Let $\hat H:Y\setminus(A\times\Delta^k)
\to\ho_\pi(X,A\times\Delta^k)$ be the adjoint
of $H$.
Choose a compact neighborhood $Y'$ of $A\times\Delta^k$ in $X$ such that
$Y'\subseteq Y$ and $\hat H(Y')\subseteq\ho_\pi(Y,A\times\Delta^k)$.
Use $i$ also to denote the inclusion $i:Y'\to X$.
From Proposition 5.2(i), it induces a fibre homotopy equivalence
$i_\ast:\ho_\pi(Y',A\times\Delta^k)\to\ho_\pi(X,A\times\Delta^k)$.
Let $R:\ho_\pi(X,A\times\Delta^k)\times I\to
\ho_\pi(X,A\times\Delta^k)$ be the fibre deformation explicitly defined in 
5.2. Thus, there is a fibre homotopy inverse
$j:\ho_\pi(X,A\times\Delta^k)\to\ho_\pi(Y',A\times\Delta^k)$ for $i_\ast$
defined by $j=R_1$. From the definition of $R$, we have 
$R(\omega,t)(u)=\omega(s)$ for some $s$.
Define $p:X\to (0,+\infty]$ by $p(x)= 1/d(x,A)$.
Define  $f:Y\to\ocyl(q)$ by
$$f(x) =\cases (\hat H(x),p(x))\in\ho_\pi(X,A\times\Delta^k)\times 
       (0,+\infty), & \text{if $x\in Y\setminus(A\times\Delta^k)$}\\
                    x, & \text{if $x\in A\times\Delta^k$.}
\endcases$$
Let $p_{Y'}:\ho_\pi(Y',A\times\Delta^k)\to Y'$ and
$p_Y^+:\ho_\pi(Y,A\times\Delta^k)\times\br\to Y$
be the evaluations
$p_{Y'}(\omega)=\omega(1)$ and
$$p_Y^+(\omega,t)=\cases \omega(1), & \text{if $t\leq 0$}\\
                         \omega(1/(1+t)), & \text{if $t\geq 0$.}
\endcases$$
Define $g:\ocyl(q)\to Y$ so that on 
$\ho_\pi(X,A\times\Delta^k)\times\br\subseteq\ocyl(q)$, 
$g$ is the composition
$$\ho_\pi(X,A\times\Delta^k)\times\br @>{j\times\id_\br}>>
\ho_\pi(Y',A\times\Delta^k)\times\br
@>{p_{Y'}\times\id_\br}>>{}$$
$$Y'\times\br @>{\hat H\times\id_\br}>>
\ho_\pi(Y,A\times\Delta^k)\times\br
@>{p_Y^+}>> Y$$
and on $A\times\Delta^k\subseteq\ocyl(q)$,\ $g$ is the identity.
Define the homotopy $F:Y\times I\to X$ by
$$F(x,t) =\cases
(\hat H[(R_{1-t}(\hat H(x)))(1)])({d(x,A)+t\over d(x,A)+1}), & \text{if $x\in 
Y\setminus(A\times\Delta^k)$}\\
x, & \text{if $x\in A\times\Delta^k$.}
\endcases$$
Define
$$\gamma:\ho_\pi(X,A\times\Delta^k)\times (0,1]\to\ho_\pi(X,A\times\Delta^k)$$
by $\gamma(\omega,t)=\hat H[\hat H(x_\omega)(t)]$
where $x_\omega=j(\omega)(1)\in Y'$.
Define $G':\ho_\pi(X,A\times\Delta^k)\times\br\times I\to\ho_\pi(X,A\times
\Delta^k)\times\br$
by
$$G'(\omega,t,s)=\cases (\gamma(\omega,{1\over 1+t}),(1-s)p[\hat H(x_\omega)
({1\over 1+t})]+st), & \text{if $s\geq t$}\\
(\gamma(\omega,1),(1-s)p[\hat H(x_\omega)
(1)]+st), & \text{if $s\geq t$.}
\endcases$$
Note that
$G'_0=fg|:
\ho_\pi(X,A\times\Delta^k)\times\br
\to\ho_\pi(X,A\times\Delta^k)\times\br$
and that $G'$ extends via the indentity on $A\times\Delta^k$ to
$G':\ocyl(q)\times I\to\ocyl(q).$
We claim that there exists a homotopy
$$G'':\ho_\pi(X,A\times\Delta^k)\times\br\times I\to
\ho_\pi(X,A\times\Delta^k)\times\br$$
such that
$$G''_0(\omega,t)=\cases \gamma(\omega,{1\over 1+t}), & \text{if $t\geq 0$}\\
\gamma(\omega,1), & \text{if $t\leq 0$.}
\endcases$$
To this end note that by
contracting $(0,1]$ to $\{ 1\}$ there is defined a homotopy
$\gamma\simeq\gamma'$ with
$$\gamma'(\omega,t) =\hat H[\hat H(w_\omega)(1)] = \hat H(x_\omega)=
\hat H(p_{Y'}(j(\omega)).$$
And it is not difficult to see that 
$\hat Hp_{Y'}:\ho_\pi(Y',A\times\Delta^k)\to\ho_\pi(Y,A\times\Delta^k)$
is homotopic to the  inclusion $i_\ast$.
Since $j$ is a homotopy inverse for $i_\ast$, the homotopy 
$G''$ exists as claimed.
We can now define the homotopy
$$G:\ho_\pi(X,A\times\Delta^k)\times\br\times I\to\ho_\pi(X,A\times\Delta^k)
\times\br$$
by
$$(\omega,t,s)\mapsto\cases G'(\omega,t,2s), & \text{if $ 0\leq s\leq 1/2$}\\
(G''(\omega,t,2s-1),t), & \text{if $1/2\leq s\leq 1$}
\endcases$$
and extending $G$ via the identity on $A\times\Delta^k$ to
get 
$$G:\ocyl(q)\times I\to \ocyl(q).$$
For the verification of the properties, see \cite{28\rm, 9.13}.          
\qed
\enddemo

\proclaim{Lemma 6.6}
Let $p:E\to B$ be a fibration with $B$ a weakly locally contractible compact
metric space. If the fibre of $p$ is finitely dominated, then there exist
a compact subspace $K\subseteq E$ and a f\.p\. homotopy
$D:E\times I\to E$ such that 
$D_0(E)\subseteq K$ and
$D_1=\id_E$.
\endproclaim

\demo{Proof}
Each $x\in B$ has an open neighborhood $U_x$ such that the inclusion
$U_x\hookrightarrow B$ is null--homotopic.
It follows that there is a fibre homotopy equivalence
$f_x:p^{-1}(U_x)\to p^{-1}(x)\times U_x$ over $U_x$.
Let $g_x:p^{-1}(x)\times U_x\to p^{-1}(U_x)$
be a fibre homotopy inverse and
$H^x:p^{-1}(U_x)\times I\to p^{-1}(U_x)$
a f\.p\. homotopy such that
$H_0^x=g_xf_x$ and $H_1^x=\id_{p^{-1}(U_x)}$.
Since $p^{-1}(x)$ is finitely dominated there exist a compact subspace
$K_x\subseteq p^{-1}(x)$ and a homotopy
$D^x:p^{-1}(x)\times I\to p^{-1}(x)$
such that 
$D_0^x(p^{-1}(x))\subseteq K_x$ and $D_1^x=\id_{p^{-1}(x)}$.
Let $\hat D^x=D^x\times\id_{U_x}:p^{-1}(x)\times U_x\times I\to
p^{-1}(x)\times U_x$.
Let $\rho_x:B\to I$ be a map such that $\rho_x^{-1}(0)\subseteq U_x$ is
a neighborhood of $x$ and $B\setminus U_x\subseteq \rho_x^{-1}(1)$.
Define a f\.p\. homotopy
$G^x:E\times I\to E$ by
$$G^x(y,t)=\cases g_x\hat D^x(f_x(y),(1-t)2\rho_x(y)+t),
& \text{if $0\leq\rho_x(y)\leq 1/2$}\\
H^x(y,t(2\rho_x(y)-1)+(1-t)),
& \text{if $1/2\leq\rho_x(y)\leq 1$.}
\endcases$$
Define a f\.p\. homotopy
$F^x:E\times I\to E$ by
$$F^x(y,t)=\cases H^x(y,t),
& \text{if $0\leq\rho_x(y)\leq 1/2$}\\
H^x(y,(1-t)(2\rho_x(y)-1)+t),
& \text{if $1/2\leq\rho_x(y)\leq 1$.}
\endcases$$
Then $F_0^x=G_1^x$ and $F_1^x=\id_E$.
Define a f\.p\. homotopy
$\tilde D^x:E\times I\to E$ by
$$\tilde D^x(y,t)=\cases G^x(y,2t),
& \text{if $0\leq t\leq 1/2$}\\
F^x(y,2t-1),
& \text{if $1/2\leq t\leq 1$.}
\endcases$$
Then $\tilde D_0^x=G_0^x$ and $\tilde D_1^x=\id_E$.
The compact subspace $C_x=g_x(K_x\times p^{-1}(\rho_x^{-1}(0)))$
of $E$
is such that $\tilde D_0^x(\rho_x^{-1}(0))\subseteq C_x$.
Let $\{ x_1,\dots,x_k\}$ be a finite subset of $B$ such that
$B=\cup_{i=1}^k\rho_{x_i}^{-1}(0)$.
Define $D:E\times I\to E$ by
$$D_t=\tilde D_t^{x_k}\circ\cdots\circ\tilde D_t^{x_1}.$$
Then $D_1=\id_E$ and
$$D_0(E)\subseteq
[\tilde D_0^{x_k}\circ\cdots\circ\tilde D_0^{x_2}(C_{x_1})]\cup
[\tilde D_0^{x_k}\circ\cdots\circ\tilde D_0^{x_3}(C_{x_2})]\cup
\cdots\cup
[\tilde D_0^{x_k}(C_{x_{k-1}})]\cup
[C_{x_k}]$$
which is compact as required.
\qed
\enddemo

\proclaim{Proposition 6.7}
Suppose $X$ is a locally compact separable metric space, $A$ is weakly
locally contractible, compact space,
$A\times\Delta^k$ is sliced forward tame in $X$ with respect to $\pi$,
and $(X,A\times\Delta^k)$ has finitely dominated local holinks.
For every open cover ${\cal U}$ of
$A\times\br\times\Delta^k$,
there exists a strong f\.p\. $\cal U$--homotopy equivalence near
$A\times\Delta^k$
$(\bar f,\bar g,X_1',X_2'):X\to\ocyl(q)$.
\endproclaim

\demo{Proof} (cf\. \cite{28\rm, 9.14}.)
Let $Y,f,g,F,G$ be as in Proposition 6.5.
By Lemma 6.6 there exist a compact subspace $K\subseteq
\ho_\pi(X,A\times\Delta^k)$ and a f\.p\. homotopy 
$D:\ho_\pi(X,A\times\Delta^k)\times I\to\ho_\pi(X,A\times\Delta^k)$
such that
$D_0(\ho_\pi(X,A\times\Delta^k))\subseteq K$
and $D_1=\id.$
Define $\hat D:\ocyl(q)\times I\to\ocyl(q)$
by
$$\hat D_s=\cases
D_s\times\id_\br, 
& \text{on $\ho_\pi(X,A\times\Delta^k)\times\br$}\\
\id,
& \text{on $A\times\Delta^k$.}
\endcases$$
Define $g':\ocyl(q)\to Y$
by $g'=g\hat D_0$.
Define $F':Y\times I\to X$ by
$$F'_s=\cases
ig\hat D_{2s}f,
& \text{if $0\leq s\leq 1/2$}\\
F_{2s-1},
& \text{if $1/2\leq s\leq 1$.}
\endcases$$
Note that $F':ig'f\simeq i$.
Define $G':\ocyl(q)\times I\to\ocyl(q)$
by
$$G'_s=\cases
G_{2s}\hat D_0,
& \text{if $0\leq s\leq 1/2$}\\
\hat D_{2s-1},
& \text{if $1/2\leq s\leq 1$.}
\endcases$$
Note that $G':fg'\simeq\id$.
As in  \cite{28\rm, 9.14} it is possible to choose a homeomorphism
$\gamma:\br\to\br$ with
$\gamma=\id$ on $(-\infty,0]$ inducing a homeomorphism
$\bar\gamma:\ocyl(q)\to\ocyl(q)$
such that
$\bar f=\bar\gamma f$ is the desired equivalence with inverse
$\bar g=\bar\gamma^{-1}g'$.
($Q$ plays the role of $p$ in  Definition 6.1.)
\qed
\enddemo

\proclaim{Theorem 6.8}
Suppose  $X$ is a locally compact separable metric space,
$(X,A\times\Delta^k)$ is a sliced homotopically stratified pair with
finitely dominated local holinks, $A$ is a compact ANR and 
$p:X\to A\times (-\infty,+\infty]\times\Delta^k$ is
a f\.p\. proper map with $p|:A\times\Delta^k = p^{-1}(A\times\{+\infty\}
\times\Delta^k)\to A\times\{+\infty\}\times\Delta^k$ the identity.
Then for every 
open cover $\cal U$ of $A\times\br\times\Delta^k$,
there exist a compact neighborhood $N$ of $A\times\Delta^k$ in $X$
and  a f\.p\. strict homotopy $p|N\simeq p':N\to
A\times (-\infty,+\infty]\times\Delta^k$ \rel $A\times\Delta^k$
such that $p'$ is a sliced $\cal U$--fibration over
$A\times (0,+\infty)\times\Delta^k$ and 
$(p')^{-1}(A\times (0,+\infty)\times\Delta^k)$ is open in $X$.
\endproclaim

\demo{Proof} 
Given the open cover $\cal U$ choose an open cover $\cal V$ such that
$\st^2({\cal V})$ refines $\cal U$.
According to Proposition 6.7 there exists a strong  f\.p\. $\cal V$--homotopy 
equivalence near $A\times\Delta^k$
$(\bar f,\bar g,X_1',X_2'):X\to\ocyl(q)$ such that
$X_1'$ is compact. Let
$p''=Q\bar f:X_1'\to A\times(-\infty,+\infty]\times\Delta^k$.
Since $(X,A\times\Delta^k)$ is sliced forward tame there exist a
compact neighborhood $N$ of $A\times\Delta^k$ in $X$ and a f\.p\.
nearly strict deformation $r$ of $N$ into $A\times\Delta^k$ with
$N\subseteq X_1'$ and
$r:N\times I\to X_1'$.
We show that there exists a f\.p\. strict homotopy
$H:p|N\simeq p''|N$ \rel $A\times\Delta^k$ as follows.
Let $\pi_1:A\times (-\infty,+\infty]\times\Delta^k\to A\times\Delta^k$
and $\pi_2:A\times (-\infty,+\infty]\times\Delta^k\to (-\infty,+\infty]$
denote the projections.
Define $H:N\times I\to A\times (-\infty,+\infty]\times\Delta^k$ by
$$\pi_1H(x,t) =\cases
pr(x,2t), & \text{if $0\leq t\leq 1/2$}\\
p''r(x,2-2t), & \text{if $1/2\leq t\leq 1$}
\endcases$$
and $\pi_2H(x,t) = (1-t)\pi_2p(x)+t\pi_2p''(x)$.
According to Proposition 6.3 there exists an $m>0$
such that $p''$ is a sliced $\cal U$--fibration
over $(A\times (m,+\infty)\times\Delta^k)$. We may assume that
$(p'')^{-1}(A\times (m,+\infty)\times\Delta^k)\subseteq N$.
We conclude the proof  by
defining an isotopy 
$G:A\times(-\infty,+\infty]\times\Delta^k\times I\to A\times (-\infty,+\infty]
\times\Delta^k$ by
$G(x,s,t,u) = (x,s-um,t)$ and
setting $p'=G_1p''$. Since $G_0=\id$, $A\times (0,+\infty]\times\Delta^k
= G_1(A\times (m,+\infty)\times\Delta^k)$ and $G_1$ is an isometry, 
it follows that
$G_up'':p''\simeq p'$, $0\leq u\leq 1$, 
and $p'$ is the desired map.
\qed
\enddemo


\head 7. Higher classification of stratified neighborhoods \endhead
Throughout this section $B$ will denote a fixed closed manifold.
We will prove Theorem 2.3, the main result of this paper, which
classifies families of neighborhoods of $B$ in stratified pairs
with $B$ as the lower stratum.
This higher classification is given in terms of families of manifold
approximate fibrations over $B\times\br$.
In fact, Theorem 2.3 asserts that the teardrop construction defines a 
homotopy equivalence between the moduli space of manifold approximate fibrations
over $B\times\br$ and the moduli space of stratified neighborhoods of $B$.
There are two aspects of the proof: existence and uniqueness.
Existence essentially means that the simplicial map between moduli
spaces is surjective on homotopy groups,
whereas uniqueness means that the
map is injective on homotopy groups.
The actual proof combines both aspects by verifying that the map is
`relatively surjective' on homotopy groups. However, the two aspects are
evident in the lead-up to the proof.

The existence problem involves showing that a family (parametrized by
$\Delta^k$) of stratified neighborhoods of $B$ is given by the
teardrop of a family of manifold approximate fibrations over
$B\times\br$.
The precise statement is Proposition 7.2. It is proved by first
appealing to Theorem 6.8 which establishes that such a family of 
neighborhoods is given by the teardrop of a family of
$\cal U$--fibrations over $B\times\br$ where $\cal U$ is an
arbitrarily small open cover of $B\times\br$. Then we use sucking phenomena
for manifold approximate fibrations, which says that if $\cal U$ is 
sufficiently fine then a $\cal U$--fibration deforms to a manifold
approximate fibration.
Sucking phenomena for approximate fibrations were
first discovered by Chapman \cite{7}, \cite{8}, but the family 
version which we require appears in \cite{24}. 
The technical version of sucking which we require is stated in
Proposition 7.1. We point out below that Proposition 7.2 together
with the material from \S4 suffices to give a proof of Theorem 2.1
(Teardrop Neighborhood Existence) even though it also follows from
Theorem 2.3.

Just as the existence aspect is based on a fundamental phenomenon
of manifold approximate fibrations, the uniqueness aspect is 
based on another such phenomenon of manifold approximate fibrations:
two families of close manifold approximate fibrations can be connected
by a close family of manifold approximate fibrations (parametrized
by $\Delta^k$). In other words, the moduli space of manifold 
approximate fibrations is locally $k$--connected for each $k\geq 0$.
This phenomenon was observed in \cite{24}.
Lemma 7.3 contains an elementary argument which shows how we get into
a situation of having two close families of manifold approximate
fibrations. Proposition 7.4 is the technical version of the local
connectivity result which we require and Proposition 7.5 sets the stage for
how it is used in the proof of the classification theorem.

We begin by quoting the version of the sucking phenomena
which we will use.

\proclaim{Proposition 7.1 (Sucking)}
Let  $n\geq 5$ and $k\geq 0$.
For every open cover $\cal U$ of $B\times\br\times\Delta^k$
there exists an open cover $\cal V$ of $B\times\br\times\Delta^k$
such that if $M$ is an $n$--manifold \rom(without boundary\rom), 
$N\subseteq M\times\Delta^k$ is a closed subset, 
$j:N\to B\times\br\times\Delta^k$ is a f\.p\. proper map such that
$j$ is a sliced $\cal V$--fibration over
$B\times (0,+\infty)\times\Delta^k$, and
$j^{-1}(B\times (0,+\infty)\times\Delta^k)$ is an open
subspace of $M\times\Delta^k$,
then $j$ is f\.p\. properly $\cal U$--homotopic 
rel $j^{-1}(B\times (-\infty,0]\times\Delta^k)$ to a map
$j':N\to B\times\br\times\Delta^k$ with $j'$ a sliced 
approximate fibration over $B\times (1,+\infty)\times\Delta^k$.
\endproclaim

\demo{Proof} See \cite{24}, \cite{29\rm, \S13}.
\qed
\enddemo

In the next result we combine the homotopy information of the previous 
section (Theorem 6.8) with the sucking result (Proposition 7.1) to
prove the existence of manifold approximate fibration teardrop
structure for manifold stratified neighborhoods.

\proclaim{Proposition 7.2}
If $n\geq 5$ and $\pi:(X,B\times\Delta^k)\to\Delta^k$ is a 
$k$--simplex of $\SN^n(B)$, then there exists a compact neighborhood
$\widehat{N}$ of $B\times\Delta^k$ in $X$ and a f\.p\. proper
strict map
$$\widehat{p}:(\widehat{N}, B\times\Delta^k)\to (B\times (-\infty,+\infty]
\times\Delta^k, B\times\{+\infty\}\times\Delta^k)
\qquad\rel B\times\Delta^k$$
such that $\widehat{p}$ is a sliced approximate fibration over 
$B\times (1,+\infty)\times\Delta^k$.
\endproclaim

\demo{Proof} 
Choose  an open cover $\cal U$ of $B\times\br\times\Delta^k$
such that 
$$\lub\{\diam(U)~|~U\in{\cal U}, U\cap(B\times [m,+\infty)\times
\Delta^k\not=\emptyset\}\to 0 ~~\text{as}~~ m\to\infty.$$
Let $\cal V$ be an open cover of $B\times\br\times\Delta^k$ given by
Proposition 7.1 which depends on $\cal U$.
Since $B\times\Delta^k$ is sliced forward tame in $X$, it follows 
that there exist a compact neighborhood $N_0$
of $B\times\Delta^k$ in $X$ and a f\.p\. retraction
$r:N_0\to B\times\Delta^k$.
We may assume that $N_0$ is contained in a trivial neighborhood of
$B\times\Delta^k$ (in the sense of Definition 5.1).
Let $N= {\inr}(N_0)$ and choose a proper map
$u:N\to (-\infty,+\infty]$ such that
$u^{-1}(+\infty) = B\times\Delta^k$.
Define $p':N\to B\times (-\infty,+\infty]\times\Delta^k$
by $p'(x) = ({\proj}_Br(x), u(x), {\proj}_{\Delta^k}r(x))$.
Note that $p'$ is a f\.p\. proper strict map and \rel $B\times\Delta^k$.
Since $(X,B\times\Delta^k)$ is a sliced manifold stratified pair,
so is $(N,B\times\Delta^k)$  (Proposition 5.2).
Theorem 6.8 implies that there exist a compact neighborhood 
$\widehat{N}$ of $B\times\Delta^k$ in $N$ and a f\.p\. proper
strict homotopy
$$p'|\widehat{N}\simeq p'':\widehat{N}\to B\times (-\infty,+\infty]\times
\Delta^k \qquad\rel B\times\Delta^k$$
such that $p''$ is a sliced $\cal V$--fibration over $B\times (0,+\infty)
\times\Delta^k$ and
$(p'')^{-1}(B\times (0,+\infty)\times\Delta^k)$ is open in $N$ (and hence
open in $X$).
Now Proposition 7.1 and the choice of $\cal V$ imply that there exists a
f\.p\. proper $\cal U$--homotopy 
$$p''|\widehat{N}\setminus(B\times\Delta^k)\simeq p''':\widehat{N}\setminus
(B\times\Delta^k)\to B\times\br\times\Delta^k$$
such that $p'''$ is a sliced approximate fibration over $B\times (1,+\infty)
\times\Delta^k$. (We are in  a product situation as required by Proposition
7.1 because $N_0$ was chosen to be in a trivial neighborhood.)
The defining property of the open cover $\cal U$ implies that the map
$p'''$ extends via the identity on $B\times\Delta^k$ to a map
$$\widehat{p}:\widehat{N}\to B\times (-\infty,+\infty]\times\Delta^k.\qed$$
\enddemo

As mentioned in \S2 we can now give a proof of Theorem 2.1 (on the
existence of teardrop neighborhoods) which avoids some of the 
machinery required for the proof of Theorem 2.3.

\demo{Proof of Theorem 2.1 \rom(Teardrop Neighborhood Existence\rom)}
If $(X,B)$ is a manifold stratified pair with $\dim(X\setminus B)
=n\geq 5$, then $(X,B)$ is a vertex of $\SN^n(B)$. It follows from
Proposition 7.2 that $B$ has a neighborhood in $X$ which is the teardrop
of a manifold approximate fibration.
The converse follows from Corollary 4.11.
\qed
\enddemo

We are now ready to begin the uniqueness aspects of the main result.
The first lemma shows how to modify two teardrop collapse maps so that
they become close near the lower stratum.

\proclaim{Lemma 7.3}
Suppose $B,K$ are compact metric spaces, $X$ is a locally compact metric
space containing $B\times K$ with a map $\pi:X\to K$ such that
$\pi|:B\times K\to K$ is projection.
Suppose  $p,q:(X,B\times K)\to (B\times (-\infty,+\infty]\times K,
B\times\{+\infty\}\times K)$ are two fibre preserving 
\rom(with respect to $\pi$\rom)
strict maps which are the identity on $B\times K$
and proper over $B\times (0,+\infty)\times K$. 
For every open 
cover $\cal V$ of $B\times\br\times K$ there exists a f\.p\. strict
isotopy $H:B\times (-\infty,+\infty]\times K\times I\to
B\times (-\infty,+\infty]\times K\times I$
rel $(B\times (-\infty,0]\times K)\cup
(B\times\{+\infty\}\times K)$
such that
$p' = H_1p$ and  $q' = H_1q$
are $\cal V$--close over $B\times (1,+\infty)\times K$
\rom(meaning if $x\in (p')^{-1}(B\times (1,+\infty)\times K)\cup
(q')^{-1}(B\times (1,+\infty)\times K)$, then there exists $V\in\cal V$
such that $p'(x),q'(x)\in V$\rom).
\endproclaim

\demo{Proof}
Assume $B\times K$ has a fixed metric, $\br$ has the standard metric
and $B\times\br\times K$ has the product metric.
For each $n=-1,0,1,2,\dots$ let $\epsilon_n >0$ be a Lebesque number for
the open cover $\{ V\cap (B\times [n,n+1]\times K)~|~ V\in{\cal V}\}$
of $B\times [n,n+1]\times K$.
We may assume that $\epsilon_{-1}<\epsilon_0<\epsilon_1<\dots$.
Using the properness of $p,q$ 
(over $B\times (0,+\infty)\times K$)
and the fact that $p,q$ are the identity on
$B\times K$, construct  (by induction) a sequence
$0 < t_{-1} < t_0 < t_ 1 <\dots$ such that $t_n\to\infty$ as $n\to\infty$,
$p,q$ are $(\epsilon_n/3)$--close over $B\times [t_n,+\infty]\times K$, and
if $x\in p^{-1}(B\times [t_n,t_{n+1}]\times K)\cup
q^{-1}(B\times [t_n,t_{n+1}]\times K)$, then $p(x),q(x)\in B\times
[t_{n-1},t_{n+2}]\times K$
for each $n=0,1,2,\dots$.
Also construct a sequence $0=y_0 < y_1 < y_2 < \dots$ refining
$\{0, 1,2,\dots\}$ such that
$y_n \geq n$ and if $n \leq y_k\leq n+1$, then $y_{k+1}-y_k < 
\epsilon_{n+1}/3$.
Define a homeomorphism $h':(-\infty,+\infty]\to (-\infty,+\infty]$
so that for each $n=0,1,2,\dots$ $h'(t_n)=y_n$, $h'$ is linear on
$[t_n,t_{n+1}]$ and is the identity on $(-\infty,0]$.
Define $h=\id_B\times h'\times\id_K:B\times (-\infty,+\infty]\times K
\to B\times (-\infty,+\infty]\times K$.
The natural isotopy $\id_{(-\infty,+\infty]}\simeq h'$ induces 
an isotopy $H:\id_{B\times(-\infty,+\infty]\times K}\simeq h=H_1$
and
one checks that $p'=H_1p$ and $q'=H_1q$ satisfy the conclusions.
\qed
\enddemo

The next result formulates the version of local connectivity
for families of manifold approximate fibrations which we require.
Then Proposition 7.5 applies it in the situation which will arise
in the proof of the main result.

\proclaim{Proposition 7.4}
Suppose that $n\geq 5$ and $K$ is a compact polyhedron. For every 
open cover $\cal U$ of $B\times\br\times K$ there exists an
open cover $\cal V$ of $B\times\br\times K$
such that if $\pi:M\to K$ is a fibre bundle projection with
$n$--manifold fibres \rom(without boundary\rom), 
$N\subseteq M$ is a closed subset,
$p_1,p_2:N\to B\times\br\times K$ are two f\.p\. proper maps 
which are $\cal V$--close over $B\times (0,+\infty)\times K$
and sliced approximate fibrations over $B\times (0,+\infty)\times K$,
and $p_i^{-1}(B\times (0,+\infty)\times K)$ is open in $M$ for $i=1,2$,
then there exists a f\.p\. proper $\cal U$--homotopy
$F:p_1\simeq p_2$ such that
$F_s:N\to B\times\br\times K$ is a sliced approximate fibration over 
$B\times (1,+\infty)\times K$ for each $0\leq s\leq 1$.
\endproclaim

\demo{Proof}
This just involves minor modifications in the arguments of
\cite{24} used to prove that spaces of manifold approximate fibrations
are locally $k$--connected for each $k\geq 0$.
\qed
\enddemo

\proclaim{Proposition 7.5}
Suppose $K$ is a compact polyhedron and
$\pi:(Y,B\times K)\to B\times K$ is a sliced manifold stratified
pair with $\dim\pi^{-1}(u) = n\geq 5$ for $u\in K$ for which there is
a f\.p\. proper strict map
$$p:(Y,B\times K)\to (B\times (-\infty,+\infty]\times K, B\times\{+\infty\}
\times K)
~~{\rel}~~ B\times K$$
which is a sliced manifold approximate fibration over $B\times\br\times K$.
Suppose $t\in\br$ and $\widehat{Y}$ is an open neighborhood of $B\times K$
in $Y$ for which there is
a f\.p\. proper strict map
$$\widehat{p}:(\widehat{Y},B\times K)\to 
(B\times (t,+\infty]\times K, B\times\{+\infty\}\times K)
~~{\rel}~~ B\times K$$
which is a sliced manifold approximate fibration over 
$B\times (t,+\infty)\times K$.
Then there exist $t_2 > t$, a compact neighborhood $X$ of $B\times K$
in $Y$ with $X\subseteq\widehat{Y}$ and a f\.p\. strict homotopy
$$F:p|X\simeq \widehat{p}|X:X\to B\times (-\infty,+\infty]
\times K ~~{\rel}~~ B\times K$$
which is proper over $B\times (t_2,+\infty]\times K$ 
and such that 
$F_s:X\to B\times (-\infty,+\infty]\times K$
is a sliced manifold approximate fibration over 
$B\times (t_2,+\infty)\times K$ for each $0\leq s\leq 1$.
\endproclaim

\demo{Proof}
Choose  an open cover $\cal U$ of $B\times\br\times K$
such that 
$$\lub\{\diam(U)~|~U\in{\cal U}, U\cap(B\times [m,+\infty)\times
\Delta^k\not=\emptyset\}\to 0 ~~\text{as}~~ m\to\infty.$$

Let $\cal V$ be the open cover of $B\times\br\times K$ given by 
Proposition 7.4 which depends on $\cal U$.
Let $W$ be a locally trivial neighborhood of $B\times K$ in $Y$
(in the sense of Definition 5.1)
and assume that $W\subseteq\widehat{Y}$.
Choose $t_0\geq t$ such that 
$$p^{-1}(B\times [t_0,+\infty]\times K)\cup 
\widehat{p}^{-1}(B\times [t_0,+\infty]
\times K)\subseteq W.$$
Let 
$$X = p^{-1}(B\times [t_0,+\infty]\times K)\cap 
\widehat{p}^{-1}(B\times [t_0,+\infty]
\times K).$$
Choose $t_1 >t_0$ such that
$$p^{-1}(B\times [t_1,+\infty]\times K)\cup \widehat{p}^{-1}(B\times (t_1,+\infty]
\times K)\subseteq X$$
and note that
$p|,\widehat{p}|:X\to B\times (t_0,+\infty]\times K$ are proper 
over $B\times (t_1,+\infty]\times K$ and sliced approximate fibrations
over $B\times (t_1,+\infty)\times K$.
Let $t_2 = t_1+1$.
Lemma 7.3 can be applied to yield a f\.p\. strict isotopy
$$\multline 
H:B\times (-\infty,+\infty]\times K\times I\to
B\times (-\infty,+\infty]\times K\times I \\
\rel (B\times (-\infty,t_1]\times K)\cup (B\times \{+\infty\}\times K)
\endmultline$$ 
such that
$p' = H_1p|X$ and $q' = H_1\widehat{p}|X$ are $\cal V$--close over
$B\times (t_2,+\infty]\times K$.
Because $H$ is \rel $B\times (-\infty,t_1]\times K$,
$p'$ and $q'$ are sliced approximate fibrations over
$B\times (t_1,+\infty]\times K$.
Proposition 7.4 can be applied to yield a f\.p\.
$\cal U$--homotopy 
$F:p'|\simeq q'|:X\setminus(B\times K)\to B\times\br\times K$
such that
$F_s:X\setminus (B\times K)\to B\times\br\times K$ is a sliced
approximate fibration over
$B\times (t_2,+\infty)\times K$ for each 
$0\leq s\leq 1$.
The choice of the open cover $\cal U$ implies that $F$ extends via
the identity $B\times K\to B\times\{+\infty\}\times K$ to a homotopy
(also denoted $F$)
$F:p'\simeq q':X\to B\times (-\infty,+\infty]\times K$.
\qed
\enddemo

We need one more lemma before proving the main result.

\proclaim{Lemma 7.6}
If $n\geq 5$ and $t\in\br$, then the restriction
$\rho:\maf^n(B\times\br)\to\maf^n(B\times (t,+\infty))$
is a homotopy equivalence.
\endproclaim

\demo{Proof}
First observe that the techniques of \cite{29\rm, \S3} show that $\rho$ is in
fact a simplicial map. There are a couple of approaches to proving 
that $\rho$
is a homotopy equivalence. One is to use geometric techniques as
presented in \cite{29\rm, \S14} in proving uniqueness of fibre germs.
The other is to use the Manifold Approximate Fibration Classification
Theorem \cite{29}, \cite{30} and observe that restriction induces a homotopy
equivalence of the classifying spaces.
\qed
\enddemo

Let  $n\ge 5$.
We prove the main theorem by showing that $\Psi : \maf^n
(B\times \br)\to \SN^n(B)$ (as constructed in
\S5) is a  homotopy equivalence.  Since both these
simplicial sets satisfy the Kan condition, it suffices
to show that $\Psi $ induces an isomorphism on
homotopy groups (including $\pi_0$).  To accomplish
this suppose that we are given the following
set-up.

\remark{Data \rom{7.7}} Suppose $k\ge 0$.
\roster
\item Let $\pi:(X,B\times\Delta^k)\to\Delta^k$
be a $k$-simplex of $\SN^n(B)$.

\item Let $p:M\to B\times\br\times\partial\Delta^k$
be a union of $(k+1)\ (k-1)$-simplices of $\maf^n(B\times\br)$.

\item Suppose for each $i = 0,\cdots k$, the
$(k-1)$-simplex $\pi|:(\pi^{-1}(\partial_i\Delta^k),B\times\partial_i
\Delta^k)\to\partial_i\Delta^k$
of $\SN^n(B)$ is the image under $\Psi $ of the $(k-1)$-simplex
$p|:p^{-1}(B\times\br\times\partial_i\Delta^k)\to 
B\times\br\times\partial_i\Delta^k$
of $\maf^n(B\times\br)$ so that
$M=\pi^{-1}(\partial\Delta^k)\setminus(B\times\partial\Delta^k)$.
\endroster
Note that if $k= 0$, then only item (1) is meaningful.
\endremark

\proclaim{Theorem 7.8}  Given Data \rom{7.7}, there is a $k$-simplex
$\widetilde{p}:\widetilde{M}\to B\times\br\times\Delta^k$
of $\maf^n(B\times \br)$ which equals $p$ over $B\times
\br\times \partial \Delta ^k$ and whose image under
$\Psi $ is homotopic in $\SN^n(B)$ to $\pi $ \rel $\partial$.
Hence, $\Psi : \maf^n
(B\times \br)\to \SN^n(B) $ induces an isomorphism
on homotopy groups and is a homotopy equivalence.
\endproclaim

\demo{Proof}
According to Proposition 7.2,
there exists a compact neighborhood
$\widehat{N}$ of $B\times\Delta^k$ in $X$ and a f\.p\. proper
strict map
$$\widehat{p}:(\widehat{N}, B\times\Delta^k)\to (B\times (-\infty,+\infty]
\times\Delta^k, B\times\{+\infty\}\times\Delta^k)
\qquad\rel B\times\Delta^k$$
such that $\widehat{p}$ is a sliced approximate fibration over 
$B\times (1,+\infty)\times\Delta^k$.
Choose $t\geq 1$ such that 
$\widehat{p}^{-1}(B\times (t,+\infty]\times\Delta^k)$ is open in $X$.
Let $Y=\partial X=\pi^{-1}(\partial\Delta^k)$ which by assumption is the
teardrop  $M\cup_p (B\times\partial\Delta^k)$.
Extend $p:M\to B\times\br\times\partial\Delta^k$
via the identity $B\times\partial\Delta^k\to B\times\{+\infty\}\times\partial
\Delta^k$ to 
$p_+:Y\to B\times B\times (-\infty,+\infty]\times\partial\Delta^k$
which is continuous since it is the teardrop collapse.
Let $\widehat{Y} = \widehat{p}^{-1}(B\times (t,+\infty]\times\partial\Delta^k)$.
Since $\widehat{Y}$ is open in $Y$, it follows that
$\widehat{p}|:\widehat{Y}\to B\times (t,+\infty]\times\partial\Delta^k$
is a sliced manifold approximate fibration over $B\times (t,+\infty)\times
\partial\Delta^k$.
It follows from Proposition 7.5 applied with $K=\partial\Delta^k$ that
there exist $t_2 > t$, a compact neighborhood $\widetilde{Y}$ of 
$B\times\partial\Delta^k$ in $Y$ with $\widetilde{Y}\subseteq
\widehat{Y}$, and a f\.p\. strict homotopy
$$F:p_+|\widetilde{Y}\simeq\widehat{p}|\widetilde{Y}:
\widetilde{Y}\to B\times (-\infty,+\infty]\times\partial\Delta^k$$
which is  proper over $B\times (t_2,+\infty]\times\partial\Delta^k$
and such that 
$F_s:\widetilde{Y}\to B\times (-\infty,+\infty]\times\partial\Delta^k$
is a sliced manifold approximate fibration over $B\times (t_2,+\infty)
\times\partial\Delta^k$ for each $0\leq s\leq 1$.
Consider $F$ as a  map 
$F:\widetilde{Y}\times I\to B\times (-\infty,+\infty]\times\partial\Delta^k
\times I$.
Choose $t_3\geq t_2$ such that 
$F^{-1}(B\times (t_3,+\infty]\times\partial\Delta^k\times I)$ is open
in $Y\times I$ and let
$W=F^{-1}(B\times (t_3,+\infty)\times\partial\Delta^k\times I)$.
Since the composition 
$$W @>F>> B\times (t_3,+\infty)
@>{\proj}>> \partial\Delta^k\times I$$
is a submersion and
$F|:W\to B\times (t_3,+\infty)\times\partial\Delta^k\times I$
is a sliced (over $\partial\Delta^k\times I$)
manifold approximate fibration, it follows from
\cite{25\rm, Lemma 4.1} that
$W\to\partial\Delta^k\times I$ is a fibre bundle projection.
Let $W_0 = p^{-1}(B\times (t_3,+\infty)\times\partial\Delta^k)$ and
$W_1= \widehat{p}^{-1}(B\times (t_3,+\infty)\times\partial\Delta^k)$.
It follows that
$F|W$ may be thought of as a homotopy in $\maf(B\times (t_3,+\infty))$
from
$p|:W_0\to B\times (t_3,+\infty)\times\partial\Delta^k$ 
to
$\widehat{p}|:W_1\to B\times (t_3,+\infty)\times\partial\Delta^k$.

Now consider the open subspace $\widehat{X} =
\widehat{p}^{-1}(B\times (t_3,+\infty]\times\Delta^k)$ of
$X$ and let 
$\widehat{M} = \widehat{X}\setminus (B\times\Delta^k) =
\widehat{p}^{-1}(B\times (t_3,+\infty)\times\Delta^k)$.
Since $\widehat{p}|:\widehat{X}\to B\times (t_3,+\infty]\times\Delta^k$
is a sliced manifold approximate fibration over $B\times (t_3,+\infty)\times
\Delta^k$, it follows using
\cite{25\rm, Lemma 4.1} again that
$\widehat{p}:\widehat{M}\to B\times (t_3,+\infty)\times\Delta^k$ is
a $k$--simplex of $\maf(B\times (t_3,+\infty))$.
Its boundary is
$\widehat{p}|= F_1|:W_1\to B\times (t_3,+\infty)\times\partial\Delta^k$.

Let $\rho:\maf(B\times\br)\to \maf(B\times (t_3,+\infty))$ be the simplicial
map induced by restriction. It is a homotopy equivalence by Lemma 7.6.
Define a simplicial map
$\Psi':\maf(B\times (t_3,+\infty))\to \SN(B)$ induced by the 
teardrop construction in analogy
to the map $\Psi:\maf(B\times\br)\to\SN(B)$.
In fact, if $q:Q\to B\times\br\times\Delta^k$ is
a $k$--simplex of $\maf(B\times\br)$,
then $\Psi'\rho(q)= q^{-1}(B\times (t_3,+\infty)\times\Delta^k)\cup_q
(B\times\Delta^k)$ is an open subspace of $\Psi(q)=
Q\cup_q(B\times\Delta^k)$
and the mapping cylinder of the inclusion induces a homotopy in
$\SN(B)$ from $\Psi'\rho(q)$ to
$\Psi(q)$ (see \S5).
In this way we construct a homotopy
$${\cal CYL}:\Psi'\rho\simeq \Psi:\maf(B\times\br)\to\SN(B).$$
Use the homotopy $F|W$ and a collar of $\partial\Delta^k$ in
$\Delta^k$ to enlarge the $k$--simplex
$\widehat{p}|:\widehat{M}\to B\times (t_3,+\infty)\times\Delta^k$
of $\maf(B\times (t_3,+\infty))$ to a $k$--simplex
$p^\ast:M^\ast\to B\times (t_3,+\infty)\times\Delta^k$ of
$\maf(B\times (t_3,+\infty))$ so that
$\partial p^\ast$ is $\rho(p)$.
Note that $F$ is a homotopy in $\maf(B\times (t_3,+\infty))$
from $\rho(p) = F_0|W_0$ to $\partial (\widehat{p}|\widehat{M}) = F_1|W_1$.
Note that since $\Psi'(\widehat{p}|\widehat{M})$ is an open
subspace of $X$, the mapping cylinder construction induces a 
homotopy
${\cal CYL}:\Psi'(\widehat{p}|\widehat{M})\simeq X$ in $\SN(B)$.
Note also that
since each 
$F^{-1}(B\times (t_3,+\infty]\times\partial\Delta^k\times\{ s\})$
is an open subspace op $\partial X$, the mapping cylinder construction
induces an extension of the
homotopy
${\cal CYL}:\Psi'(\widehat{p}|\widehat{M})\simeq X$ 
to a homotopy
${\cal CYL}:\Psi'(p^\ast)\simeq X$.

The situation now is that we have a $k$--simplex
$p^\ast$ of $\maf(B\times (t_3,+\infty))$ such that
$\rho(p)=\partial p^\ast$ and the mapping cylinder construction induces a
homotopy
${\cal CYL}:\Psi'(p^\ast)\simeq X$.
Since $\rho:\maf(B\times\br)\to\maf(B\times (t_3,+\infty))$ is a homotopy
equivalence, there exists a $k$--simplex
$\widetilde{p}$ of $\maf(B\times\br)$ 
such that $\partial\widetilde{p} = p$ and a homotopy
$G:\rho(\widetilde{p})\simeq p^\ast$ \rel $\partial\rho(\widetilde{p})
=\partial p^\ast$.
Thus $\Psi'(G)$ is a homotopy in $\SN(B)$ from
$\psi'\rho(\widetilde{p})$ to $\Psi'(p^\ast)$ \rel $\partial$.
This homotopy taken together with the homotopy
${\cal CYL}:\Psi'(p^\ast)\simeq X$,
yields a homotopy
$H:\Psi'\rho(\widetilde{p})\simeq X$ in $\SN(B)$
which restricts to
${\cal CYL}:\partial\Psi'\rho(\widetilde{p})\simeq\partial X$.
On the other hand, we have already observed that there is a homotopy
${\cal CYL}:\Psi'\rho(\widetilde{p})\simeq\Psi(\widetilde{p})$.
The concatenation 
$\Psi(\widetilde{p})\simeq\Psi'\rho(\widetilde{p})\simeq X$,
together with the fact that the two homotopies restrict to inverses
on the boundary, implies that there exists a homotopy
$\Psi(\widetilde{p})\simeq X$ \rel $\partial$.
\qed
\enddemo


\head 8. Examples of exotic stratifications \endhead
In this section we use the classification of neighborhood germs
to construct  examples of manifold stratified pairs
in which the lower stratum does not have a neighborhood given by the
mapping cylinder of a fibre bundle.
Moreover, we construct examples in which this phenomenon
persists under euclidean stabilization.

\proclaim{Theorem 8.1}
For every  integer $m\geq 6$ there exists a locally conelike
manifold stratified pair $(X,S^1)$ with $\dim (X\setminus S^1) = m$
such that $S^1$ has a manifold approximate fibration
mapping cylinder neighborhood in $X$,
but for each $i\geq 0$ $S^1\times\br^i$ does not have a 
fibre bundle mapping cylinder
neighborhood in $X\times\br^i$.
In fact, $S^1\times\br^i$ does not have a block bundle mapping cylinder
neighborhood in $X\times\br^i$.
\endproclaim

For the remainder of this section, let $F$ denote a closed connected
manifold of dimension $n$.
Let $\topbfri$ denote the simplicial set of {\it bounded homeomorphisms\/}
on $F\times\br^i$ so that a $k$--simplex of $\topbfri$ consists of
a homeomorphism $h:F\times\br^i\times\Delta^k\to F\times\br^i\times
\Delta^k$ such that $h$ is fibre preserving over $\Delta^k$ and 
bounded in the $\br^i$--direction. This latter condition means
there exists a constant $c>0$ such that $p_2h$ is $c$--close to $p_2$
where $p_2:F\times\br^i\times\Delta^k\to\br^i$ is projection.

Let $\cbfri$ denote the simplicial set of {\it bounded concordances\/}
on $F\times\br^i$ so that a $k$--simplex of $\cbfri$ consists of a
homeomorphism $h:F\times\br^i\times [0,1]\times\Delta^k\to
F\times\br^i\times [0,1]\times\Delta^k$ such that
$h$ is fibre preserving over $\Delta^k$, $h|:F\times\br^i\times\{0\}
\times\Delta^k\to F\times\br^i\times\{0\}
\times\Delta^k$ is the identity, and $h$ is bounded over $\br^i$.

A bounded concordance on $F\times\br^i$ induces a bounded homeomorphism
on $F\times\br^i$ by restricting the concordance to $F\times\br^i\times\{1\}$.
This defines a simplicial map
$$\rho:\cbfri\to\topbfri$$
by setting $\rho(h) = h|:F\times\br^i\times\{1\}\times\Delta^k
=F\times\br^i\times\Delta^k\to
F\times\br^i\times\{1\}\times\Delta^k = F\times\br^i\times\Delta^k.$

Euclidean stabilization induces a simplicial map
$$\sigma:\topbfri\to\topbf\times\br^{i+1}); \qquad h\mapsto h\times\id_\br$$
and, in particular, a group homomorphism
$\pi_0\topbfri @>\sigma>> \pi_0\topbf\times\br^{i+1})$
for each $i\geq 0$.

\proclaim{Proposition 8.2 (Anderson-Hsiang)} 
There is a homotopy fibration sequence 
$$\cbfri ~@>{\rho}>>~
\topbfri ~@>{\sigma}>>~
\topbf\times\br^{i+1}).$$
In particular, there is a short exact sequence
$$\pi_0\cbfri ~@>{\rho}>>~
\pi_0\topbfri ~@>{\sigma}>>~
\pi_0\topbf\times\br^{i+1}).$$
\endproclaim

\demo{Proof} This is essentially the fibration of Anderson-Hsiang
\cite{3\rm, 9.3}. One must use \cite{2\rm, Thm\. 4} to identify $\cbfri$ with
the fibre in \cite{3}. Similarly one needs a reinterpretation
of $\topbfri$. See  \cite{31\rm, Thm\. 1.2} for an explicit proof.
See also Lashof-Rothenberg \cite{37\rm, \S8}.
\qed
\enddemo

An {\it inertial\/} $h$--cobordism on $F$ is an $h$--cobordism
$(W;\partial_0W,\partial_1W)$ with $\partial_0W=F$ and $\partial_1W$
homeomorphic to $F$.
It is possible to define the simplicial set of $h$--cobordisms on $F$
(e\.g\. Waldhausen \cite{59}) and the simplicial set of inertial $h$--cobordisms
on $F$. However, for this paper we only need the sets of components of these
simplicial sets.
Thus, let $\pi_0\hcobf$ denote the set of equivalence classes
of $h$--cobordisms on $F$ such that 
$(W;\partial_0W,\partial_1W)$ is equivalent to
$(W';\partial_0W',\partial_1W')$ if and only if there exists a 
homeomorphism $H:W\to W'$ such that
$H|:\partial_0W =F\to \partial_0W' =F$ is the identity.
The set $\pi_0\ihcobf$ of inertial $h$--cobordisms on $F$ is the subset
of $\pi_0\hcobf$ consisting of all classes represented by
inertial $h$--cobordisms.

The $s$--cobordism theorem gives a bijection
$$\pi_0\hcobf @>\tau>> \Wh_1(\bz\pi_1 F)$$
provided $n\geq 5$, which sends an $h$--cobordism
$(W;\partial_0W,\partial_1W)$ to the Whitehead torsion
$\tau(W,\partial_0W)$ in $\Wh_1(\bz\pi_1 F)$.
In general, the image of $\pi_0\ihcobf$ in $\Wh_1(\bz\pi_1 F)$
need not be a subgroup (cf\. Hausmann \cite{19}, Ling \cite{39}).

We now recall the well-known `region between' construction
(cf\. Anderson-Hsiang \cite{2\rm,\S8}) which defines a function
$$\beta:\pi_0\topbf\times\br)\longrightarrow\pi_0\ihcobf.$$
If $h:F\times\br\to F\times\br$ is a bounded homeomorphism
representing a class $[h]\in\pi_0\topbf\times\br)$, choose a $L>0$
so large that $h(F\times\{ L\})\subseteq F\times (0,\infty)$.
Let $W= h(F\times (-\infty,L])\setminus F\times(-\infty,0)$,
$\partial_0W=F\times\{0\} = F$, and $\partial_1W=h(F\times\{ L\})$.
Then $(W;\partial_0W,\partial_1W)$ is an inertial $h$--cobordism on
$F$ representing a class $[W]\in\pi_0\ihcobf$. Set
$\beta([h]) = [W]$. 
The function $\beta$ is well-defined by the Isotopy Extension Theorem
of Edwards-Kirby \cite{11}.
One should not confuse $\tau(\beta(h))$ with the torsion of the 
homotopy equivalence
$$h_1:F = F\times\{ 0\}\hookrightarrow F\times\br @>h>> F\times\br @>{\proj}>> F.$$
To see the relationship between these two torsions let
$j:W\to [0,1]$ be any map with $j^{-1}(0) =\partial_0W$ and
$j^{-1}(1) = \partial_1W$. Since $F\times\{ 0\}\hookrightarrow F\times (-\infty,
0]$ and $h(F\times\{ L\})\hookrightarrow h(F\times [L,+\infty))$ are homotopy 
equivalences, so is the inclusion $i:W\hookrightarrow F\times\br$, and there
is a homotopy equivalence of triads
$$\gamma = (\proj\circ~ i)\times j: (W;\partial_0W,\partial_1W)\to
(F\times [0,1]; F\times\{ 0\}, F\times\{ 1\}).$$
Therefore,
$$\tau(h_1) = \tau(\gamma|\partial_1W:\partial_1W\to F\times\{ 1\})=
\tau\beta(h)-(-1)^n\ov{\tau\beta(h)}\in\Wh_1(\bz\pi_1 F),$$
where   $\ov{~\cdot~}$ is induced from the 
standard involution on $\bz\pi_1F$.
Although the composition
$$\pi_0\topbf\times\br) @>\beta>>
\pi_0\ihcobf @>\tau>>
\Wh_1(\bz\pi_1 F)$$
need not be a group homomorphism (cf\. Ling \cite{39}), it is a crossed
homomorphism; i\.e\., $\tau\beta([h\circ k]) =
h_{1\sharp}\tau\beta([k]) +\tau\beta([h])$ for $[h],[k]\in\pi_0\topbf\times\br)$
where $h_{1\sharp}$ is the homomorphism induced by the homotopy equivalence
$h_1:F\to F$.

We will need the following version of the Alexander trick in the proof
of Proposition 8.4 (cf\. Hughes \cite{23\rm, Lemma 6.4}).

\proclaim{Lemma 8.3}
If $h:F\times\br\to F\times\br$ is a bounded homeomorphism such that
$h=\id$ on $F\times (-\infty,0]$, then $h$ is boundedly isotopic to
$\id_{F\times\br}$.
\endproclaim

\demo{Proof}
For $0\leq s<1$ define $\theta_s:\br\to\br$ by
$\theta_s(t) = t -{s\over s-1}$.
Define a bounded isotopy $H: h\simeq \id_{F\times\br}$ by
$$H_s =\cases (\id_F\times\theta_s)^{-1}\circ h\circ(\id_F\times\theta_s),
& \text{if  $0\leq s<1$}\\
\id_{F\times\br}, & \text{if $s=1$.}\qed
\endcases$$
\enddemo

\proclaim{Proposition 8.4}
If $n = \dim F\geq 4$, then the sequence
$$\pi_0{\tp}(F) @>\sigma>>
\pi_0\topbf\times\br) @>\beta>>
\pi_0\ihcobf\longrightarrow 0$$
is exact in the  sense that
$\beta$ maps  the set of 
cosets of $\pi_0\topbf\times\br)/\im(\sigma)$
bijectively onto $\pi_0\ihcobf$\rom; i\.e\.,
\widestnumber\item{(ii)}
\roster
\item"(i)" If $[h_1], [h_2]\in
\pi_0\topbf\times\br)$, then
$\beta([h_1]) = \beta([h_2])$ if and only if
there exists $[g]\in\pi_0\tp(F)$ such that $[h_2^{-1}h_1] =
\sigma([g])$, and
\item"(ii)" $\beta$ is surjective.
\endroster
\endproclaim

\demo{Proof}(i)
Let $h_i:F\times\br\to F\times\br$ be  bounded homeomorphisms for $i=1,2$.
Choose $L>0$ such that
$h_i(F\times\{ L\})\subseteq F\times (0,\infty)$
so that $W_i=h_i(F\times(-\infty,L])\setminus F\times (-\infty,0)$ is an 
$h$--cobordism from $F = F\times\{ 0\}$ to $h_i(F\times\{ L\})$ and
$\beta([h_i]) = [W_i]\in\pi_0\hcobf$ for $i=1,2$.
If $[W_1] = [W_2]$, then there exists a homeomorphism
$H:W_1\to W_2$ such that
$H|:F\times\{ 0\}\to F\times\{ 0\}$ is the identity.
In particular, $Hh_1(F\times\{ 1\}) = h_2(F\times\{ L\})$.
Let $g:F\to F$ be the homeomorphism defined by
$h_2^{-1}Hh_1(x,L) = (g(x),1)\in F\times\{ L\}$ for all $x\in F$.
Extend $H$ via the identity on $F\times (-\infty,0]$
to a homeomorphism
$\widetilde{H}:(F\times (-\infty,0])\cup W_1\to
(F\times (-\infty,0])\cup W_2$.
Define a hybrid homeomorphism
$\widetilde{h}:F\times\br\to F\times\br$ by
$$\widetilde{h}(x,t) =\cases \widetilde{H}h_1(x,t), & \text{if $t\leq L$}\\
                             h_2(g(x),t),   & \text{if $t\geq L$.}
\endcases$$
According to Lemma 8.3 both $\widetilde{h}h_1^{-1}$ and
$\widetilde{h}(g^{-1}\times\id_\br)h_2^{-1}$
are boundedly isotopic to the identity.
Thus $\widetilde{h}$ is boundedly isotopic to $h_1$ and to
$h_2(g\times\id_\br)$ so that
$h_2^{-1}h_1$ is boundedly isotopic to $g\times\id_\br$
showing $[h_2^{-1}h_1] = \sigma([g])$.

Conversely, if $h_2^{-1}h_1$ is boundedly isotopic to $g\times\id_\br$
for some homeomorphism $g:F\to F$, then $h_1$ is boundedly isotopic to
$h_2(g\times\id_\br)$. If $L$ is large enough, then the isotopy
restricts to an isotopy of embeddings carrying $h_1(F\times\{ L\})$
onto $h_2(g(F)\times\{ L\}) = h_2(F\times\{ L\})$ in
$F\times (0,\infty)$. The Isotopy Extension Theorem \cite{11} shows that there
is an isotopy of $F\times\br$ to itself which is the identity on
$F\times (-\infty,0]$ and carries $h_1(F\times\{ L\})$ to $h_2(F\times\{ L\})$. 
In particular, there is a homeomorphism $H:W_1\to W_2$ such that 
$H|F\times\{ 0\}$ is the identity. Hence $[W_1] =[W_2]\in\pi_0\hcobf$.

\noindent
(ii) follows from Ling \cite{39\rm, Prop\. 3.2}
\qed
\enddemo

Anderson and Hsiang \cite{2} calculated the homotopy groups of the simplicial
set of bounded concordances. We will need their calculation of the 
group of components.

\proclaim{Proposition 8.5 (Anderson-Hsiang)}
If $n=\dim F \geq 5$, then there exists a group isomorphism
$$\alpha:\pi_0\cbfri\longrightarrow\cases
\Wh_1(\bz\pi_1 F), & \text{if $i=1$}\\
\widetilde{K}_0(\bz\pi_1 F), & \text{if $i=2$}\\
K_{2-i}(\bz\pi_1 F), & \text{if $i>2$.}
\qed
\endcases$$
\endproclaim

We  need to recall the explicit construction of the isomorphism
when $i=1$, $\alpha:\pi_0\cbf\times\br)\to\Wh_1(\bz\pi_1F)$.
If $h:F\times\br\times [0,1]\to F\times\br\times [0,1]$ is a bounded
concordance representing a class $[h]\in\pi_0\cbf\times\br)$,
choose $L>0$ so large that $h(F\times\{ L\}\times [0,1])\subseteq 
F\times (0,\infty )\times [0,1]$ and let
$W= h(F\times (-\infty,L]\times [0,1])\setminus F\times (-\infty,0)\times
[0,1]$, $\partial_0W= F\times [0,L]\times\{0\}$, and
$\partial_1W = h(F\times (-\infty,L]\times\{1\})\setminus F\times (-\infty,0)
\times\{1\}$.
Then $(W;\partial_0W,\partial_1W)$ is a relative $h$--cobordism. In 
particular, over the boundary of $\partial_0W$,
$W$ restricts to a product $h$--cobordism
$$(F\times\{0\}\times [0,1]\cup h(F\times\{ L\}\times [0,1];
F\times\{0,L\}\times\{0\}, F\times\{0\}\times\{1\}\cup
h(F\times\{ L\}\times\{1\})).$$
Define $\alpha([h])$ to be the Whitehead torsion $\tau (W,\partial_0W)
\in\Wh_1(\bz\pi_1(F\times [0,L])= \Wh_1(\bz\pi_1F)$.

Recall that $\dim F=n$. Define the {\it norm homomorphism}
$$N:\Wh_1(\bz\pi_1F)\to\Wh_1(\bz\pi_1F); \qquad x\mapsto
x+ (-1)^n\ov{x}$$ 
where $\ov{~\cdot~}$ is induced from the 
standard involution on $\bz\pi_1F$.

\proclaim{Proposition 8.6}
If\/ $n=\dim F\geq 5$, then the following diagram commutes\rom:
$$\CD 
\pi_0\cbf\times\br) @>\alpha>>  \Wh_1(\bz\pi_1F) @>N>> \Wh_1(\bz\pi_1F)\\
@V{\rho}VV @. @AA{\tau}A\\
\pi_0\topbf\times\br) @>\beta>> \pi_0\ihcobf  @>{\subseteq}>> \pi_0\hcobf. \\
\endCD$$
\endproclaim

\demo{Proof} 
If $[h]\in\pi_0\cbf\times\br)$ adopt the notation above in the explicit
description of $\alpha$ so that $\alpha([h]) = \tau(W,\partial_0W) = x$.
For $k=0,1$ let $i_k:\partial_kW\to W$ denote the inclusion and
$r_k:W\to\partial_kW$ a strong deformation retraction.
Then $x=\tau(r_0)$. 
Since $(W;\partial_0W,\partial_1W)$ is a relative $h$--cobordism
between $(n+1)$--dimensional manifolds, it follows that
$(r_0i_1)_\ast\tau(r_1) = (r_0i_1)_\ast\tau(W,\partial_1W) =
(-1)^{n+1}\ov{x}$   by the duality theorem of Milnor \cite{43\rm, p\. 394}.
Thus, $\tau(i_1) = (-1)^{n+1}\ov{\tau}(i_0)\in\Wh_1(\bz\pi_1W)$.

Let $j_1:F\times\{0\}\times\{1\}\to W$ and 
$j_2:F\times\{0\}\times\{1\}\to \partial_0W$
denote the inclusions and let
$j_3:F\times\{0\}\times\{1\}\to \partial_0W$
be the map $j_3(z,0,1) = (z,0,0)$. Since $\partial_0W= F\times [0,L]\times
\{0\}$, $\tau(j_3) = 0$.

Since $\tau\beta\rho([h])$ is the Whitehead torsion of
$(\partial_1W,F\times\{0\}\times\{1\})$ in $\Wh_1(\bz\pi_1F)$,
it suffices to show that
$$i_{1\ast}\tau(j_2) = \tau(i_0) = (-1)^n\ov{\tau}(i_0)\in
\Wh_1(\bz\pi_1W).$$
The composition formula gives
$$\tau(j_1) = \tau(i_1j_2) = i_{1\ast}\tau(j_2)+\tau(i_1).$$
Since $j_1\simeq i_0j_3$ and $\tau(j_3)=0$, the composition formula
also gives
$\tau(j_1) = \tau(i_0j_3) = \tau(i_0)$.
Thus
$$i_{1\ast}\tau(j_2) = \tau(i_0)-\tau(i_1) = \tau(i_0)+(-1)^{n+1}\tau(i_0).$$
A similar argument has been used by Siebenmann and Sondow \cite{57\rm, p\. 266}.
\qed
\enddemo

\proclaim{Lemma 8.7}

\noindent
{\rm (i)}
Suppose there is a diagram
$$\CD {} @. A   @. {} \\
         @.  @VV{\sigma_1}V   @. \\
      A' @>{\rho_1}>> B @>\beta>> C\\
         @.  @VV{\sigma_2}V   @.  \\
       {} @.  C' @. {} \\
\endCD$$
such that
\widestnumber\item{(3)}
\roster
\item"(1)" $A,B,A',C'$ are groups \rom(written additively\rom), $C$ is a set, 
and $\sigma_1,\rho_1,\sigma_2$ are group homomorphisms,
\item"(2)" $A @>{\sigma_1}>> B @>{\beta}>> C$ is exact in the sense that 
$\beta$ is surjective, and if $b_1,b_2\in B$ then $\beta(b_1)=\beta(b_2)$
if and only if $b_2 - b_1 = \sigma_1(a)$ for some $a\in A$,
\item"(3)"
$A' @>{\rho_1}>> B @>{\sigma_2}>> C'$ is an exact sequence of groups.
\endroster
If $b\in B$, then 
 $\sigma_2(b)\in\im(\sigma_2\sigma_1:A\to C')$
if and only if
$\beta(b)\in\im(\beta\rho_1:A'\to C)$.

\noindent
{\rm (ii)}
Suppose further that the diagram above is extended to a diagram
$$\CD {} @. A   @. {} \\
         @.  @VV{\sigma_1}V   @. \\
      A' @>{\rho_1}>> B @>\beta>> C @>{\tau}>> W\\
         @.  @VV{\sigma_2}V   @.  \\
       D @>{\rho_2}>> C' @. {} \\
         @.  @VV{\sigma_3}V   @.  \\
      {} @. E   @. {} 
\endCD$$
such that
\roster
\item $A,B,A',C',W,D, E$ are abelain groups and $\rho_2, \sigma_3$ are
group homomorphisms,
\item $\tau:C\to W$ is a set inclusion,
\item there is a $B$-module structure on $W$ which satisifies:
if $b_1, b_2\in B$ and $\sigma_2(b_1)=\sigma_2(b_2)$, then
$b_1w=b_2w$ for all $w\in W$,
\item $\tau\beta:B\to W$ is a crossed homomorphism with respect to 
the $B$-module structure \rom(i.e.\rom, $\tau\beta(b_1+b_2)=
b_1\tau\beta(b_2)+\tau\beta(b_1)$ for all $b_1, b_2\in B$\rom),
\item $\tau\beta(-b)=-\tau\beta(b)$ for all $b\in B$,
\item $N=\tau\beta\rho_1:A'\to W$ is a homomorphism,
\item $D @>{\rho_2}>> C' @>{\sigma_3}>> E$ is an exact sequence of groups.
\endroster
There exists a function $\tilde\beta:\im(\sigma_2)\to
W/\im(N)$ such that
if $b\in B$, then $\sigma_3\sigma_2(b)\in \im(\sigma_3\sigma_2\sigma_1)$
if and only if the class of
$\tau\beta(b)$ in $W/\im(N)$ is in  
$\tilde\beta[\im(\sigma_2)\cap\im(\rho_2)].$
\endproclaim

\demo{Proof}
(i)
Suppose first that $\sigma_2(b)= \sigma_2\sigma_1(a)$ for
some $a\in A$.
Then the exact sequence of groups implies that there exists
$a'\in A'$ such that $\rho_1(a') = b + \sigma_1(-a)$. Thus
$-b + \rho_1(a')=\sigma_1(-a)$ and exactness of the other sequence
implies $\beta(b) = \beta\rho_1(a')$.

\noindent
Conversely, suppose $\beta(b) = \beta\rho_1(a')$ for some $a'\in A'$.
Exactness implies that $b= \sigma_1(a)+\rho_1(a')$ for some $a\in A$.
Thus $\sigma_2(b) = \sigma_2\sigma_1(a)+\sigma_2\rho_1(a')= \sigma_2\sigma_1(a)$.

\noindent (ii)
Define $\tilde\beta:\im(\sigma_2)\to W/\im(N)$ by
$\tilde\beta(x) =\tau\beta\sigma_2^{-1}(x)$. 
In order to show that $\tilde\beta$ is well-defined, suppose that
$\sigma_2(y_1)=\sigma_2(y_2)$ and show that
$\tau\beta(y_1)-\tau\beta(y_2)\in\im(N)$. Since $\ker(\sigma_2) =
\im(\rho_1)$, it follows that $\tau\beta(y_1-y_2)\in\im(N)$.
Now
$\tau\beta(y_1-y_2)-[\tau\beta(y_1)-\tau\beta(y_2)]
= (-y_2)\tau\beta(y_1)\tau\beta(-y_2)-\tau\beta(y_1)+\tau\beta(y_2)
= (-y_2)\tau\beta(y_1)+\tau\beta(-y_1)
=(-y_1)\tau\beta(y_1)+\tau\beta(-y_1)
=\tau\beta(y_1-y_1) = 0$.
Thus, $\tau\beta(y_1-y_2) =\tau\beta(y_1)-\tau\beta(y_2)$ showing
$\tilde\beta$ is well-defined.

\noindent
Suppose that $\sigma_3\sigma_2(b)\in\im(\sigma_3\sigma_2\sigma_1)$,
say $a\in A$ with $\sigma_3\sigma_2b=\sigma_3\sigma_2\sigma_1(a)$.
Then $\sigma_2(b)-\sigma_2\sigma_1(a)\in\ker(\sigma_3=\im\rho_2$,
so let $d\in D$ with $\rho_2(d)=\sigma_2(b)-\sigma_2\sigma_1(a)$.
Thus, $\sigma_2(b)=\rho_2(d)+\sigma_2\sigma_1(a)$ and
$\rho_2(d)\in\im(\sigma_2)\cap\im(\rho_2)$.
It follows that $\tilde\beta\rho_2(d)=\tau\beta\sigma_2^{-1}(\rho_2(d))=
\tau\beta(b-\sigma_1(a))=\tau\beta(b)-\tau\beta(\sigma_1(a))$.
Thus, we will be done by showing that $\tau\beta(\sigma_1(a)\in\im(N)$.
By part (i), this is equivalent to showing that $\sigma_2\sigma_1(a)\in
\im(\sigma_2\sigma_1)$, which is obviously true.

\noindent
Conversely, if the class of
$\tau\beta(b)$ in $W/\im(N)$ is in  
$\tilde\beta[\im(\sigma_2)\cap\im(\rho_2)]$, choose $x\in
\im(\sigma_2)\cap\im(\rho_2)$ such that $\tilde\beta(x) = \tau\beta(b)+
\im(N)$. Thus, there exists $y\in B$ such that $\sigma_2(y)=x$ and
$\tau\beta(b)-\tau\beta(y)\in\im(N)$.
By exactness of $A @>{\sigma_1}>> B @>{\beta}>> C$ there exists $a\in A$
such that $\sigma_1(a)=b-y$, from which it follows that
$\sigma_3\sigma_2\sigma_1(a) = \sigma_3\sigma_2(b)-\sigma_3\sigma_2(y)=
\sigma_3\sigma_2(b)-\sigma_3\sigma_3(x)$. But $x\in\im(\rho_2)=\ker(\sigma_3)$,
so $\sigma_3\sigma_2\sigma_1(a) = \sigma_3\sigma_2(b)$.
\qed
\enddemo

We now recall 
the classical classification of fibre bundles over $S^1\times\br^i$
with fibre $F$, the classification of manifold approximate fibrations 
over $S^1\times\br^i$ with fibre germ $F\times\br^{i+1}\to\br^{i+1}$, 
and the relationship
between these two classifications from Hughes-Taylor-Williams \cite{29},
\cite{30}. 
\def\bun{\text{\rm Bun}(S^1}
\def\mf{\maf(S^1}
Let $\bun\times\br^i)_F$ 
denote the simplicial set of fibre bundles over $S^1\times\br^i$ with
fibre $F$, so that there exists a homotopy equivalence
$\bun\times\br^i)_F\simeq\map(S^1,\btop(F))$. Since
$\pi_0\map(S^1,\btop(F)) = \pi_0\tp(F)$, there is a 
classifying isomorphism
$c_1:\pi_0\bun\times\br^i)_F\to\pi_0\tp(F)$.
Let $\mf\times\br^i)_{F\times\br^{i+1}}$ 
denote the simplicial set of manifold approximate fibrations
over $S^1\times\br^i$ with fibre germ the projection
$F\times\br^{i+1}\to\br^{i+1}$ and assume $\dim F + i\geq 4$. 
Since $S^1\times\br^i$ is parallelizable
it follows from \cite{29} that there is a homotopy equivalence
$\mf\times\br^i)_{F\times\br^{i+1}}\simeq\map(S^1,\btop^{\text{\rm c}}(F\times
\br^{i+1}))$ 
where
$\tp^{\text{\rm c}}(F\times\br^{i+1})$ 
denotes the simplicial group of controlled
homeomorphisms on $F\times\br^{i+1}$. Since
$\tp^{\text{\rm c}}(F\times\br^{i+1}))\simeq\topbf\times\br^{i+1})$ 
by Hughes-Taylor-Williams \cite{31}
and 
$\pi_0\map(S^1,\btop^b(F\times\br^{i+1})) = \pi_0\topbf\times\br^{i+1})$, 
there is a 
classifying isomorphism
$c_2:\pi_0\mf\times\br^i)_{F\times\br^{i+1}}\to\pi_0\topbf\times\br^{i+1})$.

\proclaim{Proposition 8.8}
If\/ $\dim F+i\geq 4$, then the following diagram commutes\rom:
$$\CD \pi_0\bun\times\br^i)_F @>{\varphi}>>
\pi_0\mf\times\br^i)_{F\times\br^{i+1}}\\
@V{c_1}V{\simeq}V  @V{\simeq}V{c_2}V\\
\pi_0{\tp}(F) @>{\sigma}>>
\pi_0\topbf\times\br^{i+1})
\endCD$$
where $\varphi$ is the forgetful map and $\sigma$ is euclidean stabilization
$[h]\mapsto [h\times\id_{\br^{i+1}}]$.
\endproclaim

\demo{Proof}
This follows from Hughes-Taylor-Williams \cite{30\rm, Thm\. 0.3}.
\qed
\enddemo

If $p:M\to S^1$ is a manifold approximate fibration with fibre germ
$F\times\br\to\br$, then the {\it monodromy} of $p$ is the class
$c_2(p) = [h]\in\pi_0\topbf\times\br)$ with $h:F\times\br\to F\times\br$
a bounded homeomorphism. The monodromy induces a well-defined homotopy
equivalence $F = F\times\{ 0\} @>{h|}>>
F\times\br\to F$ which in turn induces a homomorphism
$h_\ast:\Wh_1(\bz\pi_1 F)\to\Wh_1(\bz\pi_1 F)$, also called
the monodromy of $p$.

\proclaim{Theorem 8.9}
Let $p:M\to S^1$ be a manifold approximate fibration with fibre germ
$F\times\br\to\br$ and monodromy $[h]$ with $n=\dim F\geq 4$.

\noindent
{\rm (i)} The following are equivalent\rom:
\roster
\item $p$ is controlled homeomorphic to a fibre bundle projection with
fibre $F$.
\item $\tau\beta(c_2([p])) = \tau\beta([h]) =
0\in\Wh_1(\bz\pi_1 F)$.
\endroster
{\rm (ii)} The following are equivalent\rom:
\roster
\item $p\times\id_\br$ is controlled homeomorphic to a fibre
bundle projection with fibre $F$.
\item $\tau\beta(c_2([p])) =  \tau\beta([h]) \in\im N\subseteq
\Wh_1(\bz\pi_1 F).$
\endroster
{\rm (iii)} There exist a subgroup $G$ of 
$\widetilde{K}_0(\bz\pi_1 F)$ and a function 
$$N_0:G\to\Wh_1(\bz\pi_1 F)/\im N$$ 
such that the following are equivalent\rom:
\roster
\item $p\times\id_{\br^2}$ is controlled homeomorphic to a fibre
bundle projection with fibre $F$.
\item The class of $\tau\beta(c_2([p])) =  \tau\beta([h])$
in $\Wh_1(\bz\pi_1 F)/\im N$ is in $N_0(G)$.
\endroster
\endproclaim

\demo{Proof}(i) follows from Propositions 8.4 and 8.8.

\noindent
(ii) 
Consider the diagram
$$\CD {} @. \pi_0\tp(F)   @. {} \\
         @.  @VV{\sigma_1}V   @. \\
      \pi_0\cbf\times\br) @>\rho_1>>
              \pi_0\topbf\times\br) @>\beta>>
                \pi_0\ihcobf\\
         @.  @VV{\sigma_2}V   @.  \\
       {} @. \pi_0\topbf\times\br^2) @. {} 
\endCD$$
where $\sigma_1,\sigma_2$ denote euclidean stabilization and 
$\rho_1,\beta$ have been defined above.
According to Proposition 8.8, $p\times\id_\br$ is controlled homeomorphic
to a fibre bundle with fibre $F$ if and only if
$\sigma_2c_2([p])=\sigma_2([h])\in\im(\sigma_2\sigma_1)$. By Propositions 8.2, 8.4
and Lemma 8.7, $\sigma_2([h])\in\im(\sigma_2\sigma_1)$ if and only
if $\beta ([h])\in\im(\beta\rho_1)$ if and only if
$\tau\beta ([h])\in\im(N)$. Thus, (1) and (2) are equivalent.

\noindent
(iii) The diagram above can be extended to a diagram
$$\CD {} @. \pi_0\tp(F)   @. {} \\
         @.  @VV{\sigma_1}V   @. \\
      \pi_0\cbf\times\br) @>\rho_1>>
              \pi_0\topbf\times\br) @>\beta>>
                \pi_0\ihcobf @>{\tau}>>    \Wh_1(\bz\pi_1 F)\\
         @.  @VV{\sigma_2}V   @.  \\
    \pi_0\cbf\times\br^2) @>{\rho_2}>> \pi_0\topbf\times\br^2) @. {} \\
         @.  @VV{\sigma_3}V   @.  \\
     {} @. \pi_0\topbf\times\br^3) @. {} \\
\endCD$$
As above, $p\times\id_{\br^2}$ is controlled homeomorphic to a fibre 
bundle projection with fibre $F$ if and only if $\sigma_3\sigma_2([h])\in
\im(\sigma_3\sigma_2\sigma_1)$.
Since $\pi_0\cbf\times\br^2) \cong \widetilde{K}_0(\bz\pi_1 F)$ by Proposition
8.5, the result will follow from Lemma 8.7(ii) once it is obsevered that
the action of 
$\pi_0\topbf\times\br)$ on $\Wh_1(\bz\pi_1 F)$ satisfies items (3) and
(5) of 8.7(ii). The first follows from the fact that if 
$\sigma_2([h]) =\sigma_2([h'])$, then the induced homotopy equivalences
$h_1, h_1':F\to F$ are homotopic and, hence, ${h_1}_\sharp ={h_1'}_\sharp$.
The second follows from the explicit construction of $\beta$.
\qed
\enddemo

\remark{Remark \rom{8.10}}  It follows from Hughes-Taylor-Williams \cite{32}
that condition 8.9(i)(1) holds if and only if
$p$ is homotopic to a fibre bundle projection with fibre $F$.
It seems reasonable to conjecture that condition 8.9(ii)(1) holds 
if and only if
$p\times\id_\br$ is properly homotopic to to a fibre
bundle projection with fibre $F$.
\endremark

We will now prepare for a version of Theorem 8.9(i),(ii) 
where we allow the fibre
of the fibre bundle projections to vary
(Theorem 8.13 below). The following result says that
we do not have to worry about non-manifold fibres.

\proclaim{Lemma 8.11}
\widestnumber\item{(ii)}
\roster
\item"(i)" If $p:M\to S^1$ is a manifold approximate fibration with
$m=\dim M \geq 6$ and $p$ is controlled homeomorphic to a  bundle
projection, then $p$ is controlled homeomorphic to a bundle projection
with manifold fibre.
\item"(ii)" If $p:M\to S^1\times\br$ is a manifold approximate fibration with
$m=\dim M \geq 7$ and $p$ is controlled homeomorphic to a  bundle
projection, then $p$ is controlled homeomorphic to a bundle projection
with manifold fibre.
\endroster
\endproclaim

\demo{Proof}(i)
We may assume that $p:M\to S^1$ is a bundle projection. The fibre is a
compact ANR $X$. According to \cite{32} it suffices to show that $p$ is 
homotopic to a bundle projection with manifold fibre; that is,
we need to show that the Farrell fibering obstruction of $p$ vanishes.
We will use the version of the total fibering obstruction as
exposited in Ranicki \cite{49}. Let $h:X\to X$ be the classical
monodromy of $p$ so that
the mapping torus $T(h)$ is $M$. The infinite cyclic cover of $M$ is
$X\times\br$ with generating covering translation 
$\zeta:X\times\br\to X\times\br~;~ (x,t)\mapsto (h(x), t+1).$
The mapping torus $T(\zeta)$ has a preferred finite structure and the 
fibering obstruction is the torsion of the natural homotopy equivalence
$T(\zeta)\to T(h) = M$. The preferred finite structure on $T(\zeta)$
can be defined by choosing a finite CW complex $K$ and a homotopy
equivalence $f:K\to X$ (this exists by West \cite{64}). Let $g:X\to K$
be a homotopy inverse for $f$. Then $f,g$ induce a natural homotopy
equivalence $d:K\to X\times\br~;~ x\mapsto (f(x),0)$ 
and inverse $u:X\times\br\to K~;~ (x,t)\mapsto g(x)$. In
particular, this is a finite domination of $X\times\br$ so that
$T(u\zeta d)\to T(\zeta)$ is the preferred finite structure.
Note that $T(u\zeta d) = T(ghf)$ and the composition
$T(ghf)\to T(\zeta)\to T(h) = M$ is simple.

\noindent
(ii) We may assume that $p:M\to S^1\times\br$ is a bundle projection with
fibre a compact ANR $X$. Let $W= p^{-1}(S^1\times\{ 0\})$ and
$q=p|:W\to S^1\times\{0\} = S^1$.
Note that $X\times\br^2$ is a manifold since $M$ is a manifold and $p$
is a bundle projection; however, it is unknown whether this implies 
that $X\times\br$ is a manifold (cf\. Daverman \cite{9\rm, Prob\. 625}). In particular, $W$
might not be  a manifold. On the other hand, $W\times\br$ is homeomorphic
to $M$ so that $W$ is resolvable by Quinn \cite{47, 3.2.2}; that is, there
exist a manifold $N$, $\dim N =m-1\geq 6$, and a cell--like map
$r:N\to W$. It follows as in part (i) that $qr:N\to S^1$ is homotopic
to a fibre bundle projection with manifold fibre and hence (by \cite{32})
is controlled homeomorphic to a fibre bundle projection $q':N\to S^1$ with
manifold fibre.
Since $p:M\to S^1\times\br$ is fibre preserving homeomorphic to
$q\times\id_\br:W\times\br\to S^1\times\br$, $p$ is controlled homeomorphic
to $q\times\id_\br$. Siebenmann \cite{56} implies that $r\times\id_\br
:N\times\br\to W\times\br$  can be arbitrarily closely approximated
by homeomorphisms, so that $q\times\id_\br$ is controlled homeomorphic
to $qr\times\id_\br:N\times\br\to W\times\br$. 
Finally, $q\times\id_\br$ is controlled homeomorphic to 
$q'\times\id_\br$ which is a bundle projection with manifold fibre.
\qed
\enddemo

\proclaim{Lemma 8.12}
Let $p:M\to S^1$ be a manifold approximate fibration with fibre
germ $F\times\br\to\br$, monodromy $[h:F\times\br\to F\times\br]$,
and $n = \dim F\geq 4$. Suppose $F'$ is a closed  manifold for which
there is a bounded homeomorphism $k:F\times\br\to F'\times\br$.
\widestnumber\item{(ii)}
\roster
\item"(i)" $p$ is a manifold approximate fibration with fibre germ
$F'\times\br\to\br$ and monodromy $[khk^{-1}:F'\times\br\to F'\times\br]$.
\item"(ii)" If $\beta:\pi_0\topbf\times\br)\to\pi_0\ihcobf$
and $\beta':\pi_0\topbf'\times\br)\to\pi_0\ihcobfpr$
are the  `region between' functions defined above, then
there exists $x\in\Wh_1(\bz\pi_1F)$ such that
$$\tau\beta([h]) = (k^{-1})_\ast\tau\beta'([khk^{-1}]) + x - h_\ast(x)$$
where $(k^{-1})_\ast:\Wh_1(\bz\pi_1F')\to\Wh_1(\bz\pi_1F)$ is 
induced by the composition
$F'=F'\times\{ 0\}{\buildrel k^{-1}|\over\longrightarrow}F\times\br\to F$.
Moreover, $x$ is represented by the torsion of the $h$--cobordism
associated to the bounded homeomorphism $k^{-1}:F'\times\br\to F\times\br$.
\endroster
\endproclaim

\demo{Proof} (i) If $p$ is considered to have fibre germ $F\times\br\to
\br$, then the affect of the classifying map $c_2$ is to turn
$p:M\to S^1$ into a fibre bundle over $S^1$ with fibre $F\times\br$ and
structure group $\topbf\times\br)$. The monodromy $h$ is then the
classical monodromy of this bundle. The bundle can be considered to 
be a bundle with fibre $F'\times\br$, structure group
$\topbf'\times\br)$ and monodromy $khk^{-1}$.
See \cite{29}, \cite{30}.

\noindent
(ii) Choose $L> 0$ large. Let 
$$W=(khk^{-1})(F'\times(-\infty,L])\setminus F'\times (-\infty,0)\subseteq
F'\times\br$$ 
so that 
$$(W;F'\times\{0\},khk^{-1}(F'\times\{ L\})$$ 
is an $h$--cobordism whose torsion is $\tau\beta'([khk^{-1}]).$
Let 
$$W_k = k(F\times [-L,\infty))
\setminus F'\times (0,\infty)\subseteq F'\times\br$$
so that $(W_k;k(F\times\{-L\}),F'\times\{0\})$ is an
$h$--cobordism.
Let 
$$W_{k^{-1}} = F\times (-\infty, 2L]\setminus k^{-1}(F'\times
(-\infty,L))$$
so that $(W_{k^{-1}};k^{-1}(F'\times\{ L\}), F\times\{2L\})$ is an
$h$--cobordism. 
Let 
$$U=k^{-1}W_k\cup k^{-1}W\cup hW_{k^{-1}}\subseteq F\times\br.$$
Note that
$k^{-1}W_k\cap k^{-1}W = k^{-1}(F'\times\{ 0\})$,
$k^{-1}W\cap hW_{k^{-1}} = hk^{-1}(F'\times\{ L\})$,
$k^{-1}W_k\cap hW_{k^{-1}} = \emptyset$, 
and that 
$U =h(F\times (-\infty,2L])\setminus F\times (-\infty,-L)\subseteq F\times\br$
so that 
$(U;F\times\{-L\},h(F\times\{ 2L\}))$ is an $h$--cobordism
with torsion 
$\tau (U,F\times\{ -L\}) = \tau\beta([h])\in\Wh_1(\bz\pi_1F)$.
The standard sum and composition formulae imply that
$$\multline
\tau(U,F\times\{-L\}) \\
= \tau(k^{-1}W_k,F\times\{-L\})
+ (k^{-1})_\ast\tau(W,F'\times\{0\})
+ h_\ast(k^{-1})_\ast\tau(kW_{k^{-1}},F'\times\{ L\}).
\endmultline$$
Let $x = \tau(k^{-1}W_k, F\times\{-L\})\in\Wh_1(\bz\pi_1F)$.
It is easy to see that
$$x+(k^{-1})_\ast\tau(kW_{k^{-1}},F'\times\{ L\}) = 0$$ 
so that
$\tau(U,F\times\{-L\}) = x +(k^{-1})_\ast\tau(W,F'\times\{ 0\})-h_\ast(x).$
\qed
\enddemo

\proclaim{Theorem 8.13}
Let $p:M\to S^1$ be a manifold approximate fibration with fibre germ
$F\times\br\to\br$ and monodromy $[h]$.

\noindent
{\rm (i)} If  $n=\dim F\geq 5$, then the following are equivalent\rom:
\roster
\item $p$ is controlled homeomorphic to a fibre bundle projection.
\item $\tau\beta(c_2([p])) =  \tau\beta([h]) =
0\in\Wh_1(\bz\pi_1 F)/\im(1-h_\ast)$.
\endroster
{\rm (ii)} If  $n=\dim F\geq 6$, then the following are equivalent\rom:
\roster
\item $p\times\id_\br$ is controlled homeomorphic to a fibre
bundle projection.
\item $\tau\beta(c_2([p])) =  \tau\beta([h]) =
0\in\Wh_1(\bz\pi_1 F)/
(\im N + \im(1-h_\ast)).$
\endroster
\endproclaim

\demo{Proof}(i) (1) {\it implies\/} (2):
By Lemma 8.11(i) we may  assume that $p$ is controlled homeomorphic to a 
bundle projection with fibre a closed manifold $F'$.
By uniqueness of fibre germs \cite{29} there exists a bounded homeomorphism
$k:F\times\br\to F'\times\br$. An application of Theorem 8.9(i) with 
$F'$ replacing $F$ implies that
$\tau\beta'(c_2[p]) = 0\in\Wh_1(\bz\pi_1F')$.
Now Lemma 8.12(ii) implies that
$\tau\beta(c_2[p]) =\tau\beta ([h]) = x-h_\ast(x)$ 
for some $x\in\Wh_1(\bz\pi_1F)$.

\noindent
(2) {\it implies\/} (1):
If $\tau\beta([h]) = x-h_\ast(x)$ for some  $x\in\Wh_1(\bz\pi_1F)$,
choose an $h$--cobordism $(W;F,F')$ such that
$x=\tau(W,F)$. In fact, there is a bounded homeomorphism
$k:F\times\br\to F'\times\br$ such that
$W=k^{-1}(F'\times (-\infty,L])\setminus F\times (-\infty,0)$ 
for some large $L>0$
(this is the $h$--cobordism associated to $k^{-1}$).
Lemma 8.12(i) implies that $p$ is a manifold approximate fibration 
with fibre germ $F'\times\br\to\br$ and monodromy $khk^{-1}$.
It follows from Lemma 8.12(ii) that
$(k^{-1})_\ast\tau\beta'([khk^{-1}]) = 0$.
Hence $\tau\beta'([khk^{-1}]) = 0 \in\Wh_1(\bz\pi_1F')$.
Finally, Theorem 8.9(i) implies that $p$ is controlled homeomorphic
to a bundle projection with fibre $F'$.

\noindent
(ii) (1) {\it implies\/} (2):
By Lemma 8.11(ii) we may  assume that $p\times\id_\br$ 
is controlled homeomorphic to a 
bundle projection with fibre a closed manifold $F'$.
As in (i) there exists a bounded homeomorphism
$k:F\times\br\to F'\times\br$. 
By Lemma 8.12(ii)
there exists $x\in\Wh_1(\bz\pi_1F)$ such that
$\tau\beta([h]) = k^{-1}_\ast\tau\beta'([khk^{-1}]) + x - h_\ast(x)$.
Lemma 8.12(i) implies that $p$ is a manifold approximate fibration
with fibre germ $F'\times\br\to\br$ and monodromy $[khk^{-1}]$.
Since $p\times\id_\br$ is controlled homeomorphic to a fibre bundle
projection with fibre $F'$,
Theorem 8.9 implies that
$\tau\beta'([khk^{-1}]) = N'z$ for some
$z\in\Wh_1(\bz\pi_1F')$ where
$N':\Wh_1(\bz\pi_1F')\to\Wh_1(\bz\pi_1F')$ is the norm map.
Thus $\tau\beta([h]) = k^{-1}_\ast N'z + x - h_\ast(x) = Nk^{-1}_\ast z+
x-h_\ast x$.

\noindent
(2) {\it implies} (1):
Suppose $\tau\beta([h]) = Nz+x-h_\ast x$. As in (i) there exist a closed
manifold $F'$ and a bounded homeomorphism
$k:F\times\br\to F'\times\br$ such that
$x$ is represented by the torsion associated to $k^{-1}$ via the `region
between' construction. Lemma 8.12(ii) implies that
$\tau\beta([h]) = k_\ast^{-1}\tau\beta'([khk^{-1}])+x-h_\ast x$.
Hence, $k_\ast^{-1}\tau\beta'([khk^{-1}]) = Nz$ and
$\tau\beta'([khk^{-1}]) = k_\ast Nz = N'k_\ast^{-1}z$.
Since $p$ is a manifold approximate fibration with fibre germ
$F'\times\br\to\br$ and monodromy $khk^{-1}$, Theorem 8.9
implies that $p\times\id_\br$ is controlled homeomorphic to a fibre 
bundle with fibre $F'$.
\qed
\enddemo

\remark{Remark \rom{8.14}}
\widestnumber\item{(ii)}
\roster
\item"(i)" As in Remark 8.10  it follows 
from Hughes-Taylor-Williams \cite{32}
that condition 8.13(i)(1) holds if and only if
$p$ is homotopic to a fibre bundle projection.
It seems reasonable to conjecture that condition 8.13(ii)(1) holds 
if and only if
$p\times\id_\br$ is properly homotopic to to a fibre
bundle projection.
\item"(ii)" Another way to prove 8.13(i) is to identify $\tau\beta(c_2([p]))$ with
the Farrell fibering obstruction of $p$.
\endroster
\endremark

Let $\bz_q$ denote the finite cyclic group of order $q$.

\proclaim{Proposition 8.15}
\widestnumber\item{(iii)}
\roster
\item"(i)" If\/ $\pi_1(F)= \bz_q$, $q>3$ is prime, and $\dim F = n\geq 6$ is even,
then $N:\Wh_1(\bz[\bz_q])\to\Wh_1(\bz[\bz_q])$ is not surjective,
but $\tau\beta:\pi_0\topbf\times\br)\to\Wh_1(\bz[\bz_q])$ is
surjective.
\item"(ii)" If\/ $n\geq 5$ is odd and $q>3$ is prime, then there exists a closed manifold
$F$ such that $\dim F = n$, $\pi_1(F) = \bz_q$ and
$N:\Wh_1(\bz[\bz_q])\to\Wh_1(\bz[\bz_q])$ is not surjective
(in fact, it is the $0$ homomorphism),
but $\tau\beta:\pi_0\topbf\times\br)\to\Wh_1(\bz[\bz_q])$ is
surjective.
\item"(iii)" If $q$ is prime, then 
$\widetilde{K}_0(\bz[\bz_q])$ is a finite group and
$K_{-i}(\bz[\bz_q]) = 0$ for all $i>0$.
If, in addition,  $5\leq q\leq 19$, then
$\widetilde{K}_0(\bz[\bz_q]) = 0$.
\endroster
\endproclaim

\demo{Proof}
Let $q>3$ be a prime number. It is known that 
$\Wh_1(\bz[\bz_q])$ is free abelian of finite non-zero rank and the
standard involution $\ov{~\cdot~}$ acts by the identity
(Bass \cite{4}, Bass-Milnor-Serre \cite{5}, Wall \cite{60}; 
see Oliver \cite{46} for
an exposition).
Let $F$ be a closed manifold with $\dim F = n\geq 5$ and $\pi_1(F) =
\bz_q$.
Then
$N:\Wh_1(\bz\pi_1F)\to\Wh_1(\bz\pi_1F)$ is multiplication by
$2$ if $n$ is even and multiplication by
$0$ if $n$ is odd,
and therefore not surjective.
According to Lawson \cite{38}, if $n$ is even, $\pi_0\ihcobf = \pi_0\hcobf$,
and if $n$ is odd, there exist manifolds $F$ as above such that
$\pi_0\ihcobf = \pi_0\hcobf$.
Since Proposition 8.4 implies that $\beta:\pi_0\topbf\times\br)\to
\pi_0\ihcobf = \pi_0\hcobf$ is surjective, it follows that
$\tau\beta:\pi_0\topbf\times\br)\to\Wh_1(\bz[\bz_q])$ is
surjective.
This proves (i) and (ii).

\noindent
For (iii) see Rosenberg \cite{50\rm, pp\. 23, 157}.
\qed
\enddemo

\proclaim{Theorem 8.16}
\widestnumber\item{(iii)}
\roster
\item"(i)" If\/  $q$ is a prime, $5\leq q\leq 19$, $n\geq 6$ is even, 
$F$ is any closed manifold with
$\pi_1(F)=\bz_q$ and $\dim F=n$, then there exists a manifold
approximate fibration $p:M\to S^1$ with fibre germ $F\times\br\to\br$
such that for all $i\geq 0$, 
$p\times\id_{\br^i}:M\times\br^i\to S^1\times\br^i$
is not controlled homeomorphic to a fibre bundle projection with fibre $F$.
\item"(ii)" If\/ $q>3$ is prime, $n\geq 5$ is odd, 
then there exists a closed manifold
$F$ with
$\pi_1(F)=\bz_q$ and $\dim F=n$ and a manifold
approximate fibration $p:M\to S^1$ with fibre germ $F\times\br\to\br$
such that for all $i\geq 0$, 
$p\times\id_{\br^i}:M\times\br^i\to S^1\times\br^i$
is not controlled homeomorphic to a fibre bundle projection with fibre $F$.
\item"(iii)" If\/  $n\geq 6$ is even, 
$F$ is any closed manifold with
$\pi_1(F)=\bz_5$ and $\dim F=n$, then there exists a manifold
approximate fibration $p:M\to S^1$ with fibre germ $F\times\br\to\br$
such that for all $i\geq 0$, 
$p\times\id_{\br^i}:M\times\br^i\to S^1\times\br^i$
is not controlled homeomorphic to a fibre bundle projection.
\item"(iv)" If\/ $n\geq 5$ is odd, 
then there exists a closed manifold
$F$ with
$\pi_1(F)=\bz_5$ and $\dim F=n$ and a manifold
approximate fibration $p:M\to S^1$ with fibre germ $F\times\br\to\br$
such that for all $i\geq 0$, 
$p\times\id_{\br^i}:M\times\br^i\to S^1\times\br^i$
is not controlled homeomorphic to a fibre bundle projection.
\endroster
\endproclaim

\demo{Proof} (i)
According to Proposition 8.8 we need 
a manifold approximate fibration $p:M\to S^1$ with fibre germ $F\times\br
\to\br$ and monodromy 
$[h]\in\pi_0\topbf\times\br)$
such that for all $i\geq 0$, $[h\times\id_{\br^i}]
\in\pi_0\topbf\times\br^{i+1})$ is not in
$\im(\sigma:\pi_0\tp(F)\to\pi_0\topbf\times\br^{i+1}))$.
According to Propositions 8.2, 8.5, and 8.15(iii),
$\sigma:\pi_0\topbf\times\br^{2+i})\to\pi_0\topbf\times\br^{3+i})$ is
injective for all $i\geq 0$.
Hence, it suffices to find 
a manifold approximate fibration $p$ with monodromy
$[h]\in\pi_0\topbf\times\br)$
such that  $[h\times\id_{\br}]
\in\pi_0\topbf\times\br^{2})$ is not in
$\im(\sigma:\pi_0\tp(F)\to\pi_0\topbf\times\br^{2}))$;
that is, such that $p\times\id_\br$ is not controlled homeomorphic
to a fibre bundle projection with fibre $F$.
According to Theorem 8.9(ii) this is equivalent to
$\tau\beta([h])\not= 0\in \Wh_1(\bz\pi_1F)/{\im}N$.
Such monodromies exist by Proposition 8.15(i).

\noindent
(ii) This is similar to (i) except now we know only that
$\sigma:\pi_0\topbf\times\br^{3+i})\to\pi_0\topbf\times\br^{4+i})$ is
injective for all $i\geq 0$.
Hence, it suffices to find 
a manifold approximate fibration $p$ with monodromy
$[h]\in\pi_0\topbf\times\br)$
such that  $[h\times\id_{\br^2}]
\in\pi_0\topbf\times\br^{3})$ is not in
$\im(\sigma:\pi_0\tp(F)\to\pi_0\topbf\times\br^{3}))$;
that is, such that $p\times\id_{\br^2}$ is not controlled homeomorphic
to a fibre bundle projection with fibre $F$.
According to Proposition 8.15(ii),(iii)
$\Wh_1(\bz\pi_1F)/{\im}N = \Wh_1(\bz\pi_1F)$ and is infinite
(cf\. proof of 8.15).
Hence, since 8.15(iii) implies that
$\widetilde{K}_0(\bz[\bz_q])$ 
is finite, the result follows from Proposition 8.9(iii).

\noindent
(iii) As in (i) it suffices to find 
a manifold approximate fibration $p:M\to S^1$ with fibre germ $F\times\br
\to\br$ and monodromy 
$[h]\in\pi_0\topbf\times\br)$ such that
$p\times\id_\br$ is not controlled homeomorphic
to a fibre bundle projection.
According to Theorem 8.13(ii) this is equivalent to
$\tau\beta([h])\not= 0\in \Wh_1(\bz\pi_1F)/
({\im}N+{\im}(1-h_\ast))$.
But  $\Wh_1(\bz\pi_1F)= \Wh_1(\bz[\bz_5])$ is isomorphic to
$\bz$ so that $h_\ast = \pm 1$ and $1-h_\ast=0,2$. As noted in the proof
of Proposition 8.15, $N=0$ so that
$\Wh_1(\bz\pi_1F)/({\im}N+{\im}(1-h_\ast)) \not= 0$ and the 
result follows from Proposition 8.13(ii).

\noindent
(iv) is similar to (iii).
\qed
\enddemo

\demo{Proof of Theorem \rom{8.1}}
Let $X$ be the open mapping cylinder of a manifold approximate fibration
$p:M\to S^1$ constructed in Theorem 8.16(iii) or (iv). If $S^1\times\br^i$ 
had a fibre bundle
mapping cylinder neighborhood in $X\times\br^i$, then 
according to Theorem 2.2, $p\times\id_{\br^{i+1}}$ would be controlled
homeomorphic to a fibre bundle projection, contradicting Theorem 8.16.
Since block bundles with fibre $F$ are classified by 
${\text{\rm B}}\widetilde{\tp}
(F)$, equivalence classes of block bundles over $S^1\times\br^i$ 
correspond to $\pi_0\widetilde{\tp}(F)$. Since
$\pi_0{\tp}(F)\to\pi_0\widetilde{\tp}(F)$ is surjective,
the result on block bundles follows from the fibre bundle case.
\qed
\enddemo

\remark{Remark \rom{8.17}}
(i) If $(X,B)$ is a manifold stratified pair, it is the case that 
for a large enough torus, the quotient
$X\times T/(B\times T = B)$ does have a block
structure. Moreover, the block structure on the `links' is not arbitrary:
it has some nice transfer invariance properties. In other words, for each
simplex $\Delta$ of $B$ one has a nice manifold which maps to
$\Delta\times T$, with control in the $T$ direction. (What we have shown here
is that one cannot block over simplices of $B\times T$.)
This structure is called a STIBB\footnote{An equivalent notion 
is used by Yan in \cite{68}: one has blocks over
$\Delta\times E$ where $E$ is an Euclidean space, and the data is bounded in
the $E$ direction.} in \cite{62} and is applied 
there to give a stable surgery exact sequence
for stratified spaces. Indeed, if one had block structures stably
then the $L$--cosheaves in the stable classification theorem
\cite{62\rm, \S6.2} would have to have the `$s$' 
decoration (as in the `PT category'
in \cite{62\rm, \S6.1}) rather than the $-\infty$ decoration that arises.
The differences between these decorations are accounted for by Tate cohomology
calculations rather similar to those done here.

It is not too difficult to combine Theorem 2.2 with the classification
theorem of \cite{29}, and the stabilization theorem of \cite{63} to give
a proof of the stable classification theorem for 
$S^{-\infty}(X {\rel} B)$.
Using \cite{29} the stable germ neighborhoods are computed
by maps 
$[B, \btop^b(F\times E)]$
which is the same as
$[B, {\text{\rm B}}\widetilde{\text{\rm Top}}^b(F\times E)]$ 
by \cite{63}, the last of
which is computed by bounded block surgery using
$L^{-\infty}(\ho)$.
Different structures with the same germ near the singular stratum can the be
compared using ordinary \rel $\infty$ surgery on the complement.
The result of this analysis is just a Poincar\'e duality away from the
result as expressed in \cite{62}.

\noindent
(ii) These examples are closely related to those constructed
by Anderson \cite{1}.

\noindent
(iii) Husch \cite{33} used nontrivial inertial $h$--cobordisms to
construct exotic manifold approximate fibrations over $S^1$.

\noindent
(iv) Using the the tables for relative class numbers in
Washington \cite{61, p. 412}, it is possible to construct a few more
even dimensional manifolds as in Theorem 8.16(i) for primes $q$ with
$3< q < 67$. We don't know of other calculations which give more
manifolds as in Theorem 8.16(iii) and (iv).
\endremark


\head 9. Extensions of isotopies and $h$--cobordisms \endhead
In this section we combine the geometry of teardrop neighborhoods
with manifold approximate fibration theory in order to prove
parametrized isotopy extension  and $h$--cobordism extension
theorems for manifold stratified pairs.

\subhead Extending isotopies\endsubhead
\demo{Proof of Corollary 2.4 \rom(Parametrized Isotopy Extension\rom)}
Let $(X,B)$ be a manifold stratified pair with $\dim X\geq 5$ and $B$
a closed manifold. Suppose $h:B\times\Delta^k\to B\times\Delta^k$ is
a $k$--parameter isotopy (in particular, $h|B\times\{ 0\}=\id_{B\times
\{ 0\}}$). We are required to find a $k$--parameter isotopy
$\tilde{h}:X\times\Delta^k\to X\times\Delta^k$ extending $h$ which is 
supported in a given neighborhood of $B$.
Since $B$ has a teardrop neighborhood in $X$ (Theorem 2.1) there 
exist  an open neighborhood $U$ of $B$ in $X$ (which we can take to be
contained in the given neighborhood of $B$) and a proper map
$f:U\to B\times (-\infty,+\infty]$ such that
$f|:B\to B\times\{+\infty\}$ is the identity and
$f|:U\setminus B\to B\times B\times\br$ is a manifold approximate 
fibration. We consider $\Delta^k$ embedded as a convex 
subspace of $\br^k$ with the origin the zeroth vertex (basepoint) of
$\Delta^k$.
Define a $k$--parameter isotopy
$g:B\times\br\times\Delta^k\to B\times\br\times\Delta^k$ by 
letting $g_t:B\times\br\to B\times\br, t\in\Delta^k$,
be given by 
$$g_t(x,s)=\cases (h_t(x),s), & \text{if $s\geq 0$}\\
                  (h_{(1+s)t}(x),s), & \text{if $-1\leq s\leq 0$}\\
                  (x,s), & \text{if $s\leq -1$.}
\endcases$$
Let $\cal U$ be an open cover of $B\times\br$ whose mesh goes to $0$
near $B\times\{+\infty\}$; i\.e\., if
$V\in{\cal U}$ and $V\cap (B\times[N,+\infty)\not=\emptyset$
then $\diam V<{1\over N}$ for $N=1,2,3,\dots$ (cf\. the definition of
$\Psi$ in \S5).
By the Approximate Isotopy Covering Theorem for manifold approximate
fibrations (see \cite{28\rm, 17.4} for information on how this follows from
\cite{24}) there exists a $k$--parameter isotopy
$\tilde{g}:(U\setminus B)\times\Delta^k\to (U\setminus B)\times\Delta^k$
such that for each $t\in\Delta^k$
\roster
\item $f\tilde{g}_t$ is $\cal U$--close to $g_tf|(U\setminus B)$, and
\item $\tilde{g}_t|f^{-1}(B\times (-\infty,-2]) =$ the inclusion.
\endroster
Finally, define $\tilde{h}_t:X\to X$, $t\in\Delta^k$,
by
$$\tilde{h}_t=\cases h_t, & \text{on $B$}\\
                     \tilde{g}_t, & \text{on $U\setminus B$}\\
                     \id_{X\setminus U} & \text{on $X\setminus U$.}\qed
\endcases$$
\enddemo

\subhead Stratified $h$-cobordisms\endsubhead
Throughout the rest of this section we let $(X,B)$ be a fixed manifold 
stratified pair with $B$ a closed manifold with $\dim B\geq 5$.
We now define  stratified $h$--cobordisms. The definition is a bit
more complicated than in \cite{48} because we have not allowed manifold
strata to have boundaries.

\definition{Definition 9.1} A {\it stratified $h$--cobordism on}
$(X,B)$ is denoted $(\widetilde{W};\partial_0\widetilde{W},\partial_1
\widetilde{W})$ and consists of a homotopically stratified pair
$(\widetilde{W},W)$ with finitely dominated local holinks such that 

\noindent
(i) $\widetilde{W}$ is a locally compact separable metric space,

\noindent
(ii) there is an $h$--cobordism $(W;\partial_0W,\partial_1W)$
with $\partial_0W=B$,

\noindent
(iii) there are disjoint closed subspaces
$\partial_0\widetilde{W},\partial_1\widetilde{W}\subseteq\widetilde{W}$
with $X=\partial_0\widetilde{W}$
satisfying:
\roster
\item"(a)" $\partial_i\widetilde{W}\cap W=\partial_iW$ for $i=0,1$,
\item"(b)" $\widetilde{W}\setminus W$ is a manifold with
boundary $(\partial_0\widetilde{W}\setminus\partial_0W)\cup
(\partial_1\widetilde{W}\setminus\partial_1W)$,
\item"(c)" $\partial_i\widetilde{W}$ is a stratum preserving
proper strong deformation retract of $\widetilde{W}$ for $i=0,1$.
\endroster
The stratified $h$--cobordism $(\widetilde{W};\partial_0\widetilde{W},
\partial_1\widetilde{W})$ is said to {\it extend\/}
the $h$--cobordism $(W;\partial_0W,\partial_1W)$.
Note that $(\widetilde{W}\setminus W;\partial_0\widetilde{W}\setminus
\partial_0W,\partial_1\widetilde{W}\setminus\partial_1W)$ is a
proper $h$--cobordism on $\partial_0\widetilde{W}\setminus\partial_0W$.
\enddefinition

The following result is not needed in the rest of this section, but is included
to show that stratified $h$-cobordisms keep one inside the category of
manifold stratified pairs.

\proclaim{Proposition 9.2} 
If $(\widetilde W;\partial_0\widetilde W,\partial_1\widetilde W)$
is a stratified $h$-cobordism on $(X,B)$ extending the $h$-cobordism
$(W;\partial_0W,\partial_1W)$ on $B$, then
$(\partial_1\widetilde W,\partial_1W)$ is a manifold stratified pair.
\endproclaim

\demo{Proof} By definition $(\widetilde W,W)$ is a homotopically stratified
pair with finitely dominated local holinks. Of course,
$\partial_1W$ and $\partial_1\widetilde W\setminus\partial_1W$ are manifolds.
The forward tameness of $\partial_1W$ in $\partial_1\widetilde W$ follows from
the facts that $W$ is forward tame in $\widetilde W$ and $\partial_1\widetilde W$ 
is a stratum preserving retract of $\widetilde W$.
Moreover, since $q:\ho(\widetilde W,W)\to W$ is a fibration with finitely
dominated fibre and a stratum preserving strong deformation of $\widetilde W$
to $\partial_1\widetilde W$ induces a strong deformation retraction of
$\ho(\widetilde W,W)$ to $\ho(\partial\widetilde W,\partial_1W)$
which, when restricted to $q^{-1}(\partial_1W)$ is fibre preserving
over $\partial_1W$, it follows that $\ho(\partial_1\widetilde W,\partial_1W)\to
\partial_1W$ is a fibration with finitely dominated fibre.
\qed\enddemo

We now fix some notation which will be used throughout the
rest of this section.

\definition{Notation 9.3}
Since $B$ has  a teardrop neighborhood in $X$
(Theorem 2.1) there exist an open  neighborhood $U$ of $B$ in $X$ 
and a proper map
$f:U\to B\times (-\infty,+\infty]$ such that
$f|:B\to B\times\{+\infty\}$ is the identity and
$f|:U\setminus B\to B\times\br$ is a manifold approximate 
fibration. 
\enddefinition

\definition{Definition 9.4} An {\it $h$-cobordism on $X$ rel $B$\/}
consists of: 

\noindent
(i) a proper $h$-cobordism $(V;\partial_0V,\partial_1V)$ on
$\partial_0V=X\setminus B$ (in particular, $\partial_iV$ is a proper
strong deformation retract of $V$ for $i=0,1$),

\noindent
(ii) a map of triads
$$g:(N;\partial_0N,\partial_1N)\to (B\times\br\times [0,1];
B\times\br\times\{ 0\}, B\times\br\times\{ 1\})$$
where:
\roster
\item"(a)" $N$ is an open subset of $V$ and is a neighborhood of the
end of $V$ determined by $B$ (i.e., for a proper retraction $r:V\to
X\setminus B$ there exists a neighborhood $U'$ of $B$ in $X$ such that
$r^{-1}(U'\setminus B)\subseteq N$),
\item"(b)" $\partial_iN=N\cap\partial_iV$ for $i=0,1$,
\item"(c)" $g$ is a proper approximate fibration,
\item"(d)" $\partial_0N=U$,
\item"(e)" $g|\partial_0N=f$.
\endroster
\enddefinition 

Here is some explanation for this definition.

\remark{Remarks 9.5}
\roster
\item The teardrop $V\cup_g(B\times [0,1])$ contains 
$X=\partial_0V\cup_{g|}B\times\{ 0\}$ 
so that the triad
$$(V\cup_gB\times [0,1];X,\partial_1\cup_{g|}B)$$ 
is a
stratified $h$-cobordism on $(X,B)$ extending the trivial $h$-cobordism
on $B$. The fact that the properties of Definition 9.1 are indeed satisfied
is a special case of Theorem 9.6 below.
This is why $(V;\partial_0V,\partial_1V)$ is called an $h$-cobordism on
$X$ rel $B$: because $V$ can be compactifed (if $X$ is compact) by adding 
$B\times [0,1]$ to obtain a stratified $h$-cobordism on $(X,B)$ which
is trivial on $B$.
\item Suppose $(\widetilde W;\partial_0\widetilde W,\partial_1\widetilde W)$
is any stratified $h$-cobordism on $(X,B)$ extending 
$(W;\partial_0W,\partial_1W)$. It follows that
$(\widetilde W\setminus W;\partial_0\widetilde W\setminus\partial_0W,
\partial_1\widetilde W\setminus\partial_1W)$
is an $h$-cobordism on $X$ rel $B$.
As noted above, this is obviously a proper $h$-cobordims on $X\setminus B$.
A proof of the other properties in Definition 9.4 requires the
advanced teardrop technology  from \cite{26},\cite{27}
(because $\widetilde W$ has more than two strata). Likewise, using this 
advanced teardrop technology we will be able to reformulate Definition 9.4
to be more along the lines of Definition 9.1. It is because \cite{27} has not
yet appeared that we are taking the current approach.
\item A simple example of an $h$-cobordism on $X$ rel $B$ is the trivial
one $((X\setminus B)\times [0,1]; X\setminus B\times\{ 0\},
X\setminus B\times\{ 1\})$. 
For the open set $N\subseteq (X\setminus B)\times [0,1]$ in Definition 9.4(ii)
we take $(U\setminus B)\times [0,1]$. Thus, the Teardrop Neighborhood
Existence Theorem 2.1 is required to show that the trivial $h$-cobordism
is an example.
Theorem 9.6 below, 
when applied to this trivial
$h$-cobordism, is nevertheless non-trivial. This special case (stated as
Corollary 9.7) best 
illustrates the power of the techniques of the current paper without making
motivational appeal to advanced teardrop technology.
\endroster
\endremark

The next result shows how teardrop technology can be used to
extend an $h$-cobordism on $B$ to a teardrop neighborhood of $B$ in $X$.
Moreover, the extension can be chosen so that on the complement of $B$,
it is any given $h$-cobordism on $X$ rel $B$.
The key fact that makes teardrop technology applicable to
this problem is that $h$-cobordisms on $B$ become
trivial $h$-cobordisms on $B\times\br$ after crossing with $\br$.

\proclaim{Theorem 9.6} 
Let $(X,B)$ be a manifold stratified pair
with $B$ a closed manifold, $\dim B\geq 5$.
If $(V;\partial_0V,\partial_1V)$ is an $h$-cobordism on $X$ rel $B$
and $(W;\partial_0W,\partial_1W)$ is an $h$-cobordism on $B$, then
there exists a stratified $h$-cobordism 
$(\widetilde{W};\partial_0\widetilde{W},\partial_1\widetilde{W})$ 
extending $(W;\partial_0W,\partial_1W)$
such that
$$(\widetilde{W}\setminus W;\partial_0\widetilde{W}\setminus
\partial_0W,\partial_1\widetilde{W}\setminus\partial_1W)
=(V;\partial_0V,\partial_1V).$$
\endproclaim

\demo{Proof}
As is well-known $(W;\partial_0W,\partial_1W)\times\br$ is a trivial
$h$--cobordism; i\.e\., there exists a homeomorphism
$h:W\times\br\to B\times\br\times[0,1]$ 
such that
$h|:\partial_0W\times\br=B\times\br\to B\times\br\times\{0\}$
is the identity.
Let $N\subseteq V$ and $g:N\to B\times\br\times [0,1]$
be as in Definition 9.4. 
Define
$\tilde{f}:N\to W\times\br$ to be the composition
$$\tilde{f}:N @>{g}>> B\times\br\times [0,1]
@>{h^{-1}}>>  W\times\br.$$
Form the  teardrop
$\widetilde{W} = V\cup_{\tilde{f}} W$. 
The pair $(\widetilde{W},W)$ is homotopically stratified 
with finitely dominated local holinks and
$\widetilde W$ is a locally compact separable metric space
by Corollary 4.10.
Let $\partial_i\widetilde{W} = \partial_iV\cup_{g|}
B\times\{ i\}$ for $i=0,1$ which clearly are disjoint closed
subsets of $\widetilde{W}$, and $\partial_0\widetilde{W} =X$.
Note that $\widetilde{W}\setminus W = V$
is a manifold with boundary $\partial_0V\cup\partial_1V$ as required.
In order to show that $\partial_i\widetilde{W}$ is a stratum preserving
strong deformation retract of $\widetilde{W}$ for $i=0,1$, 
one can use the fact that $\partial_iV$ is a strong deformation retract
of $V$ together with the homotopy extension theorem, to show that
it suffices to define stratum preserving strong deformation retractions
on $N\cup_{\tilde f}W$.
We concentrate
on the $i=0$ case since the $i=1$ case is similar.
Since $\partial_0W\hookrightarrow W$ is a homotopy equivalence,
there exists a strong deformation retraction
$r:W\times I\to W$
of $W$ to $\partial_0W$ (thus, $r_0=\id_W$, $r_1(W)\subseteq\partial_0W$
and $r_t|\partial_0W$ equals the inclusion for $t\in I$).
Since $\tilde{f}:N\to W\times\br$ is an
approximate fibration, there exists a homotopy
$\tilde{r}:N\times I\to N$ such that
\roster
\item $\tilde{r}_0 =\id_{N}$,
\item $\tilde{r}_t|\partial_0N =$ inclusion for each
$t\in I$,
\item $\tilde{r}_1(N)\subseteq \partial_0N$,
\item if $(x,s)\in\tilde{f}^{-1}(W\times [k,+\infty))\subseteq 
N$ and $k=1,2,3,\dots$, then for each $t\in I$
$$d(\tilde{f}\tilde{r}(x,s,t),r(\tilde{f}(x,s),t)) < 1/k.$$
\endroster
(This comes from approximately lifting the homotopy $r$ with very good
control near $W\times\{+\infty\}$. To get condition (3), first
get a homotopy as above that pulls $N$ close to
$\partial_0N$, in fact, 
so close that an additional push along a collar
will not destroy the estimates in condition (4).)
Define $R:N\cup_{\tilde f}W\times I\to\widetilde{W}$ 
by requiring
$R|W\times I=r$ and $R|N\times I =\tilde{r}$.
The continuity of $R$ follows from Lemma 3.4.
\qed
\enddemo

\proclaim{Corollary 9.7 ($h$--cobordism Extension)}
If $(W; \partial_0W,\partial_1W)$ is an $h$--cobordism
with $\partial_0W=B$, then there exists a stratified $h$--cobordism
$(\widetilde{W};\partial_0\widetilde{W},\partial_1\widetilde{W})$ with
$\partial_0\widetilde{W} = B$ extending $W$.
\endproclaim

\demo{Proof} This follows immediately from Theorem 9.6.
\qed
\enddemo

\remark{Remark 9.8}
(i) Quinn \cite{48\rm, 1.8} 
gives an $h$--cobordism theorem
for stratified spaces. He shows that if a suitable torsion vanishes
the $h$--cobordism is a product, but does not prove there is a
realization theorem for torsions (cf\. \cite{48\rm, p\. 498}). The 
realization for ${\Wh}^{\text{\rm top}}(X {\rel} B)$ 
(the set of equivalence classes of $h$-cobordisms on $X$ rel $B$)
is a natural
extension of the realization of elements of Siebenmann's proper
Whitehead group ${\Wh}^p(W)$ for a noncompact manifold $W$ with a 
tame end \cite{54}. Indeed the latter is the special case of the former
obtained by one point compactifying $W$ (see the picture on p\. 132 of
\cite{62}). What is missing from \cite{48} then is the proof that 
${\Wh}^{\text{\rm top}}(X)\to{\Wh}^{\text{\rm top}}(X {\rel} B)\times
{\Wh}(B)$ is surjective
(where ${\Wh}^{\text{\rm top}}(X)$ is the set of equivalence classes of 
stratified $h$-cobordisms on $X$).
Theorem 9.6 completes the missing step. 
F. Connolly and B. Vajiac have recently obtained related results.

\noindent
(ii) 
We suspect that there is a fibration of $h$-cobordism spaces
whose fibration sequence at $\pi_0$ contains this discusion.
We hope to return to this, as well as a discussion of stratified
$h$-cobordisms on manifold stratified spaces with more than
two strata, in a later paper.

\noindent
(iii) Jones \cite{35} proved a concordance extension theorem for locally flat 
submanifolds of topological manifolds of dimension greater than four.
His proof uses manifold approximate fibration techniques which also work
for a manifold stratified pair $(X,B)$ with $\dim X\geq 5$ such that
$B$ has a mapping cylinder neighborhood in $X$. It seems likely that his
techniques extend to arbitrary (high dimensional) manifold stratified
pairs.
At any rate, his work
is further evidence for a 
moduli space interpretation of the results of this section.
\endremark


%
%
\enddocument